\documentclass{memo-l}

\usepackage{graphicx}
\usepackage[dvipsnames]{xcolor}

\usepackage[colorlinks]{hyperref}
\newcommand\myshade{100}
\hypersetup{
  linkcolor  = RoyalBlue!\myshade!black,
  citecolor  = ForestGreen!\myshade!black,
  urlcolor   = RedOrange!\myshade!black,
  colorlinks = true,
}
\usepackage{bookmark}

\usepackage{subcaption}

\usepackage{todonotes}

\usepackage{bm}
\usepackage{bbold}
\usepackage{tikz-cd}

\usepackage{thmtools}
\usepackage[capitalise]{cleveref}
\crefformat{equation}{(#2#1#3)}
\crefname{chapter}{Chapter}{Chapters}
\crefname{section}{Section}{Sections}
\crefname{subsection}{Section}{Sections}

\newtheorem{theorem}{Theorem}[chapter]
\newtheorem{thm}[theorem]{Theorem}
\newtheorem{lem}[theorem]{Lemma}
\newtheorem{prop}[theorem]{Proposition}
\newtheorem{cor}[theorem]{Corollary}

\theoremstyle{definition}
\newtheorem{defn}[theorem]{Definition}
\newtheorem{example}[theorem]{Example}
\newtheorem{question}[theorem]{Question}

\theoremstyle{remark}
\newtheorem{rem}[theorem]{Remark}

\numberwithin{section}{chapter}
\numberwithin{equation}{chapter}

\newcommand\N{\ensuremath{\mathbb{N}}}
\newcommand\R{\ensuremath{\mathbb{R}}}

\newcommand{\el}[1]{\ensuremath{\ell_{#1}}}

\newcommand{\Lr}[1]{\ensuremath{\mathcal{L}_r(#1)}}

\newcommand{\one}{\ensuremath{\mathbb{1}}}
\newcommand{\falg}{$f\!$-algebra}
\newcommand{\fstar}{$f^{*}\!$-algebra}

\newcommand{\lalg}{$\ell$-algebra}
\newcommand{\VLA}{\operatorname{VLA}}
\newcommand{\FAFA}[1]{\operatorname{FA{\it f}\!A}(#1)}
\newcommand{\FAFAv}[1]{\operatorname{FA{\it f}\!A}[#1]}
\newcommand{\FNFA}[1]{\operatorname{FN{\it f}\!A}[#1]}
\newcommand{\FBFA}[1]{\operatorname{FB{\it f}\!A}[#1]}
\newcommand{\FBFAs}[1]{\operatorname{FB{\it f}\!A}(#1)}
\newcommand{\orth}[1]{\operatorname{Orth}(#1)}
\newcommand{\FVL}[1]{\operatorname{FVL}(#1)}
\newcommand{\FVLv}[1]{\operatorname{FVL}[#1]}
\newcommand{\FBL}[1]{\operatorname{FBL}[#1]}

\newcommand{\Sinfn}{S_\infty^n}

\newcommand{\FAFAov}[1]{\operatorname{FAFA}^1[#1]}

\newcommand{\trinor}[1]{{\vert\kern-0.25ex\vert\kern-0.25ex\vert #1 
    \vert\kern-0.25ex\vert\kern-0.25ex\vert}}
\newcommand{\bigtrinor}[1]{{\bigg\vert\kern-0.25ex\bigg\vert\kern-0.25ex\bigg\vert #1 
    \bigg\vert\kern-0.25ex\bigg\vert\kern-0.25ex\bigg\vert}}

\DeclareMathOperator{\supp}{supp}
\DeclareMathOperator{\Lat}{Lat}

\DeclareMathOperator{\barespn}{span}

\DeclareMathOperator{\range}{range}
\DeclareMathOperator{\Fin}{Fin}

\begin{document}

\frontmatter

\title{Free Banach \falg s}


\author[D. Muñoz-Lahoz]{David Muñoz-Lahoz}
\address{Instituto de Ciencias Matem\'aticas (CSIC-UAM-UC3M-UCM)\\
Universidad Autónoma de Madrid\\
C/ Nicol\'as Cabrera, 13--15, Campus de Cantoblanco UAM\\
28049 Madrid, Spain.}
\curraddr{}
\email{david.munnozl@uam.es}
\thanks{Research supported by an FPI-UAM 2023 contract funded by
Universidad Autónoma de Madrid.}

\author[P. Tradacete]{Pedro Tradacete}
\address{Instituto de Ciencias Matem\'aticas (CSIC-UAM-UC3M-UCM)\\
Consejo Superior de Investigaciones Cient\'ificas\\
C/ Nicol\'as Cabrera, 13--15, Campus de Cantoblanco UAM\\
28049 Madrid, Spain.}
\curraddr{}
\email{pedro.tradacete@icmat.es}
\thanks{Research partially supported by grants PID2020-116398GB-I00,
PID2024-162214NB-I00 and
CEX2023-001347-S funded by MCIN/AEI/10.13039/501100011033.}

\date{\today}

\subjclass[2020]{06F25, 46B42, 46A40, 46M15}

\keywords{Banach lattice, \falg, free object}


\begin{abstract}
    We construct and analyze the free Banach \falg\ $\FBFA E$
    generated by a Banach space $E$, extending recent developments
    on free Banach lattices to the setting of Banach \falg s, where
    multiplication interacts with the lattice structure. Starting
    from the explicit realization of the free Archimedean \falg\ as
    a sublattice-algebra of $\mathbb{R}^{E^*}\!\!$, we develop a new
    structure theorem for normed \falg s that allows us to identify
    the kernel of the maximal submultiplicative lattice seminorm as
    precisely those functions vanishing on the unit ball $B_{E^*}$.
    This yields a representation of the free normed \falg\ inside
    $C(B_{E^*})$. We prove that this representation extends to an
    injective map on the completion $\FBFA E$ if and only if $\FBFA
    E$ is semiprime, and we establish that $\FBFA E$ is indeed
    semiprime whenever $E$ is finite-dimensional or $E=L_1(\mu)$.
    This is closely related to approximating operators into a Banach
    \falg\ by operators into finite-dimensional Banach \falg s. We
    also use the newly constructed free objects to provide an
    example of a semiprime normed \falg\ whose norm completion is
    not semiprime. Using the tools developed for the study of
    free objects, we show the following extension
    property: if $A$ is a closed sublattice-algebra of a Banach \falg\
    $B$, then every real-valued lattice-algebra homomorphism on $A$
    extends to a real-valued lattice-algebra homomorphism on $B$.

    For a finite-dimensional Banach space $E$, a complete
    description of $\FBFA E$ is obtained: $\FBFA E$ is
    lattice-algebra isomorphic to $C([0,1]\times S_{E^{*}})$ equipped
    with pointwise order and the product given by
    \[
        (f\star g)(r,u)=rf(r,u)g(r,u).
    \]
    As a consequence, Banach spaces of the same dimension generate
    isomorphic free Banach \falg s. The interplay between the
    lattice and algebraic structures also leads to unexpected
    behavior: for instance, the free normed \falg\ is always order
    dense in $\FBFA E$, whereas for free Banach lattices this only
    holds when $E$ is finite-dimensional.
\end{abstract}

\maketitle

\tableofcontents

\chapter*{Introduction}
This paper is devoted to the construction and study of free objects
in the category of Banach \falg s. Our main motivation comes from
the recent developments in the framework of free Banach lattices,
and our aim is to explore the limits of this construction in the
more general setting of Banach \falg s. Since every Banach lattice
equipped with the identically zero multiplication becomes a Banach
\falg, this is a natural generalization. Our goal is to provide the
first steps in the foundation of free Banach \falg s and open up new
research directions on the interaction between Banach spaces and
Banach \falg s.

The reader might naturally wonder why we should focus on (Banach)
\falg s and on free objects. After addressing these questions and
explaining the motivation behind our research, we outline the
contents of the paper, its main results, and the technical and
conceptual difficulties they involve.

\section*{Why (Banach) \falg s?}

The notions of \falg\ and $f\!$-ring were introduced by
Birkhoff and Pierce back in \cite{birkhoff_pierce1956}. Since then, these structures have
been extensively studied from an algebraic point of view. A
comprehensive account of this line of research can be found in the
survey \cite{henriksen1997} and its more than 150 references.

But why were these notions introduced in the first place? To answer this question, let us consider the framework of vector lattices consisting of real-valued functions defined on a set. Most such spaces
naturally carry a multiplication operation (typically, the pointwise
product) under which the positive cone remains closed. Vector lattice algebras provide an abstraction of this situation. However, Birkhoff and Pierce soon realized that the general
theory of vector lattice algebras was not particularly tractable, as a further condition was needed to properly abstract the behavior of function spaces. This missing ingredient was the \falg\ condition: if $x\wedge y=0$ and $z\ge 0$,
then
\[
    (zx)\wedge y=0=(xz)\wedge y.
\]
This property encapsulates the idea that pointwise product of functions preserves its
supports. The mathematical fact that led
Birkhoff and Pierce to this particular condition is the following:
a vector lattice algebra is the subdirect product
of totally ordered algebras if and only if it satisfies the \falg\
condition. As a consequence, one of the most celebrated results
in the early theory of \falg s was obtained: every Archimedean \falg\ is
associative and commutative (even if associativity is not assumed as
part of the definition).

The first application of \falg s to Operator Theory
is due to Bigard and Keimel \cite{bigard_keimel1969} and
to Conrad and Diem \cite{conrad_diem1971}, independently: the set of
orthomorphisms (order bounded band preserving operators) on an arbitrary vector lattice, equipped with the usual order and with composition as product, forms an \falg. This observation
gave rise to an intensive development of the theory of \falg s from an
operator-theoretic perspective. The works of Huijsmans
\cite{huijsmans1985,huijsmans1989}, de
Pagter
\cite{depagter1981,huijsmans_depagter1982,huijsmans_depagter1984_subalgebras,depagter1984},
Bernau
\cite{bernau_huijsmans1990_unit,bernau_huijsmans1990_algebras,bernau_huijsmans1995},
van Rooij \cite{buskes_depagter_vanrooij1991}, and Wickstead
\cite{wickstead1980,wickstead1987},
among others, greatly contributed to this progress. More recent
developments and applications include:
the study of averaging and Reynolds operators
\cite{boulabiar_buskes_triki2007}, the construction of vector
lattice powers \cite{boulabiar_buskes2006}, $\mathbb{L}$-functional
analysis \cite{kikianty_etal2024}, $\mathbb{L}$-valued integration
\cite{jiang_vanderwalt_wortel2024}, tensor products
\cite{buskes_wickstead2017,jaber2020}, free objects
\cite{keimel1995,de_jeu2021}, and the parametrization of all \falg\
products on AL and AM-spaces \cite{munoz-lahoz2025_products}.

Banach \falg s were first investigated in the 1980s by Martignon
\cite{martignon1980} and Scheffold
\cite{scheffold1981,scheffold1991,scheffold1993}. It is remarkable how
little attention these objects have received, given that Banach
lattices have been extensively investigated and that Banach \falg s are
particularly well behaved. First of all, one should note that every Banach lattice equipped with the
identically zero multiplication becomes a Banach \falg\ in a standard (yet somehow trivial) way. Also recall that every Banach \falg\ with identity is lattice-algebra isomorphic to a space
of continuous functions---a fact that directly connects with the original motivation
of modelling spaces of functions. 

A relevant property in the research of \falg s is concerned with the existence of divisors of zero: an \falg\ is semiprime if and only if $x^2=0$ implies $x=0$. In our research, we will soon face the following elementary question: is the
completion of a normed semiprime \falg\ also semiprime? This innocent-looking question turned out to be surprisingly subtle. Although it is not difficult to verify that the completion of a normed \falg\ is a Banach
\falg, the answer to the question above will be negative: in \cref{ex:comp_falg} we will construct a normed
semiprime \falg\ whose completion is not semiprime.

\section*{Why free objects?}

The previous paragraph justifies how free Banach \falg s can be a useful tool to tackle certain questions. From a more general perspective, it is natural to assume that in order to understand a given category, one must understand its free objects.
Free groups, algebras, modules, lattices, and Boolean algebras, among others, play a
central role throughout algebra and order theory. In topology, the Stone--Čech
compactification and the completion of a metric or uniform space can
also be viewed as instances of free constructions. A modern and increasingly
powerful approach in Functional Analysis is to study free objects in
categories related to Banach spaces. Lipschitz free Banach spaces \cite{godefroy_kalton2003},
holomorphic free Banach spaces \cite{mujica1991}, and free Banach
lattices \cite{aviles_rodriguez_tradacete2018} are notable examples
(see the survey \cite{garcia-sanchez_dehevia_tradacete2024} for a
general point of view on these). It is therefore natural to attempt
extending the theory of Banach \falg s in this direction.

The main inspiration for developing free Banach \falg s comes from the
theory of free Banach lattices. The free Banach lattice generated by a
set was first introduced in \cite{de_pagter_wickstead2015}. Although this foundational work
established several key properties, many questions remained open, as
the authors lacked an explicit description of the norm in such objects.
This gap was later filled in \cite{aviles_rodriguez_tradacete2018}, where both an explicit construction and
a norm were provided for the more general case of the free Banach lattice
generated by a Banach space. This breakthrough resolved several of the
open problems from \cite{de_pagter_wickstead2015} and also answered a question posed by J.\ Diestel
concerning weakly compactly generated Banach lattices. Following these
two seminal papers, the theory of free Banach lattices has flourished,
with significant developments in various directions (for instance, free complex
Banach lattices \cite{dehevia_tradacete2023},
free Banach lattices generated by a lattice
\cite{aviles_rodriguez-abellan2019,aviles_etal2022} and free Banach lattices
satisfying convexity conditions \cite{jardon-sanchez_etal2022}) and applications in the study of
norm-attaining lattice homomorphisms \cite{bilokopytov_etal2025,dantas_etal2022}, relations between linear and
lattice embeddings \cite{aviles_etal2022_embeddings}, constructions of push-outs \cite{aviles_tradacete2023} and complemented subspaces of Banach lattices
\cite{dehevia_tradacete2025}.
The study of free Banach lattices themselves has also become
a subject of great interest
\cite{oikhberg_etal2022,aviles_tradacete_villanueva2019,garcia-sanchez_tradacete2019,aviles_plebanek_rodriguez-abellan2018,dantas_etal2021}, as they serve as the canonical tool for
understanding the interaction between Banach spaces and Banach lattices.

The study of free Banach \falg s represents a further step in this
direction. Although these objects (and more generally, free Banach lattice
algebras) have been considered before
\cite{de_jeu2021,wickstead2017_questions}, no successful construction is yet
known of a free object where the norm, lattice, and algebraic
structures interact. The main obstacle lies in the presence
of a product, which renders most known techniques ineffective. Indeed,
the existing methods for free Banach lattices often rely on the fact that equations in a vector lattice are
positively homogeneous---something that no longer holds in \falg s.
New approaches are therefore required. In what follows, we outline the
nature of these difficulties, the techniques we have employed, and the
results we have obtained.

\section*{Outline and main results}

The core of the paper lies in Chapters~\ref{sec:FBFA}
and~\ref{sec:representation}. In \cref{sec:FBFA} we describe the
abstract construction of the free Banach \falg\ generated by a Banach
space $E$, and try to make this construction as explicit as possible. More
precisely, in \cref{sec:archimedean} we start from the explicit description of the free
Archimedean \falg\ $\FAFAv E$ as a sublattice-algebra of $\R^{E^{*}}$
and work our
way to the free normed and Banach \falg s using a known
general argument. This argument, however, does not provide anything
close to an explicit description of the objects, because one of the steps in
the construction
involves quotienting out $\FAFAv E$ by the kernel of a maximal seminorm
that is not given explicitly.

Using a novel structure theorem for Banach \falg s developed in
\cref{sec:structure}, we are able to compute this kernel:
it consists of the functions that vanish on the
unit ball $B_{E^{*}}$ (\cref{thm:kernel}). This elegant result is
particularly useful, as this
set is precisely the kernel of the restriction map $\FAFAv E\to
C(B_{E^{*}})$. Quotienting it out we obtain an injective lattice-algebraic
representation of the free normed \falg\ $\FNFA E$ inside
$C(B_{E^{*}})$ (\cref{cor:iota}). \Cref{sec:structure} also contains
two further developments of independent interest. For every semiprime
Banach \falg\ $A$, we construct a contractive representation $A\to
C(K_A)$, where $K_A$ denotes the compact Hausdorff space of
lattice-algebra homomorphisms on $A$, and characterize it by a natural
extension property (\cref{prop:char_KA2}). Building on this, we
establish a remarkable extension result: every lattice-algebra
homomorphism $\phi \colon B\to \R$ defined on a closed
sublattice-algebra $B$ of a Banach \falg\ extends to a lattice-algebra
homomorphism on the whole space (\cref{thm:extension}). Such an
extension does not, in general, preserve the norm; in fact, we show
that there are no non-zero $\lambda $-injective Banach \falg s for any
$\lambda \ge 1$ (\cref{thm:no_injectives}).

The reader familiar with free Banach lattices may not be surprised. After all, it is a known fact that the free Banach lattice generated by the Banach space $E$ can be represented
inside $C(B_{E^{*}})$. Yet several new difficulties arise in our setting. The elements of the free vector
lattice are positively homogeneous, so the restriction map $\FVLv E\to
C(B_{E^{*}})$ is already injective. In contrast, the elements of $\FAFAv E$ are not necessarily
homogeneous, and therefore $\FAFAv E$ cannot be represented inside
$C(B_{E^{*}})$ via this map. It is a remarkable coincidence that, in
constructing $\FNFA E$, we quotient out exactly the kernel of this
non-injective restriction map.

The representation of $\FNFA E$ inside $C(B_{E^{*}})$ has immediate applications. For instance, the free normed
\falg\ generated by a Banach space is always semiprime
(\cref{prop:FNFAsemiprime}). A much more substantial application is the
description, up to isomorphism, of the free Banach \falg\ generated by
a finite-dimensional Banach space $E$. This space is isomorphic to
$C([0,1]\times S_{E^{*}})$ equipped with pointwise order and multiplication given by
\[
    (f\star g)(r,u)=rf(r,u)g(r,u),
\]
for $(r,u)\in [0,1]\times S_{E^{*}}$ and $f,g \in C([0,1]\times
S_{E^{*}})$, together with the map $\eta \colon E\to C([0,1]\times
S_{E^{*}})$ defined by 
\[
\eta (x)(r,u)=u(x)
\]
for $x \in E$ (\cref{thm:finite}). This characterization has important consequences:
the free Banach \falg\ generated by a finite-dimensional Banach space
is semiprime (\cref{cor:finite_semiprime}), and Banach spaces of the
same dimension generate isomorphic free Banach \falg s
(\cref{cor:dim_isomorphic}). It is worth noting here that unlike in the case of free Banach
lattices, the latter result is not immediate from the definition,
as free Banach \falg s allow only the extension of contractive operators. Indeed, a
contractive isomorphism between finite-dimensional Banach spaces need
not extend to an isomorphism between the free Banach \falg s they
generate (\cref{ex:extension_of_bijection}).

The proof of the above theorem relies heavily on the study of formal
expressions involving the lattice, linear and multiplication operations. These expressions, called LLA expressions, are analyzed in
\cref{sec:archimedean}. The main result, \cref{thm:falgYudin}, states
that every LLA expression that vanishes on $\R$ must also vanish on
every Archimedean \falg. This result, already known in the context of
universal algebras, is here reproved in a language more familiar to
functional analysts. As part of the proof, we establish analytic
properties of LLA expressions (\cref{lem:lla_limit}) that will later
prove useful in several other results. Using these facts, we obtain an
explicit description of the free Archimedean \falg\ generated by a set
(\cref{thm:FAFA}) and study its properties
(Chapters~\ref{sec:lat_prop_FAFA} and~\ref{sec:alg_prop_FAFA}).

Returning to free Banach \falg s, note that the above characterization
for finite-dimensional Banach spaces provides only limited information
about the free norm. One of the most striking features of the free
Banach lattice is that its norm can be computed by considering only
operators with range in $\el 1^{n}$. In \cref{thm:norm_fin_dim} we
prove an analogous result: the norm in $\FBFA{\el 1^{n}}$ can be
computed by restricting to operators into finite-dimensional and semiprime Banach
\falg s. After some additional results, we extend this to show that the
same holds for $\FBFA{L_1(\mu)}$ for any measure $\mu$
(\cref{cor:tau_is_norm_example}).

The representation of the free normed \falg\ $\FNFA E$ inside $C(B_{E^{*}})$ plays a
central role in understanding this space. It is then natural to ask:
is the extension of this map to its completion (that is, to the free
Banach \falg) also injective? \Cref{sec:representation} is devoted to
this question. When the answer is affirmative, we say that $\FBFA E$
is representable in $C(B_{E^{*}})$. This property is significant,
since many results about free Banach lattices rely only on the fact
that they map injectively into $C(B_{E^{*}})$, and one may expect similar
phenomena for the free Banach \falg. However, the question has no
immediate answer, as the representation of $\FNFA E$ in $C(B_{E^{*}})$
is not an isomorphic embedding.

\Cref{thm:semiprime_iff_representable} clarifies this issue by showing
that $\FBFA E$ is representable in $C(B_{E^{*}})$ if and only if it is
semiprime. The next natural question is then whether $\FBFA E$ itself
is semiprime. We know that $\FNFA E$ is, but the completion of a
semiprime normed \falg\ need not be semiprime
(\cref{ex:comp_falg}). In general, it remains open whether $\FBFA E$
is semiprime, though the results of \cref{sec:FBFA} show that it is
semiprime whenever $E$ is finite-dimensional or $E=L_1(\mu)$ for an
arbitrary measure $\mu$.

\Cref{sec:properties} examines the lattice and algebraic properties of
the free Banach \falg\ generated by a Banach space. As we mentioned above, much
more can be said when these objects are representable in
$C(B_{E^{*}})$ (that is, when they are semiprime). It is remarkable
that, although the free vector lattice generated by a Banach space $E$ is order dense in
the corresponding free Banach lattice only when $E$ is finite-dimensional, the free normed \falg\ is always order dense in the
corresponding free Banach \falg\ whenever the latter is semiprime
(\cref{cor:ord_dense_FBFA}). We also analyze the behavior of an
operator $T\colon E\to F$ with respect to its natural extension
$\bar{T}\colon \FBFA E\to \FBFA F$
(\cref{prop:extension_properties}). These properties are subtler in this context than
in the free Banach lattice setting. As mentioned earlier, a
contractive isomorphism $T\colon E\to F$ need not extend to an
isomorphism $\bar{T}\colon \FBFA E\to \FBFA F$, even when $E$ and $F$
are finite-dimensional (\cref{ex:extension_of_bijection}).

For isometries, however, the situation is more accessible. Isometric
Banach spaces have lattice-algebra isometric free Banach \falg s, and
in \cref{sec:isometries} we prove the converse, provided that the
underlying Banach spaces have smooth duals (\cref{thm:isometries}).
Finally, \cref{sec:identity} is devoted to the construction of free
(Archimedean, normed, Banach) \falg s with an algebraic identity. The
discussion of these objects is considerably simpler, since every
Banach \falg\ with a norm-one identity is lattice-algebra isometric to
a space of continuous functions.

\section*{Acknowledgments}

The authors would like to thank Eugene Bilokopytov for valuable
discussions and ideas, and Jochen Glück for suggesting
\cref{ex:comp_falg_simple}.

\mainmatter
\chapter{Background}

In order to make the text as self-contained as possible, this chapter
contains the necessary background on \falg s, Banach \falg s, and free
Banach lattices to follow the exposition. For the basic facts
regarding vector and Banach lattices, and standard notation, we refer
the reader to
\cite{aliprantis_burkinshaw2006,lindenstrauss_tzafriri1979,abramovich_aliprantis2002}.

\section{A primer on \falg s}\label{sec:background_falg}

The theory of \falg s has been studied from several points of view.
Each of them has its own tools, strengths and weaknesses. In this
section we will motivate and introduce \falg s from three points of
view: the algebraic, the operator-theoretic, and the analytic.

\subsection{The algebraic approach}

The notions of $f\!$-ring and \falg\ were originally introduced by Birkhoff and
Pierce \cite{birkhoff_pierce1956} as part of their seminal work on
lattice-ordered rings and algebras. Here we will focus only on
lattice-ordered algebras, even though many of the results stated below
hold for general lattice-ordered rings.

\begin{defn}
    A \emph{partially ordered algebra} is a real algebra $A$ together with a
    partial order $\le $ for which, for all $x,y,z \in A$,
    \begin{enumerate}
        \item $x\le y$ implies $x+z \le y+z$,
        \item $x\le y$ implies $\lambda x\le \lambda y$ for all
            $\lambda \in \R_+$,
        \item $x,y\ge 0$ implies $xy \ge 0$.
    \end{enumerate}
    A \emph{lattice-ordered algebra} or \lalg\ is a partially ordered
    algebra which is a lattice under $\le $ (i.e., given any two
    elements, their supremum and infimum exist).
\end{defn}

Thus an \lalg\ is in particular a vector lattice. For this reason, and
because of the well-established term Banach lattice algebra (see
\cref{sec:background_norms}), we will usually refer to \lalg s as \emph{vector lattice
algebras}.

A subset of a vector lattice algebra that is closed under the linear
and lattice operations (i.e., a vector sublattice) and the product
(i.e., a subalgebra) is called an \emph{$\ell$-subalgebra} or
\emph{sublattice-algebra}. A map between vector lattice algebras that
is both a lattice and algebra homomorphism is called an
\emph{$\ell$-algebra homomorphism} or
\emph{lattice-algebra homomorphism}. As usual, the direct product of vector
lattice algebras is again a vector lattice algebra with the
operations and order defined coordinate-wise, and the associated coordinate
projections are lattice-algebra homomorphisms. Another important
construction, which comes from the theory of universal algebras, is
the following.

\begin{defn}
    Let $\{A_\gamma\}_\gamma $ be an indexed family of vector lattice
    algebras. Let $\prod_\gamma A_\gamma $ be its direct product,
    and let $p_\alpha \colon \prod_\gamma A_\gamma \to A_\alpha$ be the
    associated coordinate projections. A \emph{subdirect product} of
    $\{A_\gamma \}_\gamma $ is a sublattice-algebra $A$ of the direct
    product $\prod_\gamma A_\gamma $ such that, for every $\gamma $ and
    every $x_\gamma  \in A_\gamma $, there exists an $x \in A$ with
    $p_\gamma (x)=x_\gamma $.
\end{defn}

The following was originally proved in \cite{birkhoff_pierce1956} (see
also the more easily accessible \cite[Section VII.5]{birkhoff1967}).

\begin{thm}\label{thm:subdirect_product}
    A vector lattice algebra $A$ is a subdirect product of totally ordered
    algebras if and only if
    \[
    a\wedge b=0\text{ implies }(ac)\wedge b=0=(ca)\wedge b\text{ for
    all
    }a,b,c \in A_+.
    \]
\end{thm}

\begin{defn}
    A vector lattice algebra satisfying the conditions of the previous
    theorem is called an \emph{\falg}.
\end{defn}

The name \falg\ stands for ``function algebra;'' it will become clear
in
\cref{sec:background_analytic} why these
structures deserve to be considered as such.
\Cref{thm:subdirect_product} is very useful to prove equalities and
inequalities in \falg s, for it implies that they only need to be
checked in totally ordered algebras.

\begin{prop}
    Let $A$ be an \falg, and let $a,b,c \in A$ with $c\ge 0$. Then the
    following hold:
    \begin{enumerate}
        \item $c(a\vee b)=(ca)\vee (cb)$ and $(a\vee b)c=(ac)\vee
            (bc)$;
        \item $c(a\wedge b)=(ca)\wedge (cb)$ and $(a\wedge b)c=(ac)\wedge
            (bc)$;
        \item $|ab|=|a| |b|$;
        \item if $a\wedge b=0$, then $ab=0$ (in particular,
            $a^2=|a|^2\ge 0$).
    \end{enumerate}
\end{prop}

Similarly, one can check that $n|ab - ba|\le a^2+b^2$ for every $a$
and $b$ in an \falg\ and every $n \in \N$. This readily implies
another classical property of \falg s.

\begin{thm}
    Every Archimedean \falg\ is commutative.
\end{thm}

Another characterization of \falg s, first due to Fuchs \cite[Section
IX.3]{fuchs1963}, is the following.

\begin{thm}\label{thm:bands_ring_ideals}
    A vector lattice algebra is an \falg\ if and only if all its bands
    are ring ideals.
\end{thm}

It is straightforward to check that, if an \falg\ $A$ has an algebraic identity $e$, then $e$ must be positive and a weak order unit. This implies that $a \wedge ne \uparrow a$ for
every $a \in A_+$; in fact, the convergence is uniform, as B.\ de Pagter
proved in \cite[Lemma 10.6]{depagter1981}.

\subsection{The operator-theoretic approach}

Before presenting further properties of \falg s, it is necessary to
bring operators to the table. This is the approach to \falg s of
\cite[Section 2.3]{aliprantis_burkinshaw2006}, the reference we will
mainly follow.

\begin{defn}
    An operator $T\colon X\to X$ on a vector lattice $X$ is \emph{band
    preserving} if $T(B)\subseteq B$ for every band $B$ of $X$.
\end{defn}

\begin{prop}\label{prop:band_preserving}
    For an operator $T\colon X\to X$ on a vector lattice $X$, the
    following are equivalent.
    \begin{enumerate}
        \item $T$ is band preserving.
        \item $x \perp y$ implies $Tx \perp y$.
        \item For each $x \in X$, $Tx$ is contained in the principal
            band generated by $x$.
    \end{enumerate}
\end{prop}

Not every band preserving operator is order bounded. An example, due
to M.\ Meyer, of
a non-order bounded band preserving operator can be
found in \cite[Example 2.38]{aliprantis_burkinshaw2006}.

\begin{defn}
    An \emph{orthomorphism} is an order bounded and band preserving
    operator.
\end{defn}

By \cref{prop:band_preserving}, an orthomorphism $T\colon X\to X$
preserves disjointness (i.e., $x\perp y$ implies $Tx \perp Ty$) and
therefore every positive orthomorphism is also a lattice homomorphism.
Note also that every multiplication operator in an \falg\ is an
orthomorphism. Let $\orth X$ denote the set of orthomorphisms on a vector
lattice $X$. The following remarkable result is due independently to
Bigard and Keimel \cite{bigard_keimel1969}, and Conrad and Diem
\cite{conrad_diem1971}.

\begin{thm}
    Let $X$ be an Archimedean vector lattice. Then $\orth X$, under
    the pointwise linear operations and order, and the product given by composition, is an Archimedean \falg. Moreover, if $T \in \orth X$,
    then
    \[
    |T|(|x|)=|T(|x|)|=|Tx|\quad\text{for all }x \in X.
    \]
    In particular, the lattice operations on $\orth X$ can be computed pointwise on
    positive vectors.
\end{thm}

Next we are recalling that orthomorphisms enjoy some nice properties
that have immediate implications in the theory of \falg s.

\begin{thm}
    Every orthomorphism on an Archimedean vector lattice is order
    continuous.
\end{thm}

\begin{cor}
    Multiplication in an Archimedean \falg\ is order continuous.
\end{cor}

\begin{thm}
    Let $X$ be an Archimedean vector lattice and let $T \in \orth X$.
    Then
    \[
        \ker T = [T(X)]^{d}.
    \]
    In particular, the kernel of an orthomorphism is a band.
\end{thm}

Thus, if two orthomorphisms coincide on a set, then they also coincide
on the band generated by that set. This was used by Zaanen
\cite{zaanen1975} to show, without applying \cref{thm:subdirect_product}, that every Archimedean \falg\ is commutative. It can also be used to prove the following.

\begin{cor}
    Let $X$ be an Archimedean vector lattice and let $e>0$. Then there
    exists at most one product on $X$ that makes $X$ an \falg\ having
    $e$ as its multiplicative unit.
\end{cor}

Let $A$ be an Archimedean \falg. For every $a \in A$, the left
multiplication operator $L_a\colon A\to A$, defined by $L_a(b)=ab$ for all
$b \in A$, is an orthomorphism. We therefore have a natural map $L\colon A\to
\orth A$ given by $L(a)=L_a$. The kernel of this map is the
\emph{annihilator of $A$}:
\[
N(A)=\{\, a \in A : ab=0\text{ for all }b \in A \, \}.
\]
The annihilator of an Archimedean \falg\ has a remarkably simple
characterization.

\begin{thm}
    Let $A$ be an Archimedean \falg. Its annihilator $N(A)$ is a band,
    and it satisfies
    \begin{align*}
        N(A)&=\{\, a \in A : a^{k}=0\text{ for some }k \in \N \, \}\\
            &=\{\, a \in A : a^2=0 \, \}.
    \end{align*}
    Moreover, if $A$ has the principal projection property, then
    $N(A)$ is a projection band.
\end{thm}

When $N(A)=\{0\}$ (i.e., when $a^2=0$ implies $a=0$), we say that $A$
is \emph{semiprime}. Note that every \falg\ with identity is
semiprime. In this case the map $L\colon A\to \orth A$ is
injective, and moreover $L(A)$ is order dense in $\orth A$
\cite[Proposition 2.1]{huijsmans_depagter1984_bidual}. Thus we can think of $\orth A$
as a unitarization of $A$. It is obvious that, if $L$ is surjective,
then $A$ has an identity. The converse is also true.

\begin{thm}
    Let $A$ be an Archimedean \falg\ with identity. Then the map $L$
    is bijective, and therefore $A$ is isomorphic to $\orth A$.
\end{thm}

Thus the possible spaces of orthomorphisms of Archimedean
vector lattices are precisely the \falg s with identity. From this
theorem it also follows that $\orth{\orth X}$ can be
canonically identified with $\orth X$ for every Archimedean vector lattice $X$.

Recall that $X$ has an (essentially unique) order (or Dedekind)
completion $X^{\delta }$. This $X^{\delta }$ is characterized by the
fact that it is an order complete vector lattice in which $X$ embeds
as an order dense and majorizing vector sublattice. It also has an
(essentially unique) universal completion $X^{u}$. This $X^{u}$ is
characterized by the fact that it is a universally complete vector
lattice (i.e., it is
order complete and every set of pairwise disjoint positive vectors
has an upper bound) in which $X$ embeds as an order dense sublattice.
With the usual identifications, we have $X\subseteq X^{\delta
}\subseteq X^{u}$.

The vector lattice $X^{u}$ can be explicitly constructed as a vector
lattice of the form $C^{\infty }(K)$. Here $K$ is an extremally
disconnected compact Hausdorff space, and $C^{\infty }(K)$ is the set
of continuous functions $f\colon K\to [-\infty ,\infty ]$ that are
finite on a dense set. It can be shown that $C^{\infty }(K)$ is a
universally complete vector lattice under pointwise linear and lattice
operations. Moreover, under pointwise multiplication, it is an
Archimedean \falg\ with the constant one function as a multiplicative
identity. The details of this construction can be found in \cite[Section 50]{luxemburg_zaanen1971}.

This illustrates another way in which \falg s appear naturally when
studying vector lattices. And there is more. Since orthomorphisms are
order continuous, every positive
orthomorphism $T \in \orth{X}_+$ can be extended uniquely to an
orthomorphism, first on $X^{\delta }$, and then on
$X^{u}$. Conversely, every orthomorphism on $X^{u}$ that leaves $X$
invariant restricts to an orthomorphism on $X$. Thus, with these
identifications, we have natural inclusions as sublattice-algebras:
\[
    \orth X \subseteq \orth{X^{\delta }} \subseteq \orth{X^{u}}.
\]
Since $X^{u}$ is an \falg\ with identity, $\orth{X^{u}}$
is precisely the set of multiplication operators by elements of $X^{u}$.
Hence every $T \in \orth X$ is of the form $T(x)=yx$, for some $y
\in X^{u}$, where the multiplication is taken in $X^{u}$. This culminates
the intimate relation between orthomorphisms and
\falg s. The omitted details may be found in \cite[Theorem
2.63]{aliprantis_burkinshaw2006}.

As an application of these facts, we have the following.

\begin{thm}[{\cite[Theorem 2.64]{aliprantis_burkinshaw2006}}]
    \label{thm:repr_identity}
    Let $A$ be an Archimedean \falg\ with identity. Then $A$ is a
    sublattice-algebra of $A^{u}$.
\end{thm}

In particular, every Archimedean \falg\ with identity can be
represented as an \falg\ of extended functions on some extremally
disconnected compact Hausdorff space. This remarkable fact leads to
our final approach: \falg s as spaces of functions.

\subsection{The analytic approach}\label{sec:background_analytic}

Consider a family of (extended) real-valued functions on a set. When this family
is closed under pointwise addition and scalar multiplication, it is a
vector space. When it is also closed under pointwise supremum and
infimum, it is an Archimedean vector lattice. Conversely, every abstract
Archimedean vector
lattice can be realized as a space of extended real-valued functions through its uniform
completion, as explained above. Suppose that, additionally, this
family of functions is closed under pointwise product. Then it is not only a real
algebra, but also an \falg, as one may readily check. The idea is that
the condition defining \falg s encapsulates the fact that pointwise
product preserves supports (this is also closely related to \cref{thm:bands_ring_ideals}).

As it happens for
Archimedean vector lattices, abstract Archimedean \falg s can be
realized as concrete spaces of functions. For \falg s with identity,
this was shown in \cref{thm:repr_identity} (see also
\cite{henriksen_johnson1961} for an alternative construction). For
a semiprime Archimedean \falg\ $A$, we can embed $A$ inside $\orth A$,
which is an Archimedean \falg\ with identity, and then apply
\cref{thm:repr_identity}.

In the general case, when $A$ is not necessarily semiprime, we have to consider other products in spaces of functions beyond the pointwise product. The
reason is that the pointwise product is always semiprime. For a general
$A$, we have the following.

\begin{thm}[{\cite[Corollary 2.5]{buskes_wickstead2017}}]
    Let $A$ be an Archimedean \falg, and let $A^{u}=C^{\infty }(K)$ be
    its universal completion, where $K$ is an extremally disconnected
    compact Hausdorff space. Identify $A$ as a vector sublattice of
    $C^{\infty }(K)$, and let $\ast$ denote the product in $A$. Then
    there exists $w \in C^{\infty }(K)_+$ such that
    \[
        (a \ast b)(t)=w(t)a(t)b(t)
    \]
    for all $t \in K$ for which $w$, $a$ and $b$ are finite.
\end{thm}

\section{About $f$-algebras and norms}\label{sec:background_norms}

Once we see \falg s as spaces of functions, it is natural to add a
(complete) norm to this structure, just as it is done for vector lattices (to
yield Banach lattices) and for real algebras (to yield Banach
algebras). Even though Banach lattices and Banach algebras have been
extensively studied, the interaction between both structures (i.e.,
Banach lattice algebras) is not that well understood
\cite{wickstead2017_open,wickstead2017_questions}.

\begin{defn}
    A \emph{Banach lattice algebra} is a vector lattice algebra
    equipped with a complete lattice norm that is also
    submultiplicative. A \emph{Banach lattice algebra with identity}
    is a Banach lattice algebra together with a positive algebraic
    identity.
\end{defn}

In some contexts we ask for the algebraic identity to have norm one.
This is closely related to positiveness, as the following theorem shows.

\begin{thm}[{\cite[Theorem 2.2]{wickstead2017_questions}}]
    In a Banach lattice algebra $A$ with multiplicative identity $e$
    the following are equivalent:
    \begin{enumerate}
        \item $e\ge 0$,
        \item there is an equivalent Banach lattice algebra norm on
            $A$, $\trinor{\cdot }$, for which $\trinor{e}=1$.
    \end{enumerate}
\end{thm}

Examples of Banach lattice algebras include: $C(K)$, for some compact
Hausdorff space $K$, with pointwise operations; the space of regular
operators $\Lr X$ on an order complete Banach lattice $X$; and $L_1(G)$,
for a locally compact group $G$, where multiplication is given by
convolution. See \cite[Example 3.3]{munoz-lahoz2025} for further
details and examples.

The following facts about Banach lattice algebras are relevant
to our exposition. These results appear scattered through the literature
\cite{martignon1980,scheffold1984,huijsmans1995}; a full proof can be found in \cite[Theorem
3.4]{munoz-lahoz2025}. Given a Banach lattice algebra $A$ with
identity $e$, $A_e$ denotes the principal ideal generated by $e$
in $A$, and $\|{\cdot }\|_e$ denotes the associated norm.

\begin{thm}\label{thm:BLA}
    Let $A$ be a Banach lattice algebra with identity $e$.
    \begin{enumerate}
        \item The space $(A_e,\|{\cdot }\|_e)$, with the order and
            product of $A$, is a Banach lattice algebra lattice-algebra
            isometric to $C(K)$, for a certain compact Hausdorff space
            $K$. Moreover, if
            $\|e\|=1$, then $\|x\|_e=\|x\|$ for every $x \in A_e$.
        \item The ideal $A_e$ is a principal projection band in $A$.
        \item The space $A_e$ is an inverse closed subalgebra of $A$.
    \end{enumerate}
\end{thm}

Let us now focus on the relation between norms and \falg s.

\begin{defn}
    A \emph{normed \falg} is an \falg\ together with a lattice norm
    that is submultiplicative. A normed \falg\ is a \emph{Banach
    \falg} if the norm is complete.
\end{defn}

Banach \falg s are first mentioned as such in the work of Martignon
\cite{martignon1980}. It is surprising that Banach \falg s have
received relatively little attention in the literature, at least compared, on the one hand, to the theory of Banach lattices and algebras and, on the other
hand, to the theory of \falg s without norms. We hope that by now it
is clear to the reader why it is natural to further understand these
objects. The purpose of this work is to show that not only is this natural, but also interesting, since the interaction between the
norm and the \falg\ structure will prove to be rich and non-trivial. Let us first provide some examples of Banach \falg s.

\begin{example}\leavevmode
    \begin{enumerate}
        \item For every compact Hausdorff space $K$, the space $C(K)$,
            with pointwise order and product, is a Banach \falg. More
            generally, if $w \in C(K)_+$ with $\|w\|_\infty \le 1$,
            then the product
            \[
                (f \ast g)(t)=w(t)f(t)g(t),\quad t \in K,\;f,g \in
                C(K),
            \]
            turns $C(K)$ into a Banach \falg.
        \item Any Banach lattice with the zero product is a Banach
            \falg.
        \item Let $I$ be any set, and let $w \in \el \infty (I)_+$ be
            such that $\|w\|_\infty \le 1$. Then for $1\le p \le\infty$, $\el p(I)$ with
            product
            \[
                (x \ast y)(i)=w(i)x(i)y(i),\quad i \in I,\; x,y \in
                \el p(I)
            \]
            is a Banach \falg.
        \item Let $H$ be a Hilbert space, and let $D$ be a commuting
            set of bounded Hermitian operators on $H$. Its second
            commutant $D''$ is a Banach \falg\ with positive cone the
            one inherited from $B(H)$ and composition as product.
    \end{enumerate}
\end{example}

The following is a direct consequence of \cref{thm:BLA}.

\begin{cor}
    Every Banach \falg\ with identity is lattice-algebra isomorphic
    to $C(K)$, for some compact Hausdorff space $K$. If, moreover, the
    identity has norm $1$, then the isomorphism is an isometry.
\end{cor}

In view of this result, one may be led to think that Banach \falg s
are rather trivial. However, in the case they do not have an identity, several elementary questions arise that are surprisingly far from trivial. For instance:
given a semiprime normed \falg, is its completion (which is a Banach
\falg) also semiprime? This seemingly
innocent question will play a central role in our exposition (see
\cref{sec:representation}). We will be able to provide
in \cref{ex:comp_falg} a semiprime normed \falg\ whose completion is not
semiprime.

Biduals are also an important tool to deal with \falg s. Recall
that the bidual $X^{* *}$ of a Banach lattice $X$ is again a Banach
lattice, with the supremum and infimum of $x^{* *},y^{* *} \in X^{*
*}$ given by the Riesz--Kantorovich formulae:
\[
    (x^{* *}\vee y^{* *})(z^{*})=\sup \{\, x^{* *}(z^{*}_1)+y^{*
    *}(z^{*}_2) : z^{*}_1,z^{*}_2\ge 0,\, z_1^{*}+z_2^{*}=z^{*} \, \},
\]
\[
    (x^{* *}\wedge y^{* *})(z^{*})=\inf \{\, x^{* *}(z^{*}_1)+y^{*
    *}(z^{*}_2) : z^{*}_1,z^{*}_2\ge 0,\, z_1^{*}+z_2^{*}=z^{*} \, \},
\]
for all $z^{*}\in X^{*}_+$. This Banach lattice structure makes the
canonical embedding $j\colon X\to X^{* *}$ a lattice homomorphism.
Furthermore, if $T\colon X\to Y$ is a lattice homomorphism between
Banach lattices $X$ and $Y$, its second adjoint $T^{* *}\colon X^{* *}\to
Y^{* *}$ is again a lattice homomorphism.

If $A$ is a Banach \falg, then we can also endow its bidual $A^{* *}$
with either of the two Arens products \cite{huijsmans1989}. Since $A$ is commutative, both
Arens products coincide. Recall that the first Arens product is
constructed iteratively as follows. First, for every $x \in A$ and
$x^{*}\in A^{*}$, the product $x^{*}x$ is defined by
\[
    (x^{*}x)(y)=x^{*}(xy),\quad y \in A.
\]
For every $x^{*} \in A^{*}$ and $x^{* *}\in A^{* *}$, the product
$x^{* *}x^{*} \in A^{*}$ is defined by
\[
    (x^{* *}x^{*})(x)=x^{* *}(x^{*} x),\quad x \in A.
\]
Finally, for every $x^{* *},y^{* *} \in A^{* *}$, the Arens product
is defined by
\[
    (x^{* *}y^{* *})(x^{*})=x^{* *}(y^{* *}x^{*}),\quad x^{*} \in A^{*}.
\]
This product makes the canonical embedding $j\colon A\to A^{* *}$ an
algebra homomorphism. It was shown by Scheffold \cite{scheffold1991}
that $A^{* *}$, when equipped with the lattice structure indicated
above and the Arens product, becomes again a Banach \falg.
Furthermore, if $T\colon X\to Y$ is a lattice-algebra homomorphism
between Banach \falg s $X$ and $Y$, its double adjoint $T^{* *}\colon
X^{* *}\to Y^{* *}$ is again a lattice-algebra homomorphism.

\section{Free Banach lattices}\label{sec:background_FBL}

From \cref{ex:comp_falg}, and many other results in this paper, it is clear
that free objects play an important role in the correct understanding
of Banach \falg s. More generally, to properly understand any category, one must
understand its free objects.

The study of free objects in Functional Analysis has sparked a great
interest in recent years. The free Banach lattice generated by a
Banach space \cite{aviles_rodriguez_tradacete2018} is the main
inspiration for our work. For this reason, we will recall its
construction and main properties. A detailed exposition of the topic
can be found in \cite{oikhberg_etal2022}.

\subsection{Construction and basic properties}

Let $E$ be a Banach space. Let $H[E]$ denote the vector sublattice of
$\R^{E^{*}}$ consisting of all positively homogeneous functions
$f\colon E^{*}\to \R$ (that is, functions satisfying $f(\lambda
x^{*})=\lambda f(x^{*})$ for all $\lambda \ge 0$ and $x^{*} \in
E^{*}$). Given $f \in H[E]$, define the quantity
\[
\|f\|=\sup \left\{\, \sum_{k=1}^{n}|f(x_k^{*})| : n \in \N,
x_1^{*},\ldots ,x_n^{*}\in E^{*}, \sup_{x \in
B_E}\sum_{k=1}^{n}|x_k^{*}(x)|\le 1 \,\right\}.
\]
One can readily check that
\[
    H_1[E]=\{\, f :  \|f\|<\infty \, \}
\]
is a sublattice of $H[E]$, and that $\|{\cdot }\|$ defines a complete
lattice norm on $H_1[E]$.

For every $x \in E$, define $\delta _x \in H[E]$ by $\delta
_x(x^{*})=x^{*}(x)$ for $x ^{*}\in E^{*}$. Note that $\|\delta
_x\|=\|x\|$ for every $x \in E$. Let $\FVLv E$ denote the sublattice
generated by $\{\, \delta _x : x \in E \, \} $ in $H[E]$, and let
$\delta_E \colon E\to \FVLv E$ be the linear map defined by $\delta_E (x)=\delta _x$ for
all $x \in E$. The space
$\FVLv E$ consists of all expressions that can be built from
finitely many elements of the form $\delta _x$ using finitely many
linear and lattice operations. The pair $(\FVLv E,\delta _E)$ has the universal property of the
free vector lattice generated by $E$: for every linear map $T\colon
E\to X$ into a vector lattice $X$, there exists a unique lattice
homomorphism $\hat{T}\colon \FVLv E\to X$ such that $T=\hat{T}\circ\delta_E
$.

Let $\FBL E$ be the norm closure of $\FVLv E$ in $H_1[E]$. The map $\delta
_E$ naturally extends to a linear isometry $\phi _E\colon E\to \FBL
E$. The space $\FBL E$, together with the map $\phi _E$, satisfies the
universal property of the free Banach lattice:

\begin{thm}[{\cite[Theorem 2.1]{oikhberg_etal2022}}]
    Let $X$ be a Banach lattice, and let $T\colon E\to X$ be an
    operator. There is a unique lattice homomorphism $\hat{T}\colon
    \FBL E\to X$ such that $\hat{T}\circ \phi _E=T$ and
    $\|\hat{T}\|= \|T\|$.
\end{thm}

The results we will cite from \cite{oikhberg_etal2022} are always
stated in greater generality, because in that paper the object of study is the free $p$-convex Banach lattice. In this text we will not work with convexity
conditions, so our interest will focus on the case $p=1$.

The elements of $\FBL E$ are weak$^*$ continuous on bounded sets.
Since they are also positively homogeneous, we can isometrically embed
$\FBL E$ as a sublattice of $C(B_{E^{*}})$. Most of the time we will
regard the elements of $\FBL E$ as functions on $B_{E^{*}}$ using this
identification. Moreover, if $E$ is finite-dimensional, then the
restriction of the functions of $\FBL E$ to the unit sphere
$S_{E^{*}}$ is onto. Hence, when $E$ is finite-dimensional, the free
Banach lattice generated by $E$ can be identified with $C(S_{E^{*}})$.
This is the only case in which the free Banach lattice has a strong
unit.

\begin{prop}[{\cite[Propositions 9.1 and 9.4]{oikhberg_etal2022}}]\label{prop:FBL_strong_unit}
    Let $E$ be a Banach space.
    \begin{enumerate}
        \item $E$ is finite-dimensional if and only if $\FBL E$ has a
            strong unit.
        \item $E$ is separable if and only if $\FBL E$ has a
            quasi-interior point.
    \end{enumerate}
\end{prop}

The following are some basic properties of the free Banach lattice.

\begin{prop}[{\cite[Proposition 2.11]{oikhberg_etal2022}}]\label{prop:FBL_basic}
    Let $E$ be a Banach space.
    \begin{enumerate}
        \item For every $x \in E$, $x\neq 0$, $|\delta _x|$ is a weak
            order unit in $\FBL E$.
        \item If $E$ has dimension strictly greater than one, then
            $\FBL E$ has no non-trivial projection bands.
        \item If $E$ has dimension strictly greater than one, then
            $\FBL E$ is not $\sigma $-order complete and has no atoms.
    \end{enumerate}
\end{prop}

\subsection{Operators and their extensions}

Every bounded linear operator $T\colon E\to F$ between Banach spaces
extends uniquely to a lattice homomorphism $\bar{T}\colon \FBL E\to
\FBL F$ making the following diagram commute:
\[
\begin{tikzcd}
    \FBL E\arrow[r, "\bar{T}"]&\FBL F\\
    E\arrow[u, "\phi _E"]\arrow[r, "T"]&F\arrow[u, "\phi _F"]
\end{tikzcd}
\]
To be more precise, $\bar{T}=\widehat{\phi _F T}$; in particular,
$\|\bar{T}\|=\|T\|$. One can even provide an explicit expression of
$\bar{T}$.

\begin{lem}[{\cite[Lemma 3.1]{oikhberg_etal2022}}]\label{lem:FBL_extension}
    Given an operator $T\colon E\to F$ between Banach spaces, the
    extension $\bar{T}\colon \FBL E\to \FBL F$ is given, for $f \in
    \FBL E$, by
    \[
        \bar{T}(f)=f\circ T^{*}.
    \]
\end{lem}

There exists a close relation between the properties of $T$ and
$\bar{T}$.

\begin{prop}[{\cite[Proposition 3.2]{oikhberg_etal2022}}]\label{prop:FBL_extension_properties}
    Let $T\colon E\to F$ be an operator, and let $\bar{T}\colon \FBL
    E\to \FBL F$ be its extension. Then
    \begin{enumerate}
        \item $T$ is injective if and only if $\bar{T}$ is injective.
        \item $T$ has dense range if and only if $\bar{T}$ has dense
            range.
        \item $T$ is onto if and only if $\bar{T}$ is onto.
    \end{enumerate}
\end{prop}

If $T\colon E\to F$ is an isomorphism (resp.\ isometry) between Banach
spaces, then $\bar{T}\colon \FBL E\to \FBL F$ is a lattice isomorphism
(resp.\ isometry). Thus isomorphic (resp.\ isometric) Banach spaces
generate isomorphic (resp.\ isometric) free Banach lattices. Whether
the converse holds is an important open question. For the
isometric version, the answer is affirmative whenever $E$ and $F$ have
smooth duals.

\begin{thm}[{\cite[Theorem 10.18]{oikhberg_etal2022}}]\label{thm:FBL_isometries}
    Let $E$ and $F$ be Banach spaces with smooth dual. A map $T\colon
    \FBL E\to \FBL F$ is a surjective lattice isometry if and only if
    $T=\bar{U}$, for some isometry $U\colon E\to F$.
    Consequently, $E$ and $F$ are isometric if and only if $\FBL E$ is
    lattice isometric to $\FBL F$.
\end{thm}

This result uses the following description of lattice homomorphisms
between free Banach lattices.

\begin{prop}[{\cite[Proposition 10.1]{oikhberg_etal2022}}]
    Let $E$ and $F$ be Banach spaces, and let $T\colon \FBL E\to \FBL
    F$ be a lattice homomorphism. Then there exists a map $\Phi
    _T\colon F^{*}\to E^{*}$ such that $Tf=f\circ \Phi _T$ for every
    $f \in \FBL E$. Moreover, $\Phi _T$ satisfies the following
    properties:
    \begin{enumerate}
        \item For every $y^{*}\in F^{*}$ and $x \in E$, $\Phi
            _T(y^{*})(x)=(T\delta _x)(y^{*})$.
        \item $\Phi _T$ is positively homogeneous.
        \item $\Phi _T$ is weak$^*$ to weak$^*$ continuous on bounded
            sets.
        \item For $y^{*}\in F^{*}$, we have $\|\Phi _Ty^{*}\|\le \|T\|
            \|y^{*}\|$.
        \item For every $y_1^{*},\ldots ,y_m^{*}\in F^{*}$ we have
            \[
            \sup_{x \in B_E} \sum_{k=1}^{m}|\Phi _T(y_k^{*})(x)| \le
            \|T\| \sup_{y \in B_F} \sum_{k=1}^{m} |y_k^{*}(y)|.
            \]
    \end{enumerate}
\end{prop}

\subsection{Projective and injective Banach lattices}\label{sec:background_projective}

Projective Banach lattices were also introduced for the first time in
\cite{de_pagter_wickstead2015}.

\begin{defn}
    Let $\lambda >1$ be a real number. A Banach lattice $P$ is
    \emph{$\lambda $-projective} if whenever $X$ is a Banach lattice,
    $J$ a closed lattice ideal in $X$, and $q\colon X\to X/J$
    the quotient map, then for every lattice homomorphism
    $T\colon P\to X/J$ there is a lattice homomorphism
    $\hat{T}\colon P\to X$ with $T=q \hat{T}$ and $\|\hat{T}\|\le
    \lambda \|T\|$.
    A Banach lattice is \emph{projective} if it is $(1+\varepsilon
    )$-projective for every $\varepsilon >0$.
\end{defn}

Projective Banach lattices can be characterized using free Banach
lattices as follows.

\begin{thm}[{\cite[Theorem 10.3]{de_pagter_wickstead2015}}]
    A Banach lattice $P$ is projective if and only if for every
    $\varepsilon >0$ there exist an
    index set $\Gamma$, a closed sublattice $H$ of $\FBL{\el 1(\Gamma
    )}$, a lattice isomorphism $T\colon P\to H$ with
    $\|T\|,\|T^{-1}\|\le 1+\varepsilon $, and a lattice homomorphism $R\colon
    \FBL{\el 1(\Gamma )}\to H$ that is a projection and has norm
    $\|R\|\le 1+\varepsilon $.
\end{thm}

In particular, $\FBL{\el 1(\Gamma )}$ is projective for every set
$\Gamma $. Further natural examples of projective Banach lattices are
the following.

\begin{thm}[{\cite[Theorem 11.1]{de_pagter_wickstead2015}}]
    Every finite-dimensional Banach lattice is projective.
\end{thm}

More generally, we have the following characterization of the $C(K)$
spaces that are projective.

\begin{thm}[{\cite[Theorem 1.4]{aviles_martinez-cervantes_abellan2020}}]
    Let $K$ be a compact Hausdorff space. Then $C(K)$ is projective if
    and only if $K$ is an absolute neighbourhood retract.
\end{thm}

The free Banach lattices that are projective must be close to $\el 1$.

\begin{thm}[{\cite[Theorem 1.3]{aviles_martinez-cervantes_abellan2020}}]
    Let $E$ be a Banach space. If $\FBL E$ is $\lambda $-projective,
    for some $\lambda >1$, then $E$ has the Schur property.
\end{thm}

The notion dual to projectivity is that of injectivity. However, as
already noted in \cite[Footnote 1]{de_pagter_wickstead2015}, there are
no injective objects in the category of Banach lattices and lattice
homomorphisms. For this reason, injective Banach lattices have
traditionally been defined with respect to positive operators
\cite{haydon1977}. More recent work \cite{aviles_tradacete2023} has also
considered injective objects in the category whose objects are
separable Banach lattices and whose arrows are lattice homomorphisms.

\chapter{Free Archimedean \texorpdfstring{$f$-algebras}{f-algebras}}
\label{sec:archimedean}

This chapter is devoted to the construction and study of the following
object.

\begin{defn}
    Let $S$ be a set. The \emph{free (Archimedean) \falg\ generated
    by $S$} is an (Archimedean) \falg\ $\FAFA S$ together with a map
    $\delta \colon S\to \FAFA S$ such that, for every
    (Archimedean)
    \falg\ $A$ and every map $T\colon S\to A$, there exists a
    unique lattice-algebra homomorphism $\hat{T}\colon \FAFA S\to A$
    satisfying $\hat T \circ \delta =T$.
\end{defn}

The existence of free \falg s has been long known in the field of
universal algebras \cite[Section 5.2]{keimel1995}: it follows from the fact that \falg s form an
equationable class \cite[Theorem VII.8]{birkhoff1967}. Archimedean
\falg s, however, do not form an equationable class \cite[Section 8]{de_jeu2021}.
Henriksen and Isbell showed in
\cite{henriksen_isbell1962} that every LLA expression (that
is, every expression involving finitely many variables, linear and
lattice operations, and a product) that vanishes on the reals, also
vanishes on every Archimedean \falg. It then follows from a standard argument due to
Birkhoff \cite{birkhoff1935} that the free
Archimedean \falg\ exists.

In \cref{sec:construction_FAFA} we revisit the result of Henriksen and Isbell through
a completely new lens. Their proof and formulation depended heavily on
universal algebra language and techniques that may be unfamiliar to
the modern analyst. Here we shall use a more operator-centric
approach together with a representation theorem
due to Henriksen and Johnson \cite{henriksen_johnson1961}.

Admittedly, this may not seem very original. Yet this approach will be
used in \cref{sec:structure} to show that a refined version of the
result holds for Banach \falg s. This result will, in turn, be central
in the construction of free Banach \falg s, the main topic of this
paper.

In \cref{sec:lat_prop_FAFA,sec:alg_prop_FAFA} the lattice and algebraic
properties of the free Archimedean \falg\ are explored. This work is
the analogue of what Baker \cite{baker1968} and Bleier
\cite{bleier1973} did for free vector lattices.
The reader familiar with free vector lattices
will find the properties of free Archimedean
\falg s surprisingly similar and, at the same time, different enough
to deserve a separate study.

Even though free \falg s date back to 1956, several
open conjectures regarding free (Archimedean) \falg s remain (see
\cite{keimel1995} for a survey on the topic). Some of these
have received recent attention, especially the Birkhoff--Pierce
conjecture (see for instance
\cite{madden2011,mahe2007,wagner2010,lucas_schaub_spivakovsky2015}).
Free \falg s also play a prominent role in the
study of semialgebraic geometry (see \cite[Section 5]{keimel1995} for
an introduction and further references).

\section{Construction of free Archimedean \falg
s}\label{sec:construction_FAFA}

Informally, a lattice-linear-algebraic (LLA) expression is any
expression constructed using finitely many variables, the lattice,
linear and multiplication operations. An LLA expression is said to vanish
on a subset $B$ of an \falg\ if ``evaluating'' the expression at arbitrary
elements of $B$ always gives zero. This section is devoted to proving
the following.

\begin{thm}\label{thm:falgYudin}
    Let $\Phi $ be an LLA expression. If $\Phi $ vanishes on $\R$,
    then it also vanishes on every Archimedean \falg.
\end{thm}

Once this result is established, it will follow from a standard
argument that the free Archimedean \falg\ generated by a set $S$ is
the vector lattice algebra of $\R^{\R^{S}}$ generated by the
evaluations $\delta _s$, where $\delta _s(x)=x(s)$ for all $x \in
\R^{S}$ and $s \in S$.
But before proceeding to the proof of \cref{thm:falgYudin}, a more
precise definition of LLA expressions is needed.

\subsection{LLA expressions}

Let $\Omega $ be a set of formal operations and let $T$ be a set of
formal variables. We are going to describe the procedure to construct formal
expressions on the variables of $T$ using the operations from
$\Omega $. The definition of such expressions is inductive on the
complexity of the expression. Define every element of $T$ and the
0-ary operations
to be expressions of complexity $1$. If $\omega \in \Omega $ is a
$k$-ary operation, for some $k\ge 1$, and $\Phi _1,\ldots ,\Phi _k
$ are expressions of complexity $n-1$, then $\omega (\Phi _1,\ldots
,\Phi _k)$ is an expression of complexity $n$. By construction, every formal
expression is a word in the alphabet $T\cup \Omega $
containing only finitely many elements from $T$. The
fact that at most $t_1,\ldots ,t_n \in T$ appear in $\Phi $ is made
explicit by writing $\Phi [t_1,\ldots ,t_n]$. Denote by
$\mathcal{E}(T,\Omega )$ the set of formal expressions on the
variables of $T$ and using the operations from $\Omega $. The formal
operations with words in $\mathcal{E}(T,\Omega )$ make it into an
$\Omega $-algebra (see \cite[Section 3]{de_jeu2021} for a more
detailed introduction to the language of universal algebras).

Let $A$ be an $\Omega $-algebra. Given a map $\phi \colon T\to A$ it
extends, by induction on complexity, to a unique $\Omega $-algebra
homomorphism $\tilde \phi \colon \mathcal{E}(T,\Omega )\to A$. If
$\phi (t)=a_t$, denote $\tilde \phi (\Phi [t_1,\ldots ,t_n])$ by
$\Phi (a_{t_1},\ldots ,a_{t_n})$.

When $\Omega $ consists of a 0-ary operation $0$, a unary operation
$\lambda $ for every $\lambda \in \R$, and two binary operations
$+$ and $\vee$, the elements of $\mathcal{E}(T,\Omega )$ are called
\emph{lattice-linear (LL) expressions} in $T$. Note that, for this
particular $\Omega $,
every vector lattice is an $\Omega $-algebra. If $\Omega $ contains an
additional binary operation $\cdot $, the elements of
$\mathcal{E}(T,\Omega )$ are called \emph{lattice-linear-algebraic (LLA)
expressions}. In this case, every vector lattice algebra is an $\Omega
$-algebra. Observe that every LL expression can also be regarded as an
LLA expression where $\cdot $ does not appear. Conversely, from every
LLA expression $\Phi [t_1,\ldots ,t_n]$ one can produce an LL expression
$\Phi _0[t_1,\ldots ,t_n]$ in the following way: repeat the same steps
as in the inductive construction of $\Phi $ except that, every time the
product operation is performed, the resulting expression is set to
$0$. This way, only operations $0$, $\lambda $, $+$ and $\vee $ appear
in $\Phi _0$, which is therefore a lattice-linear expression. An LLA
expression $\Phi[t_1,\ldots ,t_n]$ is said to \emph{vanish} on a
subset $B$ of a vector lattice algebra $A$ if $\Phi
(a_1,\ldots ,a_n)=0$ for every $a_1,\ldots ,a_n \in B$.

\subsection{\Cref{thm:falgYudin} for \falg s with identity}

The proof of \cref{thm:falgYudin} is divided into several steps: first,
we show the result is true when an identity is present; then we show
it is true when the product is identically zero; and finally we prove
the general case.

\begin{lem}\label{lem:calc_identity}
    Let $\Phi $ be an LLA expression. If $\Phi $ vanishes on $\R$,
    then it also vanishes on every Archimedean \falg\ with identity.
\end{lem}
\begin{proof}
    By \cref{thm:repr_identity}, we may identify $A$ with an \falg\
    in $C^{\infty }(K)$, for a certain compact Hausdorff space $K$.
    For $f \in C^{\infty }(K)$, let $D(f)$ denote the open dense set
    of points at which $f$ is real-valued. Let $a_1,\ldots ,a_n \in
    A$. By definition of the operations in $C^{\infty }(K)$:
    \[
    \Phi (a_1,\ldots ,a_n)(t)=\Phi (a_1(t),\ldots
    ,a_n(t))=0\quad\text{for every }t \in D(a_1)\cap \cdots \cap
    D(a_n).
    \]
    Since $D(a_1)\cap \cdots \cap D(a_n)$ is dense in $K$, it follows
    that $\Phi (a_1,\ldots ,a_n)=0$.
\end{proof}

\subsection{\Cref{thm:falgYudin} for \falg s with the zero product}

It is a well-known fact that every lattice-linear expression vanishing on $\R$
must also vanish on every vector lattice (see, for instance,
\cite{baker1968}). This will be key in showing that
\cref{thm:falgYudin} holds for \falg s with the zero product. The
proof of this fact is done in two steps: first, it is
shown that an LLA expression vanishing on $\R_0$ must vanish on
every vector lattice with the zero product
(\cref{lem:calc_zero_step1}); then it is shown that every LLA
expression vanishing on $\R$ must also vanish on $\R_0$
(\cref{lem:calc_zero_step2}).

\begin{lem}\label{lem:calc_zero_step1}
    Let $\Phi$ be an LLA expression and let $X$ be a
    vector lattice. If $\Phi$ vanishes on $\R_0$, then it also vanishes
    on $X_0$.
\end{lem}
\begin{proof}
    At the core of this proof is the fact that $\Phi (x_1,\ldots ,x_n)=\Phi _0 (x_1,\ldots ,x_n)$ for
    every $x_1,\ldots ,x_n \in X_0$. To show this, proceed by induction on $k$, the
    complexity of $\Phi $. If $k=1$, then $\Phi =\Phi_0$ and the result
    is clear. Suppose it is true for LLA expressions of complexity up to
    $k$, and let $\Phi $ have complexity $k+1$. If either
    $\Phi =\Psi + \lambda \Upsilon  $ or $\Phi =\Psi \vee \Upsilon  $ for some LLA
    expressions $\Psi $ and $\Upsilon $ of complexity up to $k$ and $\lambda \in
    \R$, then $\Phi _0=\Psi _0+\lambda \Upsilon_0$ or $\Phi _0=\Psi
    _0\vee \Upsilon_0$, and the desired result follows from applying
    the induction hypothesis to $\Psi $ and $\Upsilon$.
    If instead $\Phi =\Psi \Upsilon $, then by definition $\Phi _0=0$
    and, since the product in $X_0$ is identically zero, also
    \[
    \Phi (x_1,\ldots ,x_n)=\Psi (x_1,\ldots ,x_n)\Upsilon (x_1,\ldots
    ,x_n)=0.
    \]

    In particular, this applies to $\R_0$. Hence $\Phi _0 (\lambda
    _1,\ldots ,\lambda _n)=\Phi (\lambda _1,\ldots ,\lambda _n)=0$ for
    every $\lambda _1,\ldots ,\lambda _n \in \R_0$. Since $\Phi _0$ is
    a
    lattice-linear expression that vanishes on \R, it must also vanish
    on every vector lattice. It
    follows that $0=\Phi _0 (x_1,\ldots ,x_n)=\Phi (x_1,\ldots ,x_n)$ for
    every $x_1,\ldots ,x_n \in X_0$.
\end{proof}

\begin{lem}\label{lem:calc_zero_step2}
    Let $\Phi $ be an LLA expression. If $\Phi $ vanishes on a
    neighbourhood of $0$ in $\R$, then it also vanishes on $\R_0$.
\end{lem}
\begin{proof}
    Note that $\Phi =\Phi [t_1,\ldots ,t_n]$ naturally determines a
    continuous function $\R^{n}\to \R$ defined by $(\lambda _1,\ldots
    ,\lambda _n)\mapsto \Phi (\lambda _1,\ldots ,\lambda _n)$. The
    same is true of $\Phi _0$. We are going to prove that
    \begin{equation}\label{eq:lla_limit}
    \lim_{\varepsilon  \to 0^{+}}\frac{\Phi (\varepsilon \lambda
    _1,\ldots ,\varepsilon \lambda _n)}{\varepsilon }=\Phi _0 (\lambda
    _1,\ldots ,\lambda _n)
    \end{equation}
    holds for every $(\lambda _1,\ldots ,\lambda _n)\in \R^{n} $
    by induction on the complexity of $\Phi $.

    If $\Phi $ has complexity $1$, then \eqref{eq:lla_limit} is clear.
    Assume that the result is true if $\Phi $ has complexity up to $k-1$.
    Let $\Phi
    $ have complexity $k$. There exist LLA expressions $\Psi
    [t_1,\ldots ,t_n]$ and $\Upsilon[t_1,\ldots ,t_n]$ of complexity
    up to $k-1$ such that either:
    \begin{enumerate}
        \item $\Phi =\lambda \Psi + \Upsilon $ for some $\lambda \in
            \R$, in which case $\Phi _0=\lambda \Psi _0+\Upsilon _0$ and
            \begin{align*}
                \lim_{\varepsilon  \to 0^{+}} &\frac{\lambda \Psi
            (\varepsilon \lambda _1,\ldots ,\varepsilon \lambda
        _n)+\Upsilon (\varepsilon \lambda _1,\ldots ,\varepsilon \lambda
                _n)}{\varepsilon }\\&=\lambda\lim_{\varepsilon  \to 0^{+}} \frac{
        \Psi (\varepsilon \lambda _1,\ldots ,\varepsilon \lambda
    _n)}{\varepsilon }+\lim_{\varepsilon  \to 0^{+}} \frac{ \Upsilon
            (\varepsilon \lambda _1,\ldots ,\varepsilon \lambda
    _n)}{\varepsilon }
                    \\&=\lambda \Psi _0(\lambda _1,\ldots ,\lambda
    _n)+\Upsilon _0(\lambda _1,\ldots ,\lambda _n).
            \end{align*}
        \item $\Phi =\Psi \vee \Upsilon $, in which case $\Phi _0 =\Psi_0
            \vee \Upsilon_0 $ and
        \begin{align*}
            \lim_{\varepsilon  \to 0^{+}} &\frac{
            \Psi (\varepsilon \lambda _1,\ldots ,\varepsilon \lambda
        _n)\vee \Upsilon (\varepsilon \lambda _1,\ldots ,\varepsilon \lambda
    _n)}{\varepsilon }\\&=\bigg(\lim_{\varepsilon  \to 0^{+}} \frac{ \Psi
            (\varepsilon \lambda _1,\ldots ,\varepsilon \lambda
    _n)}{\varepsilon }\bigg)\vee \bigg(\lim_{\varepsilon  \to 0^{+}}
    \frac{ \Upsilon
            (\varepsilon \lambda _1,\ldots ,\varepsilon \lambda
    _n)}{\varepsilon }\bigg)
                    \\&= \Psi _0(\lambda _1,\ldots ,\lambda
    _n)\vee \Upsilon _0(\lambda _1,\ldots ,\lambda _n).
            \end{align*}
        \item $\Phi =\Psi \cdot \Upsilon$, in which case $\Phi _0 =0$
            and
    \begin{align*}
        \lim_{\varepsilon  \to 0^{+}} &\frac{\Psi
            (\varepsilon \lambda _1,\ldots ,\varepsilon \lambda
        _n)\cdot \Upsilon (\varepsilon \lambda _1,\ldots ,\varepsilon \lambda
    _n)}{\varepsilon }\\&=\bigg(\lim_{\varepsilon  \to 0^{+}} \frac{ \Psi
            (\varepsilon \lambda _1,\ldots ,\varepsilon \lambda
    _n)}{\varepsilon }\bigg) \bigg(\lim_{\varepsilon  \to 0^{+}} \Upsilon
            (\varepsilon \lambda _1,\ldots ,\varepsilon \lambda
        _n)\bigg)=0,
            \end{align*}
        where in the last equality we are using that $(\lambda _1,\ldots
        ,\lambda _n)\mapsto \Upsilon(\lambda _1,\ldots ,\lambda _n)$ defines a continuous
        function on $\R^{n}$, and that $\Upsilon(0,\ldots ,0)=0$ by
        construction.
    \end{enumerate}

    Now suppose there exists a $\delta >0$ such that $\Phi (\mu
    _1,\ldots ,\mu  _n)=0$ for every $(\mu  _1,\ldots ,\mu  _n)
    \in [-\delta ,\delta ]^{n}$. Then the limit in \eqref{eq:lla_limit} must vanish, and
    therefore $\Phi _0 (\lambda  _1,\ldots ,\lambda  _n)=0$ for
    every
    $\lambda _1,\ldots ,\lambda _n \in \R$. It was shown in
    the proof of \cref{lem:calc_zero_step1} that, when evaluated in $\R_0$,
    $\Phi (\lambda _1,\ldots , \lambda _n)=\Phi _0(\lambda
    _1,\ldots ,\lambda _n)=0$. Hence $\Phi $ vanishes on $\R_0$.
\end{proof}

Equation \eqref{eq:lla_limit} will be of great use in later sections.
For convenience, we collect it in a separate lemma and show that the
limit is uniform on compact subsets of $\R^{n}$.

\begin{lem}\label{lem:lla_limit}
    For every LLA expression $\Phi [t_1,\ldots ,t_n]$ and $(\lambda
    _1,\ldots ,\lambda _n) \in \R^{n}$:
    \[
    \lim_{\varepsilon  \to 0^{+}}\frac{\Phi (\varepsilon \lambda
    _1,\ldots ,\varepsilon \lambda _n)}{\varepsilon }=\Phi _0 (\lambda
    _1,\ldots ,\lambda _n).
    \]
    Moreover, the limit is uniform on compact subsets of $\R^{n}$.
\end{lem}
\begin{proof}
    The proof is based on the following elementary fact:
    Let $K$ be a compact metric space and let $(f_n)$ and $(g_n)$ be
    sequences of real-valued continuous functions defined on $K$ such
    that $f_n\to f$ and $g_n\to g$ uniformly, for some $f,g \in C(K)$.
    Let $h\colon \R^2\to \R$ be a continuous function. Then $h\circ
    (f_n\times g_n)\to h\circ (f\times g)$ uniformly, where $f\times g
    \colon K\times K\to \R^2$ is defined by $(f\times
    g)(k_1,k_2)=(f(k_1),g(k_2))$, $(k_1,k_2)\in K\times K$.

    Let us check this. By keeping only a tail of the sequences, we may
    assume that $\|f_n-f\|_\infty, \|g_n-g\|_{\infty
    }<1$ hold for all $n \in \N$. Let $[a,b]$ be an interval in
    $\R$ containing the images of $f$ and $g$. Then
    $f_n(K),g_n(K)\subseteq [a-1,b+1]$ for all
    $n \in \N$. Being $[a-1,b+1 ]^2$ a compact
    subset of $\R^2$, the continuous function $h$ is uniformly
    continuous on it.
    Fix an arbitrary $\varepsilon >0$.
    Let $\delta >0$ be such that $|h(x,y)-h(x',y')|<\varepsilon $
    whenever $|x-x'|,|y-y'|<\delta $ and $x,x',y,y' \in [a-1,b+1]$.
    Let $N \in \N$ be such that $\|f_n-f\|_\infty, \|g_n-g\|_{\infty
    }<\delta $ whenever $n\ge N$. It follows that for every $k_1,k_2 \in K$:
    \[
        |h(f_n(k_1),g_n(k_2))-h(f(k_1),g(k_2))|<\varepsilon
    \]
    whenever $n\ge N$. This proves the fact.

    Now let $K\subseteq \R^{n}$ be a compact subset. We will show that
    the limit is uniform on $K$ by induction on the complexity of $\Phi $. If
    $k=1$, then $\Phi =\Phi _0$ and the limit is uniform on $K$. Suppose the
    limit is uniform on $K$ whenever the LLA expression has complexity up to $k$.
    Let $\Phi $ be a LLA expression of complexity $k+1$. Then there
    exist LLA expressions $\Psi $ and $\Upsilon $ of complexity up to $k$, and
    a continuous function $h\colon \R^2\to \R$, such that
    \[
    \Phi (\lambda _1,\ldots ,\lambda _n)=h(\Psi (\lambda _1,\ldots
    ,\lambda _n), \Upsilon (\lambda _1,\ldots ,\lambda _n)).
    \]
    Of course, this continuous function $h$ can only be a linear function,
    the supremum or the multiplication. Let $(r_m)$ be an arbitrary
    sequence of positive real numbers decreasing to zero. If $h$ is
    positively homogeneous (i.e., if it is a linear function or the
    supremum), then
    \[
        \frac{\Phi (r_m\lambda _1,\ldots ,r_m\lambda
        _n)}{r_m}=h\bigg(\frac{\Psi (r_m\lambda _1,\ldots
        ,r_m\lambda _n)}{r_m}, \frac{\Upsilon (r_m\lambda _1,\ldots
        ,r_m\lambda _n)}{r_m}\bigg).
    \]
    By induction hypothesis
    \[
    \bigg(\frac{\Psi (r_m\lambda _1,\ldots
        ,r_m\lambda _n)}{r_m}\bigg)_m\text{ and } \bigg(\frac{\Upsilon (r_m\lambda _1,\ldots
        ,r_m\lambda _n)}{r_m}\bigg)_m
    \]
    converge to $\Psi _0(\lambda _1,\ldots ,\lambda _n)$ and
    $\Upsilon_0(\lambda _1,\ldots ,\lambda _n)$, respectively, uniformly
    on $(\lambda _1,\ldots ,\lambda _n) \in K$. By the initial
    observation, the sequence
    \[
         \bigg(\frac{\Phi (r_m\lambda _1,\ldots ,r_m\lambda
        _n)}{r_m}\bigg)_m
    \]
    also converges to $h(\Psi _0(\lambda _1,\ldots ,\lambda _n),
    \Upsilon_0(\lambda _1,\ldots ,\lambda _n))=\Phi _0(\lambda_1,\ldots
    ,\lambda _n)$ uniformly on $(\lambda
    _1,\ldots ,\lambda _n) \in K$.

    It remains to check the case when $h$ is the multiplication. In
    this case, $\Phi _0=0$. The sequence
    \[
        (\Upsilon (r_m\lambda _1,\ldots ,r_m\lambda _n))_m
    \]
    converges to zero uniformly on $(\lambda _1,\ldots ,\lambda _n) \in
    K$. Indeed, $\Psi $, when seen as a function defined on $\R^{n}$, is
    continuous. This means that, fixed an $\varepsilon >0$, there exists
    a $\delta >0$ such that $|\Upsilon (x_1,\ldots ,x_n)|<\varepsilon $
    whenever $\max_{i=1,\ldots ,n}|x_i|<\delta $. Let
    \[
    M=\max_{(\lambda _1,\ldots ,\lambda _n)\in K} |\lambda _i|.
    \]
    Then $|\Upsilon (r_m\lambda _1,\ldots ,r_m\lambda _n)|<\varepsilon $
    as long as $|r_m|<\delta /M$. This proves the claim.

    In the expression
    \[
    \frac{\Phi (r_m\lambda _1,\ldots ,r_m\lambda
    _n)}{r_m}=\frac{\Psi (r_m\lambda _1,\ldots ,r_m\lambda
        _n)}{r_m}\Upsilon (r_m\lambda _1,\ldots ,r_m\lambda _n)
    \]
    the left term of the product converges to $\Psi_0(\lambda _1,\ldots
    ,\lambda _n)$ while the right term converges to $0$, both uniformly on
    $K$. Again by the initial observation, it follows that the product
    converges uniformly to $0$.
\end{proof}

The desired result now follows directly from the previous lemmas.

\begin{lem}\label{lem:calc_zero}
    Let $\Phi $ be an LLA expression and let $X$ be a vector lattice. If $\Phi $
    vanishes on a neighbourhood of $0$ in $\R$, then it also vanishes on $X_0$.
\end{lem}

\subsection{\cref{thm:falgYudin} for Archimedean \falg
s}\label{sec:proof_falgYudin}

The only ingredient left to cook \cref{thm:falgYudin} is the
structural role played by the annihilator $N(A)$, already discussed
in \cref{sec:background_falg}: it is both an algebraic ideal and a band, and coincides with
the set of nilpotent elements.

\begin{proof}[Proof of \cref{thm:falgYudin}]
    The product in an Archimedean \falg\ $A$ extends to its order
    completion, making it an Archimedean \falg\ (see \cite[Section
    7]{huijsmans1991}). Since $A$ lattice-algebra embeds in its order
    completion, it suffices
    to prove the result for order complete Archimedean \falg s.

    From now on, assume $A$ is order complete. In this case, the
    annihilator $N(A)$ is a
    projection band. Let $P\colon A\to A$ be the band projection onto
    $N(A)$, and let $P^{d}$ be the band projection onto its disjoint
    complement, $N(A)^{d}$. The
    \falg\ condition implies that the product is band
    preserving, so $P^{d}xP^{d}y \in N(A)^{d}$ for every $x,y \in A$.
    From the identity
    \[
        xy=(Px+P^{d}x)(Py+P^{d}y)=P^{d}xP^dy
    \]
    follows both that $P(xy)=0=Px Py$ and
    $P^{d}(xy)=P^{d}xP^{d}y$. This shows that $P$ and $P^{d}$ are
    lattice-algebra homomorphisms from $A$ onto the \falg s $N(A)$ and
    $N(A)^{d}$, respectively.

    Thus, given $a_1,\ldots, a_n \in A$,
    \begin{align*}
        \Phi (a_1,\ldots ,a_n)&=P(\Phi (a_1,\ldots ,a_n))+ P^{d}(\Phi
        (a_1,\ldots ,a_n))\\&=\Phi (Pa_1,\ldots ,Pa_n)+\Phi (P^{d}a_1,\ldots
    ,P^{d}a_n),
    \end{align*}
    and we just need to show that $\Phi (Pa_1,\ldots ,Pa_n)=0$ in
    $N(A)$ and $\Phi (P^{d}a_1,\ldots, P^{d}a_n)=0$ in $N(A)^{d}$.
    Since the product in $N(A)$ is identically zero, the first follows
    immediately from \cref{lem:calc_zero}.

    It is direct to check that $N(N(A)^{d})=\{0\}$, so that the map
    $L\colon N(A)^{d}\to \orth{N(A)^{d}}$ is an injective
    lattice-algebra homomorphism. But $\orth{N(A)^{d}}$ is an
    Archimedean \falg\
    with identity; from \cref{lem:calc_identity} it follows that
    \[
    L(\Phi (P^{d}a_1,\ldots ,P^{d}a_n))=\Phi (LP^{d}a_1,\ldots
    ,LP^{d}a_n)=0.
    \]
    Since $L$ is injective, $\Phi (P^{d}a_1,\ldots ,P^{d}a_n)=0$, and
    this finishes the proof.
\end{proof}

\subsection{The free Archimedean \falg\ generated by a set}

Now that \cref{thm:falgYudin} has been established, we can construct
the free Archimedean \falg\ generated by a set $S$. A word on
notation: when we say that a concrete space with a concrete map is
\emph{the free object} we mean that they satisfy the required
universal property. This is reasonable because any other pair
satisfying the universal property will be isomorphic to the given one.
Once we have an explicit construction, we will always identify the
free object with that particular realization, and consequently denote
both of them the same way.

The function space $\R^{\R^{S}}$, equipped with pointwise
linear operations, order and product, becomes an Archimedean \falg. For every $s
\in S$, define $\delta _s\colon \R^{S}\to \R$ by $\delta _s(x)=x(s)$
for all $x \in \R^{S}$. Denote $\FAFA S=\VLA \{\, \delta _s : s \in S
\, \} \subseteq \R^{\R^{S}}$, and define the map $\delta \colon S\to
\FAFA S$ by $\delta (s)=\delta _s$.

\begin{thm}\label{thm:FAFA}
    The free Archimedean \falg\ over a set $S$ is $\FAFA S$ together
    with the map $\delta $.
\end{thm}
\begin{proof}
    This fact follows from \cref{thm:falgYudin} using a standard
    argument due to Birkhoff \cite{birkhoff1935}. For
    convenience of the reader, we sketch it here.

    Let $A$ be an Archimedean \falg, and let $T\colon S\to A$ be a
    map. For every $f \in \FAFA S$ there exists an LLA expression
    $\Phi [t_1,\ldots ,t_n]$ such that $f=\Phi (\delta _{s_1},\ldots
    ,\delta _{s_n})$ for some $s_1,\ldots ,s_n \in S$. Define
    $\hat{T}\colon \FAFA S\to A$ by $\hat{T}(f)=\Phi (Ts_1,\ldots
    ,Ts_n)$. We need to check that $\hat{T}$ is indeed well-defined. Suppose
    that also $f=\Psi (\delta _{t_1},\ldots ,\delta _{t_m})$ for some
    other LLA expression $\Psi $ and some $t_1,\ldots ,t_m \in S$.
    Then for every $x \in \R^{S}$:
    \[
    \Phi (x(s_1),\ldots ,x(s_n))=f(x)=\Psi (x(t_1),\ldots ,x(t_m)).
    \]
    Since $(x(s_1),\ldots ,x(s_n),x(t_1),\ldots ,x(t_m))$ takes all
    possible values in $\R^{n+m}$ as $x$ ranges through $\R^{S}$, it
    follows that the LLA expression $\Phi -\Psi $ vanishes on $\R$.
    By \cref{thm:falgYudin}, it must also vanish on $A$. In
    particular, $\Phi (Ts_1,\ldots ,Ts_n)=\Psi (Tt_1,\ldots ,Tt_m)$.
    Hence $\hat{T}$ is well-defined, and
    it is then clear from the definition
    that it is the unique lattice-algebra homomorphism satisfying
    $\hat{T}\delta =T$.
\end{proof}

\section{Lattice properties}\label{sec:lat_prop_FAFA}

In this section, we consider the properties of $\FAFA S$ as a vector
lattice and compare them with those of the free vector lattice $\FVL
S$.

Before proceeding to the study of such properties, a preliminary
result on the simplification of LLA expressions in \falg s is needed.

Recall that, whenever $A$ is a real algebra, and $B\subseteq A$ is a
subset, $\R[B]$ denotes the subalgebra of $A$ generated by $B$. In other words, $\R[B]$
contains the polynomials in the elements of $B$ with no constant term.
Whenever $X$ is a vector lattice, and $Y\subseteq X$ is a subset,
$\Lat Y$ denotes the sublattice generated by $Y$. It is standard
that every element $x \in \Lat Y$ can be written as
\[
x=\bigvee_{i=1} ^{n}y_i - \bigvee_{j=1} ^{m} z_j,
\]
where $y_1,\ldots ,y_n,z_1,\ldots ,z_m \in \barespn Y$.

A \emph{$d$-algebra} $A$ is a
vector lattice algebra in which $a(b\vee c)=(ab)\vee (ac)$ and $(b\vee
c)a=(ba)\vee (ca)$ hold for all $b,c \in A$ and $a \in A_+$. Every \falg\ is
in particular a $d$-algebra.

\begin{prop}\label{prop:simplificationLLA}
    If $A$ is a $d$-algebra and $S\subseteq A$ is such that $A=\VLA(S)$, then
    \[
    A=\Lat(\R[s_+, s_- : s \in S]).
    \]
\end{prop}
\begin{proof}
    The proof consists in showing that $B=\Lat(\R[s_+, s_- : s \in
    S])$ is a sublattice-algebra; since it contains $S$, and
    $A=\VLA(S)$, it will then follow that $A=B$.
    Obviously, $B$ is a sublattice, so it is only left to show that
    $B$ is closed under taking products. Observe that, given $x \in \R[s_+, s_- : s
    \in S]$, separating its positive and negative coefficients one can
    write $x=x_p-x_n$, where
    $x_p,x_n\ge 0$ and $x_p,x_n \in \R[s_+, s_- : s
    \in S]$. Now let $y \in B$, say
    \[
    y=\bigvee_{i=1} ^{m}u_i-\bigvee_{j=1} ^{m}v_j\quad\text{with
    }u_i,v_j \in \R[s_+,s_-:s \in S].
    \]
    Using the $d$-algebra condition:
    \begin{align*}
        xy&=(x_p-x_n)\bigg(\bigvee_{i=1} ^{m}u_i-\bigvee_{j=1}
        ^{m}v_j\bigg)\\
          &=x_p\bigvee_{i=1} ^{m}u_i - x_n\bigvee_{i=1} ^{m}u_i -
          x_p\bigvee_{j=1} ^{m}v_j+x_n\bigvee_{j=1} ^{m}v_j\\
          &=\bigvee_{i=1} ^{m}x_pu_i - \bigvee_{i=1} ^{m}x_nu_i -
          \bigvee_{j=1} ^{m}x_pv_j+\bigvee_{j=1} ^{m}x_nv_j.
    \end{align*}
    Since $x_pu_i,x_pv_j,x_n u_i, x_n u_j \in
    \R[s_+, s_- : s \in S]$, it follows that $xy \in B$.

    Now suppose $x=\bigvee_{k=1} ^{p} w_k$ for some $w_k \in
    \R[s_+, s_- : s \in S]$. The element $x$ can be rewritten as
    \[
    x=w_{1,p}\vee \bigvee_{k=2} ^{p} (w_k+w_{1,n}) -
    w_{1,n}=w-w_{1,n},
    \]
    where $w=w_{1,p}\vee \bigvee_{k=2} ^{p} (w_k+w_{1,n}) \in B$ is positive
    because $w_{1,p}\ge 0$. Then
    \begin{align*}
        xy&=w\bigvee_{i=1} ^{m}u_i-w\bigvee_{j=1} ^{m}v_j-w_{1,n} y\\
            &=\bigvee_{i=1} ^{m}wu_i-\bigvee_{j=1} ^{m}wv_j-w_{1,n} y,
    \end{align*}
    where $w u_i, wv_j, w_{1,n} y \in B$ because we have already
    checked that the product of an element of $\R[s_+, s_- : s \in S]$ by an
    element of $B$ is in $B$. Finally, suppose that $x$ is an
    arbitrary
    element of $B$, say $x=\bigvee_{k=1}
    ^{p} w_k-\bigvee_{l=1} ^{p}z_l$ for some $w_k, z_l \in
    \R[s_+, s_- : s \in S]$. Then
    \[
    xy=\bigg(\bigvee_{k=1} ^{p}w_k\bigg)y-\bigg(\bigvee_{k=1}
    ^{p}z_k\bigg)y \in B
    \]
    because each of the terms is an element of $B$ by the previous case.
    So $B$ is closed under multiplication, and therefore is a
    sublattice-algebra. Since it contains $S$, it must be $B=A$.
\end{proof}

Recall that in $\FVL S$, $\delta _s$ is a weak unit for every $s \in S$, and
$\FVL S$ has a strong unit if and only if $S$ is finite (in which
case $\sup_{s \in S}|\delta _{s}|$ is a strong unit). In $\FAFA S$
the situation is as follows.

\begin{prop}
    Let $S$ be a non-empty set.
    \begin{enumerate}
        \item The element $\delta _s$ is a weak unit of $\FAFA
    S$ for every $s \in S$.
        \item $\FAFA S$ does not have a strong unit.
    \end{enumerate}
\end{prop}
\begin{proof}
    \begin{enumerate}
        \item Suppose $f \in \FAFA S$ were such that $|f|\wedge
            |\delta _s|=0$. Write $f=\Phi (\delta _s,\delta
            _{s_1},\ldots ,\delta _{s_n})$, for a certain LLA
            expression $\Phi $ and $s_1,\ldots ,s_n \in S$. By
            evaluating the expression $|f|\wedge
            |\delta _s|=0$ at $x \in \R^{S}$ it follows that $0=f(x)=\Phi (x(s),x(s_1),\ldots
            ,x(s_n))$ whenever $x(s)\neq 0$. This
            implies that the continuous function
            \[
            \begin{array}{cccc}
            & \R^{n+1} & \longrightarrow & \R \\
                    & (\lambda _0,\lambda _1,\ldots ,\lambda _n) &
            \longmapsto & \Phi (\lambda _0,\lambda _1,\ldots ,\lambda
            _n) \\
            \end{array}
            \]
            vanishes as long as $\lambda _0\neq 0$. Hence it must be
            identically zero, and therefore $f=0$ by
            \cref{thm:falgYudin}.
        \item Suppose $e \in \FAFA S$ is a strong unit. The
            space of real sequences $\R^{\N}$ equipped with
            pointwise linear operations, order and product is an
            Ar\-chi\-me\-dean \falg. Fix
            $s_0\in S$ and let $T\colon S\to \R^{\N}$ be defined
            by $Ts_0=(n)_{n \in \N}$ and $Ts=0$ if $s\neq s_0$. Denote
            $a=Ts_0$. By construction, $\hat{T}(e) \in
            \VLA\{a\}=\Lat\{\R[a]\}$. More precisely,
            $\hat{T}(e)\in \Lat\{a,a^2,\ldots ,a^{k}\}$, for a certain
            $k \in \N$. Being finitely generated, $a+a^2+ \cdots
            +a^{k}$ is a strong unit of $\Lat\{a,a^2,\ldots ,a^{k}\}$,
            so $\hat{T}(e)\le \lambda
            (a+a^2+ \cdots +a^{k})$ for a certain $\lambda$. And since $a\le
            a^2\le \cdots \le a^{k}$, one arrives at $\hat{T}(e)\le k
            \lambda a^{k}$.
            By assumption, $e$ is a strong unit, so there exists $\mu > 0$
            such that $|\delta _{s_0}|^{k+1}\le \mu e$. Evaluating $\hat{T}$
            at both sides of the inequality yields $a^{k+1}\le k
            \lambda \mu a^{k}$ or, equivalently, $a\le k \lambda \mu
            $. This contradicts the fact that $a$ is unbounded.\qedhere
    \end{enumerate}
\end{proof}

If the set $S$ has more than one element, then $\FVL S$ does not have
non-trivial projection bands. In particular, it
does not have atoms and it is not $\sigma $-order complete. This is based on a
topological argument (see \cite{bleier1973}) which also
works for $\FAFA S$. We sketch it below, after introducing the
necessary notation.

\begin{defn}
    The \emph{support} of $f \in \R^{\R^{S}}$ is defined to be the set
    \[
    \supp(f)= \{\, x \in \R^{S} : f(x)\neq 0 \, \}.
    \]
    For a subset $F\subseteq \R^{\R^{S}}$ its \emph{support} is
    defined as
    \[
    \supp(F)=\bigcup_{f \in F} \supp(f).
    \]
\end{defn}

Equip $\R^{S}$ with the product topology. Then every element of $\FAFA S$ is a continuous function and therefore its support is an open set.

\begin{prop}\label{prop:lat_prop_not_one}
    Let $S$ be a set with more than one element. Then $\FAFA S$ does
    not have non-trivial projection bands. In particular, it does not
    have atoms and it is not $\sigma $-order complete.
\end{prop}
\begin{proof}
    Suppose $\FAFA S=B\oplus B^{d}$, with $B$ a non-trivial projection
    band. Denote by $P$ the projection onto $B$, and by $P^{d}$ its
    complementary projection. The first step is to show that
    $\supp(B)\cup \supp(B^{d})=\R^{S}\setminus\{0\}$. Let $x \in
    \R^{S}$, $x\neq 0$. Choose $s \in S$ such that $x(s)\neq 0$. Since
    $0\neq |\delta _s|(x)= |P\delta _s|(x)+|P^{d}\delta _s|(x)$,
    either $|P\delta _s|(x)\neq 0$ or $|P^{d}\delta _s|(x)\neq 0$. In
    any case, $x \in \supp(B)\cup \supp(B^{d})$. This proves
    $\supp(B)\cup \supp(B^{d})=\R^{S}\setminus\{0\}$.

    The second step is to show that $\supp(B)\cap 
    \supp(B^{d})=\emptyset$. If $x \in \supp(B)\cap 
    \supp(B^{d})$, there exist positive $f \in B$ and $g \in B^{d}$
    such that $f(x)>0$ and $g(x)>0$. Then $(f\wedge g)(x)=f(x)\wedge
    g(x)>0$, which contradicts the fact that $f\wedge g\in B\cap
    B^{d}=\{0\}$.

    We have thus shown that $\supp(B)$ and $\supp(B^{d})$ are a
    disconnection of $\R^{S}\setminus\{0\}$. This is absurd, since
    $|S|>1$.
\end{proof}

When the set $S$ has cardinality one, $\FVL 1$ can be identified with
$\R^2$, and therefore it is both order complete and atomic. For $\FAFA 1$ the
situation is not so straightforward.
Recall that $\FAFA 1=\VLA\{I\}\subseteq \R^{\R}$, where $I(x)=x$ for
all $x \in \R$.

\begin{prop}
    \begin{enumerate}
        \item The unique decomposition into non-trivial projection
            bands of $\FAFA 1$ is $\VLA\{I_+\}\oplus
            \VLA\{I_-\}$.
        \item $\FAFA 1$ has no atoms.
        \item $\FAFA 1$ is not $\sigma $-order complete.
    \end{enumerate}
\end{prop}
\begin{proof}
    \begin{enumerate}
        \item Let $B_+=\VLA\{I_+\}$ and $B_-=\VLA\{I_-\}$. It suffices
            to check that $B_+$ and $B_-$ are subspaces with
            $B_+\subseteq B_-^{d}$ satisfying $\FAFA
            1=B_+\oplus B_-$ as vector spaces (see \cite[Lemma
            1.2.8]{meyer-nieberg1991}). Certainly, $B_+$ and $B_-$ are
            subspaces. Since $B_+$ and $B_-$ have disjoint supports,
            $B_+\cap B_-=\{0\}$ and $B_+ \subseteq B_-^{d}$.
            According to \cref{prop:simplificationLLA}, every $f \in
            \FAFA 1$ can be written as
            \[
                f=\bigvee_{i=1}
                ^{n}[p_i^{+}(I_+)+p_i^{-}(I_-)]-\bigvee_{j=1}^{m}[q_j^{+}(I_+)+q_j^{-}(I_-)],
            \]
            where $p_i^{+},p_i^{-},q_j^{+},q_j^{-}$ are polynomials.
            Again, since the supports of $I_+$ and $I_-$ are disjoint,
            $f$ can be rewritten as
            \[
            f=\bigg[\bigvee_{i=1} ^{n}p_i^{+}(I_+)-\bigvee_{j=1}
            ^{m}q_j^{+}(I_+)\bigg]+\bigg[\bigvee_{i=1} ^{n}p_i^{-}(I_-)-\bigvee_{j=1}
            ^{m}q_j^{-}(I_-)\bigg].
            \]
            This proves $f \in B_+ + B_-$.

            If $\FAFA 1=B\oplus B^{d}$, where $B$ is a non-trivial
            band projection, then, as in \cref{prop:lat_prop_not_one},
            $\supp(B)$ and $\supp(B^{d})$ form a disconnection of
            $\R\setminus\{0\}$. So we may assume
            $\supp(B)=(0,\infty)$. Then $I_+ \in B$ and, since band
            projections are closed under multiplication,
            $B_+\subseteq B$. It follows easily that $B=B_+$, and
            similarly $B^{d}=B_-$.
        \item Let $a> 0$, say $a=\Phi (I)$ for some LLA expression
            $\Phi $. It is not difficult to
            show, by induction on the complexity of $\Phi $, that there
            exists an $m \in \N$ such that
            \[
            \lim_{x \to 0^{+}}\frac{\Phi (x)}{x^{m}}=\infty .
            \]
            On the other hand
            \[
            \lim_{x \to 0^{+}}\frac{x^{m}\wedge \Phi (x)}{x^{m}}=1,
            \]
            so it cannot be the case that $(I_+)^{m}\wedge a$ is a
            scalar multiple of $a$. Hence, $a$ is not an atom.
        \item Even though the sequence $f_n=(I_+)^{n}\wedge I_+$ is
            bounded below by $0$, it does not have an infimum. Indeed,
            suppose $f$ is a positive lower bound. Then $f(t)=0$ for $t\in
            [0,1)$ and, since $f$ is continuous, there must exist
            $\delta >0$ such that $0\le f(t)\le 1/2$ for $t \in
            [1 ,1+\delta ]$. Choose $N \in \N$ big enough so
            as to have
            \[
                1/2<(1+\delta )^{N}-(1+\delta ).
            \]
            By continuity, there exists $0<\delta '<\delta $ such that
            \[
                1/2<(1+\delta ')^{N}-(1+\delta ')<1.
            \]
            It then follows that $[((I_+)^{N}-I_+)_+\wedge I_+]\vee f$ is
            a lower bound for $(f_n)$ that is strictly bigger than
            $f$. Hence $(f_n)$ does not have an infimum.\qedhere
    \end{enumerate}
\end{proof}

\begin{rem}\label{rem:CCC}
    The free vector lattice $\FVL S$ has the \emph{countable chain
    condition}: a disjoint collection in $\FVL S$ is at most
    countable. The reason is purely topological: the supports of a disjoint collection of functions form a
    disjoint collection of open subsets of $\R^{S}$, and such a collection
    is at most countable (see \cite{ross_stone1964}). The very same argument also
    shows that $\FAFA S$ has the countable chain condition.
\end{rem}

\begin{rem}\label{rem:FVLinFAFA}
    Let $T$ be a subset of $S$. The free Archimedean \falg\ $\FAFA T$ may be viewed as the
    sublattice-algebra of $\FAFA S$ generated by $\delta _t$ for $t
    \in T$. The unique lattice-algebra homomorphism $P\colon \FAFA S\to \FAFA
    S$ determined by $P\delta _s=\delta _s$ if $s \in T$ and $P\delta
    _s=0$ if $s \not\in T$, defines a projection onto $\FAFA T$.

    Similarly, the free vector lattice $\FVL S$ may be viewed as the sublattice of $\FAFA S$
    generated by $\delta _s$ for $s \in S$. Equip $\FVL S$ with
    the zero product, which as usual we denote by ${\FVL S}_0$. The map $\delta \colon S\to {\FVL S}_0$
    extends to a unique lattice-algebra homomorphism $P\colon \FAFA S\to {\FVL S}_0$. Forgetting
    about the product on ${\FVL S}_0$, the map
    \[
    P\colon \FAFA S\to \FVL S \subseteq \FAFA S
    \]
    is a lattice homomorphic projection. The explicit expression of
    this projection is
    \[
        (Pf)(x)=
        \lim_{t\to 0^{+}}\frac{f(tx)}{t},\quad\text{where }f \in
        \FAFA S.
    \]
    Indeed, suppose $f=\Phi (\delta _{s_1},\ldots ,\delta _{s_n})$
    for some LLA expression $\Phi $ and $s_1,\ldots ,s_n \in S$. By
    \cref{lem:lla_limit},
    \[
    \lim_{t\to 0^{+}}\frac{f(tx)}{t}=\lim_{t\to 0^{+}} \frac{\Phi
    (tx(s_1),\ldots ,tx(s_n))}{t}=\Phi _0(\delta _{s_1},\ldots ,\delta
    _{s_n})(x).
    \]
    So the limit exists and, since $\Phi _0$ is a lattice-linear
    expression, $Pf \in \FVL S$. From the formula defining $Pf$ it is
    immediate to check that the map $f\mapsto Pf$ is a lattice
    homomorphic projection extending $\delta $. By uniqueness, it must
    be the projection introduced above.
\end{rem}

For lattice-linear expressions, there is continuity with respect to
the expressions in the following sense: if $(\Phi _k[t_1,\ldots
,t_n])_k$ is a sequence of lattice-linear expressions that converges
uniformly to $\Phi [t_1,\ldots ,t_n]$ when seen as functions in
$C([-1,1]^{n})$, then $\Phi _k(x_1,\ldots ,x_n)$ converges in norm to
$\Phi (x_1,\ldots ,x_n)$ for every $x_1,\ldots ,x_n$ in a Banach
lattice $X$. This is no longer true for LLA expressions.

\begin{example}
    Consider the sequence of LLA expressions
    \[
        \Phi_k [t]=t\bigg( \sum_{j=0}^{k} \frac{t^{2j}}{(2j)!} \bigg)
        ^2-t\bigg( \sum_{j=0}^{k} \frac{t^{2j+1}}{(2j+1)!} \bigg) ^2.
    \]
    For every $x \in [-1,1]$:
    \[
    \lim_{k\to \infty }\Phi _k(x)=x \cosh^2x - x\sinh^2x=x
    \]
    and this limit is uniform on $[-1,1]$. In other words, if $\Phi [t]=t$,
    and we see $\Phi _k$ and $\Phi $ as elements in $C([-1,1])$, $\Phi
    _k\to \Phi $ in the uniform norm. But, in $\R_0$ (or any other
    Banach lattice with the zero product), $\Phi _k(x)=0$ whereas
    $\Phi (x)=x$ for every $x \in \R_0$. Hence $\Phi _k(x)$ does not
    converge to $\Phi (x)$ if $x\neq 0$.
\end{example}

\section{Algebraic properties}\label{sec:alg_prop_FAFA}

In this section the algebraic properties of $\FAFA S$ are considered.
Recall that an \falg\ $A$ is \emph{semiprime} if $N(A)=\{0\}$.

\begin{prop}\label{prop:FAFA_no_identity}
    Let $S$ be a set.
    \begin{enumerate}
        \item $\FAFA S$ is semiprime.
        \item $\FAFA S$ does not have an identity element.
    \end{enumerate}
\end{prop}
\begin{proof}
    The result is trivial when $S$ is empty. For the rest of the
    proof, assume that $S$ is non-empty.
    \begin{enumerate}
        \item Suppose $f \in N(\FAFA S)$, and let $x \in \R^{S}$ be
            different from zero. There exists $s \in S$ such that
            $x(s)\neq 0$. Then $0=(f\delta _s)(x)=f(x)x(s)$, so
            $f(x)=0$. This shows that $f =0$.
        \item Suppose $1 \in \FAFA S$ were an identity. Fix $s_0\in S$
            and define $T\colon S\to C[-1,1]$ by $Ts=0$ for $s\neq
            s_0$ and $(Ts_0)(t)=t_+$ for $t \in [-1,1]$. Then
            $\hat{T}(1)=\hat{T}(1)^2$. Since
            $[-1,1]$ is connected, $\hat{T}(1)$ is either $0$ or
            $\one$, where $\one$ is the constant one
            function. But $\hat{T}(1)\hat{T}(\delta
            _{s_0})=\hat{T}(\delta _{s_0})=Ts_0\neq 0$, so
            $\hat{T}(1)=\one$. This means that $\one\in
            \hat{T}(\FAFA S)=\VLA\{Ts_0\}$. This is impossible because every element of $\VLA\{Ts_0\}$ has support contained
            in $[0,1]$.\qedhere
    \end{enumerate}
\end{proof}

It is natural to wonder whether a chain condition for the product
exists. Namely, if $F\subseteq \FAFA S$ is an uncountable family, is
it true that there exist $f,g \in F$ such that $fg\neq 0$? To answer
this question, it is useful to reformulate semiprimeness in terms of
what Scheffold \cite{scheffold1981} calls \fstar s: a vector lattice
algebra $A$ is called an \fstar\ when
\[
ab=0\text{ if and only if }|a|\wedge |b|=0\text{ for all }a,b \in A.
\]
The following is part of \cite[Théorème
9.3.1]{bigard_klaus_wolfenstein1977}.

\begin{lem}
    An \falg\ $A$ is an \fstar\ if and only if it is semiprime.
\end{lem}

\begin{rem}
    In particular, $\FAFA S$ is an \fstar. Therefore,
    families $F\subseteq \FAFA S$ such that $ab=0$ for all distinct $a,b \in F$ are the same as families of pairwise disjoint elements,
    and therefore must be countable.
\end{rem}

\chapter{A Structure Theorem for Banach \texorpdfstring{$f$-algebras}{f-algebras}}
\label{sec:structure}

\section{A structure theorem}

Before getting to free Banach \falg s, the main topic of this paper,
we need to introduce a new tool in the theory of Banach \falg s.
Recall that every vector lattice $X$, when equipped with the
identically $0$ product, becomes an \falg. We are denoting this \falg\
by $X_0$. In particular, when $X$ is a Banach lattice, $X_0$ is a
Banach \falg. Another basic example of Banach \falg\ is the space
$C(K)$ of continuous functions on a compact Hausdorff space $K$, with
the linear operations, order and product defined pointwise. Given
Banach \falg s $A$ and $B$ their direct sum $A\oplus_\infty B$,
equipped with the maximum norm and coordinatewise operations, is also a
Banach \falg.

The next theorem shows that every Banach \falg\ can be seen as a (not
necessarily closed) sublattice-algebra of
$X_0\oplus_\infty C(K)$, for appropriate $X$ and $K$.

\begin{thm}\label{thm:structure}
    For every Banach \falg\ $A$ there exist a Banach lattice $X$,
    a compact Hausdorff space $K$ and a contractive injective lattice-algebra
    homomorphism $R\colon A\to X_0\oplus_\infty C(K)$.
\end{thm}

The approach to the proof will be very similar to that of
\cref{thm:falgYudin}. First we consider the case that $A$ has trivial
annihilator. The following was already proved in \cite[Proposition
1.7]{martignon1980}. Since the proof uses some ideas that are worth
keeping in mind, we sketch it below.

\begin{lem}\label{lem:zeroan}
    For every Banach \falg\ $A$ with trivial
    annihilator there exists a contractive injective lattice-algebra
    homomorphism $A\to C(K)$, for a certain compact Hausdorff space
    $K$.
\end{lem}
\begin{proof}
    Since $A$ is a Banach lattice, $\orth A$ coincides with the center
    $\mathcal{Z}(A)$ of $A$ (see \cite[Theorem
    3.29]{abramovich_aliprantis2002}). The map $L\colon A\to \mathcal{Z}(A)$
    sending $a \in A$ to $L_a$, the operator of left multiplication by
    $a$, is an injective lattice-algebra homomorphism (see the discussion in
    \cref{sec:proof_falgYudin}). The center
    $\mathcal{Z}(A)$, being a Banach \falg\ with strong unit $I$, is
    lattice-algebra isometric to $C(K)$ for some compact Hausdorff
    space $K$ (see \cite[Proposition 1.4]{martignon1980}).
\end{proof}

The map in the previous lemma need not be an embedding, as the following
example shows.

\begin{example}\label{ex:not_embed}
    Consider $\el 1$ with the usual Banach lattice structure and
    pointwise product. This is a Banach \falg. Note that, given $x \in
    \el 1$, the norm of $L_x$ is
    \[
    \|L_x\|=\sup_{y \in B_{\el 1}} \|xy\|_1=\|x\|_\infty.
    \]
    Hence the map $L\colon \el 1\to \mathcal{Z}(\el 1)$ is not an
    embedding.
\end{example}

In view of the previous example, one could be led to think that the
embedding fails because we are taking a norm that is ``much bigger''
than a supremum norm. The following example shows that even for norms
defined as a supremum the map $L$ can fail to be an embedding.

\begin{example}
    Let $(a_n)\subseteq \R$ be a sequence such that $a_n\ge 1$ and
    $a_n\to \infty $. Consider the space of sequences
    \[
    A=\{\, (x_n)\in \R^{\N} : x_na_n\to 0 \, \} .
    \]
    This is a sublattice of $\R^{\N}$ that is closed under the
    coordinatewise product. Equipped with the norm
    \[
    \|(x_n)\|=\sup_n |x_n|a_n
    \]
    it becomes a Banach lattice. Since
    \[
    \|xy\|=\sup_n |x_n| |y_n|a_n\le \sup_n |x_n| |y_n| a_n^2\le \|x\|
    \|y\|,
    \]
    $A$ is a Banach \falg. The sequence $(e_n/a_n)$ is such that
    $\|e_n/a_n\|=1$. If $x \in A$ is in the unit ball (that is,
    $|x_n|\le 1/a_n$ for all $n \in \N$), then
    \[
    \|e_n/a_n x\|=\frac{1}{a_n} \|x_n e_n\|=|x_{n}|\le \frac{1}{a_n}
    \]
    which tends to $0$ as $n\to \infty $. Therefore $\|L_{e_n/a_n}\|$
    tends to zero while $\|e_n/a_n\|=1$.
\end{example}

\begin{proof}[Proof of \cref{thm:structure}]
    By embedding $A$ in its bidual (which is again a Banach
    \falg\ when equipped with the Arens product, see \cite{scheffold1991}), we may assume $A$ is order complete.
    Let $P\colon A\to A$ denote the band projection onto
    its annihilator $N=N(A)$,
    and let $P^{d}$ denote the band projection onto its disjoint
    complement. Both $P$ and $P^{d}$ are
    multiplicative (see the proof of \cref{thm:falgYudin}). In particular, both $N$ and $N^{d}$
    are Banach \falg s on their own, with the structure inherited from
    $A$.

    The map
    \[
    \begin{array}{cccc}
    T\colon& A & \longrightarrow & N\oplus_\infty N^{d} \\
            & f & \longmapsto & (Pf,P^{d}f) \\
    \end{array}
    \]
    defines an injective lattice-algebra homomorphism.
    Since both $P$ and $P^{d}$ are contractive, so is $T$.
    Note that $N$ is a Banach lattice with the zero product.
    On the other hand, $N^{d}$ is a Banach \falg\ with trivial
    annihilator. By \cref{lem:zeroan} there exists an injective and contractive lattice-algebra
    homomorphism $S\colon N^{d}\to C(K)$. The composition
    \[
        R=(I_N\oplus S)T\colon A\to N\oplus_\infty C(K)
    \]
    defines the desired injective and contractive lattice-algebra homomorphism.
\end{proof}

\begin{rem}
    \begin{enumerate}
        \item In \cref{thm:structure}, $R$ need not be an
            embedding, as \cref{ex:not_embed} shows.
        \item Passing to the bidual in the proof of \cref{thm:structure} is not
            superfluous, for the annihilator need not be a projection
            band when the Banach \falg\ is not order complete. For instance, consider in the Banach lattice
            $C[0,1]$ the element
            \[
                p(t)=(2t-1) \chi _{[1/2,1]},\quad t \in [0,1],
            \]
            and define, for $f,g \in C[0,1]$, the product $f\star
            g=fgp$, where juxtaposition denotes the usual pointwise
            product. It is clear that $C[0,1]$ with this product
            becomes a Banach \falg, with annihilator
            \[
                N=\{\, f \in C[0,1] : \supp(f)\subseteq [0,1/2] \, \}.
            \]
            This is certainly a band but it is
            not a projection band.
        \item \Cref{thm:structure} is not true for general \falg s, even
            if they are Archimedean and
            uniformly complete. For instance, $\R^{\N}$ with the
            pointwise order and product is an Archimedean and
            uniformly complete \falg\ with trivial annihilator, yet
            there is no positive map from $\R^{\N}$ to a
            space of continuous functions. Indeed, if $R\colon
            \R^{\N}\to C(K)$ were positive, where $K$ is some compact
            Hausdorff space, then
            the sequence $x=(n/\|Re_n\|_\infty )_n$ would satisfy
            \[
            Rx\ge \frac{n}{\|Re_n\|_\infty } Re_n
            \]
            so that $\|Rx\|_\infty \ge n$ for every $n \in \N$. This
            is a contradiction.
        \item However, \cref{thm:structure} does hold for normed
            \falg s: just embed them in their norm completion, which
            is a Banach \falg\ (see \cite[Proposition
            2.1]{jaber2020}).
        \item Given a Banach \falg\ $A$ and a map $R\colon A\to X_0\oplus_\infty
        C(K)$ as in \cref{thm:structure}, we shall denote $Ra=(a_N,a_C)$ for
        $a \in A$. Note that the maps $a\mapsto a_N$ and $a\mapsto a_C$ are
        also contractive lattice-algebra homomorphisms.
    \end{enumerate}
\end{rem}

The following theorem is a refinement of \cref{thm:falgYudin} for
normed \falg s.

\begin{thm}\label{thm:falgYudin_ball}
    Let $\Phi $ be an LLA expression. If $\Phi $
    vanishes on $[-1,1]$, then it also vanishes on the unit ball of
    every normed \falg.
\end{thm}
\begin{proof}
    Since $\Phi $ vanishes on a neighbourhood of $0$, it must vanish
    on $\R_0$ (\cref{lem:calc_zero_step2}). Since $\Phi $ vanishes on
    $\R_0$, it must also vanish on every
    vector lattice with the zero product (\cref{lem:calc_zero_step1}).
    Let $A$ be a Banach \falg, and let $a_1,\ldots ,a_n \in A$ be such
    that $\|a_i\|\le 1$.
    Let $R\colon
    A\to X_0\oplus_\infty C(K)$ be a map as in \cref{thm:structure}.
    Then
    \[
    \Phi (a_1,\ldots ,a_n)_N=\Phi ((a_1)_N,\ldots ,(a_n)_N)=0.
    \]
    Since $R$ is contractive, $\|(a_i)_C\|_\infty
    \le \|Ra_i\|\le \|a_i\|\le 1$, so $|(a_i)_C(t)|\le 1$ for every $t
    \in K$. By assumption:
    \[
    \Phi ((a_1)_C(t),\ldots ,(a_n)_C(t))=0\quad\text{for every }t \in
    K.
    \]
    This yields $0=\Phi ((a_1)_C,\ldots ,(a_n)_C)=\Phi (a_1,\ldots ,a_n)_C$.
    Thus $R(\Phi (a_1,\ldots ,a_n))=(0,0)$. Injectivity of $R$ implies
    $\Phi (a_1,\ldots ,a_n)=0$.
\end{proof}

\begin{rem}\label{rem:LLA vanishing}
    The previous result is meaningless for LL expressions: if an
    LL expression vanishes on $[-1,1]$, then by positive homogeneity it
    vanishes on $\R$. This is no longer
    true for LLA expressions since, in general, they are not
    positively homogeneous. For instance,
    \[
        \Phi [t]=(t_+^2 - t_+)_+
    \]
    is an LLA expression that vanishes on $[-1,1]$ but not on $\R$. By
    the previous theorem, the identity $(x_+^2-x_+)_+=0$ holds whenever $x$ is
    contained in the unit ball of a normed \falg.
\end{rem}

\section{Minimal representation of semiprime $f$-algebras}

In the case that a Banach \falg\ $A$ is semiprime,
\cref{thm:structure} yields a contractive lattice-algebra homomorphism
$A\to C(K)$, for some compact Hausdorff space $K$. We are going to
show that there exists a unique $K$ that is minimal in a certain sense,
and such that the map $A\to C(K)$ has good extension
properties. This minimal compact Hausdorff space will be the set of
functionals on $A$ that are lattice-algebra homomorphisms.

\begin{defn}
    Let $A$ be a Banach \falg. We will denote by $K_A\subseteq A^{*}$
    the set of functionals on $A$ that are lattice-algebra
    homomorphisms.
\end{defn}

\begin{lem}\label{lem:compact_sep_points}
    Let $A$ be a semiprime Banach \falg.
    The set $K_A$, equipped with the weak$^*$ topology, is a compact
    Hausdorff space that separates the points of $A$.
\end{lem}
\begin{proof}
    If $\phi \in K_A$ and $a \in A$, then
    \[
    |\phi (a)|^2=\phi (a^2)\le \|a\| \phi (|a|)=\|a\| |\phi (a)|,
    \]
    where we have used that $a^2\le \|a\| |a|$ (recall that this holds
    because left multiplication by $|a|$ is a central operator with
    norm no bigger than $\|a\|$). Hence $\|\phi \|\le 1$, and
    $K_A\subseteq B_{A^{*}}$. It is then immediate to check that $K_A$
    is closed in the weak$^*$ topology.

    To show that $K_A$ separates the points of $A$, let $a,b \in A$ be
    different points. By the arguments of \cref{thm:structure}, there
    exists a contractive and injective lattice-algebra homomorphism
    $R\colon A\to C(K)$, for some compact Hausdorff space $K$. Since
    $Ra \neq Rb$, there exists a point $t \in K$ at which $(Ra)(t)\neq
    (Rb)(t)$. Then $R^{*}\delta _t \in K_A$, where $\delta _t$ denotes
    the point evaluation at $t$, is such that $(R^{*}\delta _t)(a)\neq
    (R^{*}\delta _t)(b)$.
\end{proof}

\begin{prop}
    Let $A$ be a semiprime Banach \falg. Then the map $R\colon A\to
    C(K_A)$, that sends $a \in A$ to its evaluation functional, is a
    contractive and injective lattice-algebra homomorphism. Moreover, both
    the vector lattice and the algebra generated by $\{\, Ra : a \in A
    \, \} \cup \{\one\}$ in $C(K_A)$ are dense.
\end{prop}
\begin{proof}
    Injectivity follows from \cref{lem:compact_sep_points}. Density of both the vector lattice and the algebra generated by $\{\, Ra : a \in A
    \, \} \cup \{\one\}$ is a direct consequence of the Stone--Weierstrass theorem.
\end{proof}

The following example illustrates this construction and shows that
adding the constant one function is necessary in order to obtain density.

\begin{example}
    Consider $c_0$ with the usual order and product. The
    lattice-algebra homomorphisms on $c_0$ are the $\delta _i$, the
    projection onto the $i$ coordinate, and the zero functional.
    It is immediate to check that $K_{c_0}$ is homeomorphic to $\N\cup
    \{\infty \}$, and therefore that $C(K_{c_0})$ can be identified
    with $c$, the space of convergent sequences. The canonical map
    $R\colon c_0\to c$ is the usual inclusion.
\end{example}

\section[Extensions and injective objects]{Extension of lattice-algebra homomorphisms and injective objects}

In this section we show the remarkable fact that, if $A$ is a closed
sublattice-algebra of a Banach \falg\ $B$, then every
lattice-algebra homomorphism $\phi \colon A\to \R$ extends to a
lattice-algebra homomorphism on $B$. Notice that, for such a
real-valued lattice-algebra homomorphism $\phi \colon A\to \R$, $\phi
(x)=0$ for every $x \in N(A)$. Hence $\phi |_{N(A)}=0$ and, to truly
understand $\phi $ and its extensions, we need only understand its
restriction to the semiprime Banach \falg\ $N(A)^{d}$. For this
reason, our focus is in the case that $A$ is semiprime.

So let $A$ be a semiprime Banach \falg, and let $R\colon A\to C(K_A)$
be the map that sends each element of $A$ to its evaluation functional.
Since the non-zero lattice-algebra homomorphisms on $C(K_A)$ are precisely the
point evaluations, it follows from the construction that each
non-zero lattice-algebra homomorphism $\phi \colon A\to \R$ extends uniquely to
a lattice-algebra homomorphism $\tilde \phi \colon C(K_A)\to \R$, in
the sense that $\phi =\tilde \phi R$. The zero homomorphism extends in precisely two ways: as the zero homomorphism, and as the point evaluation at $0\in K_A$. The following shows that this extension property characterizes $K_A$.

\begin{prop}\label{prop:char_KA1}
    Let $A$ be a semiprime Banach \falg. Let $R\colon A\to C(K)$ be a
    contractive and injective lattice-algebra homomorphism, for some
    compact Hausdorff space $K$. Suppose that every lattice-algebra
    homomorphism $\phi \colon A\to \R$ extends uniquely to a non-zero
    lattice-algebra homomorphism $\tilde \phi \colon C(K)\to \R$, in
    the sense that $\phi =\tilde \phi R$. Then $K$ is homeomorphic to
    $K_A$.
\end{prop}
\begin{proof}
    Consider the map $K\to K_A$ that sends each $t \in K$ to
    $R^{*}\delta _t \in K_A$. This map is certainly continuous. It is
    also bijective, because for every $\phi \in K_A$, the unique
    non-zero extension $\tilde\phi \colon C(K)\to \R$ is of the form
    $\tilde \phi =\delta _t$ for some $t \in K$, and $t$ is precisely
    the unique preimage of $\phi $ through this map. Since $K$ and
    $K_A$ are compact, this map is a homeomorphism.
\end{proof}

\begin{prop}\label{prop:char_KA2}
    Let $A$ be a semiprime Banach \falg. Let $R\colon A\to C(K)$ be a
    contractive and injective lattice-algebra homomorphism, for some
    compact Hausdorff space $K$. Suppose that every lattice-algebra
    homomorphism $\phi \colon A\to \R$ extends to a non-zero
    lattice-algebra homomorphism $\tilde \phi \colon C(K)\to \R$, in
    the sense that $\phi =\tilde \phi R$. Suppose also that the vector
    lattice algebra generated by $\{\, Ra : a \in A \, \} \cup
    \{\one\}$ is dense in $C(K)$. Then $K$ is homeomorphic to $K_A$.
\end{prop}
\begin{proof}
    By \cref{prop:char_KA1}, we need only show that the (non-zero) extension
    of every lattice-algebra homomorphism $\phi $ is unique. Suppose
    $\psi $ and $\psi '$ are two non-zero lattice-algebra
    homomorphisms that satisfy $\psi R=\psi 'R$. Since every non-zero
    lattice-algebra homomorphism on $C(K)$ is a point evaluation,
    $\psi (\one)=1=\psi '(\one)$. Thus $\psi $ and $\psi '$ coincide
    on the vector lattice algebra generated by $\{\, Ra : a \in A \,
    \} \cup \{\one\}$. By density, $\psi =\psi '$.
\end{proof}

It is remarkable that the center of a Banach \falg\ provides
an instance of a $C(K)$ space to which every functional that is a
lattice-algebra homomorphism extends.

\begin{prop}\label{prop:KA_center}
    Let $A$ be a Banach \falg, and let $L\colon A\to
    \mathcal{Z}(A)$ be the map that sends every element of $A$ to its left
    multiplication operator. For every lattice-algebra homomorphism
    $\phi \colon A\to \R$ there exists a
    lattice-algebra homomorphism $\tilde \phi \colon \mathcal{Z}(A)\to
    \R$ satisfying $\phi =\tilde \phi L$. Moreover, if $A$ is
    semiprime and the
    closed sublattice-algebra generated by $\{\, La: a \in A  \, \}
    $ in $\mathcal{Z}(A)$ does not contain $I_A$, then $C(K_A)$ can
    be identified with the closed unital sublattice-algebra of
    $\mathcal{Z}(A)$ generated by $\{\, La : a \in A \, \} $.
\end{prop}
\begin{proof}
    The zero homomorphism clearly satisfies the statement, so assume that $\phi \neq 0$.
    Then there exists $x_0 \in A_+$ such that $\phi (x_0)=1$. Define
    the map
    \[
    \begin{array}{cccc}
    \tilde \phi \colon& \mathcal{Z}(A) & \longrightarrow & \R \\
            & T & \longmapsto & \phi (Tx_0) \\
    \end{array}.
    \]
    First note that, for every $a \in A$,
    \[
    \tilde \phi (L_a)=\phi (L_a(x_0))=\phi (ax_0)=\phi (a) \phi
    (x_0)=\phi (a),
    \]
    so $\phi =\tilde \phi L$. It is immediate that $\tilde \phi $ is
    linear. It is also a lattice homomorphism, because
    \[
    \tilde \phi (|T|)=\phi (|T|x_0)=\phi (|Tx_0|)=|\phi (Tx_0)|=|\tilde \phi (T)|,
    \]
    where in the second equality we are using that the modulus of a
    central operator is computed pointwise, and in the third that
    $\phi $ is a lattice homomorphism. We also have $\tilde \phi
    (I)=1$. Since $\mathcal{Z}(A)$ is lattice-algebra isometric to a
    space of continuous functions on a compact Hausdorff space, every
    functional that is a lattice homomorphism and sends the identity
    to $1$ is a point evaluation, and therefore an algebra
    homomorphism. Hence $\tilde \phi $ is the desired extension.

    Suppose that $\mathcal{A}$, the closed sublattice-algebra
    generated by $\{\, La: a \in A  \, \} $ in $\mathcal{Z}(A)$,
    does not contain $I_A$. Then there exists a point evaluation on
    $\mathcal{Z}(A)$ that vanishes on $\mathcal{A}$. This point
    evaluation is a non-zero extension of the zero functional on
    $\mathcal{A}$. Hence, by \cref{prop:char_KA2}, the closed
    sublattice-algebra generated by $\{\, La : a \in A \, \} \cup
    \{I\}$ in $\mathcal{Z}(A)$, which is again a space of continuous
    functions, can be identified with $C(K_A)$.
\end{proof}

These properties can now be used to prove the main extension result.

\begin{thm}\label{thm:extension}
    Let $B$ be a Banach \falg, and let $A\subseteq B$ be a closed
    sublattice-algebra. Then
    for every lattice-algebra homomorphism $\phi \colon A\to \R$ there
    exists a lattice-algebra homomorphism $\psi \colon B\to \R$ such that
    $\psi |_A=\phi $.
\end{thm}
\begin{proof}
    As usual, denote by $L_b\colon B\to B$ the multiplication operator
    defined by $b
    \in B$. Note that, for every $a \in A$, the operator $L_a$ leaves
    $A$ invariant and its restriction to $A$ is precisely the
    multiplication operator defined by $a$ in $A$. Let $\mathcal{A}$
    be the unital sublattice-algebra generated by $\{\,
    L_a|_A : a \in A \, \}$ inside $\mathcal{Z}(A)$.

    Let $\mathcal{A}'$ be the unital sublattice-algebra generated by
    $\{\, L_a : a \in A \, \} $ inside $\mathcal{Z}(B)$. We claim that
    every element of $\mathcal{A}'$ leaves $A$ invariant, and that the
    elements of $\mathcal{A}$ are precisely the restriction of
    elements of $\mathcal{A}'$ to $A$. Indeed, if $T,S \in
    \mathcal{Z}(B)$ leave $A$ invariant, then:
    \begin{enumerate}
        \item $\alpha T+\beta S$ leaves $A$ invariant, for every
            $\alpha ,\beta  \in \R$, and $(\alpha T+\beta S)|_A=\alpha
            T|_A+\beta S|_A$.
        \item $TS$ leaves $A$ invariant, and $(TS)|_A=T|_AS|_A$.
        \item $|T|$ leaves $A$ invariant, and $(|T|) |_A=|(T|_A)|$.
            Indeed, if $x \in A_+\subseteq B_+$, then since the
            lattice operations on central operators are computed
            pointwise both on $A_+$ and $B_+$:
            \[
                |T| |_A(x)=|T|(x)=|Tx|=|T|_Ax|=|T|_A|(x).
            \]
    \end{enumerate}

    In other words, we have shown that the map
    \[
    \begin{array}{cccc}
    \mathcal{R}\colon& \mathcal{A}' & \longrightarrow & \mathcal{A} \\
            & T & \longmapsto & T|_A \\
    \end{array},
    \]
    is a lattice-algebra homomorphism.
    Moreover, if we equip $\mathcal{A}$ with the norm of
    $\mathcal{Z}(A)$ (i.e., with the uniform norm $\|{\cdot }\|_{I_A}$
    induced by the identity $I_A$ on $A$), and $\mathcal{A}'$ with the
    norm of $\mathcal{Z}(B)$ (i.e., with the uniform norm $\|{\cdot
    }\|_{I_B}$ induced by the identity $I_B$ on $B$), then
    $\mathcal{R}$ is contractive. Indeed, if $|T|\le \|T\|_{I_B} I_B$, then
    taking restrictions to $A$, we have 
    \[
    |T|_A|=|T| |_A\le \|T\|_{I_B} I_B|_A=\|T\|_{I_B}I_A.
    \]
    Hence $\|T|_A\|_{I_A}\le \|T\|_{I_B}$, and
    $\mathcal{R}$ extends naturally to the closures $\mathcal{R}\colon
    \overline{\mathcal{A}'}\to \overline{\mathcal{A}}$.

    Let $\phi \colon A\to \R$ be a lattice-algebra homomorphism. From
    \cref{prop:KA_center} it follows that there exists a
    lattice-algebra homomorphism $\psi \colon
    \overline{\mathcal{A}}\to \R$ such that $\phi =\psi L$. Then $\psi  \circ \mathcal{R}\colon \overline{\mathcal{A}'}\to \R$ is a
    lattice-algebra homomorphism satisfying
    \[
        (\psi \circ \mathcal{R})(L_a)=\psi (L_a|_A)=\phi
        (a).
    \]
    Since $\overline{\mathcal{A}'}$ is a closed unital subalgebra of
    $\mathcal{Z}(B)$, which can be identified with a space of
    continuous functions, $\psi \circ \mathcal{R}$ further extends to
    a lattice-algebra homomorphism $\tilde \psi $ on $\mathcal{Z}(B)$.
    The composition $\tilde \phi =\tilde \psi L\colon B\to \R$ is the
    claimed extension.
\end{proof}

The extensions in \cref{thm:extension} do not, in general, preserve
the norm.

\begin{example}\label{ex:extension_not_preserves_norm}
    Consider $\el 1^2$ with pointwise order and product. This is a
    Banach \falg. Let $A=\barespn\{(1,1)\}$. It is clear that $A$ is a
    sublattice-algebra of $\el 1^2$. The map $\phi \colon A\to \R$
    defined by $\phi (\lambda,\lambda )=\lambda $, $\lambda \in \R$,
    is a lattice-algebra homomorphism. If $(\lambda ,\lambda )\in A$,
    and $\|(\lambda ,\lambda )\|_1\le 1$, then $|\lambda |\le 1/2$ and
    therefore $\|\phi \|\le 1/2$.

    The non-zero lattice-algebra homomorphisms on $\el 1^2$ are the
    projections onto the coordinates. Both of these functionals are
    extensions of $\phi $, yet they have norm one.
\end{example}

It is pertinent now to introduce, in the category of Banach \falg s,
the analogue of the notion of $\lambda $-injective Banach lattice
mentioned in \cref{sec:background_projective}.

\begin{defn}
    Let $\lambda \ge 1$.
    A Banach \falg\ $I$ is \emph{$\lambda $-injective} if, for every Banach
    \falg\ $A$, every closed sublattice-algebra $B$ of $A$, and every
    lattice-algebra homomorphism $T \colon B\to I$, there exists an
    extension $\tilde T \colon A\to I$ that is also a
    lattice-algebra homomorphism and satisfies $\|\tilde T\|\le
    \lambda \|T\|$.
\end{defn}

\Cref{ex:extension_not_preserves_norm} shows that $\R$ is not
$1$-injective in the category of Banach \falg s. This example can be
generalized to show the following.

\begin{thm}\label{thm:no_injectives}
    For $\lambda \ge 1$, there are no non-zero $\lambda $-injective
    Banach \falg s.
\end{thm}

For the proof of this theorem we will need the following
elementary lemma.

\begin{lem}\label{lem:lat_hom_sum}
    Let $\Gamma $ be a set and let $\{X_\gamma \}_{\gamma \in \Gamma
    }$ be a family of Banach lattices. Let $X$ be
    their $\el 1$-sum, and let $p_\gamma \colon X\to X_\gamma $ denote the
    projection onto the $\gamma  \in \Gamma $ coordinate. Let $\phi \colon
     X\to \R$ be a lattice homomorphism. Then there exists a $\delta
     \in \Gamma $ and a lattice homomorphism $\phi '\colon X_\delta
     \to \R$ such that $\phi =\phi ' p_\delta $.
\end{lem}
\begin{proof}
    The result is obvious for $\phi =0$, so assume $\phi \neq 0$. Let
    $(x_\gamma )\in X_+$ be such that $\phi ( (x_\gamma) )>0$. For $x\in X_\gamma$, denote
    by $x e_\gamma $ the element of $X$ that takes the value $x$ at
    $\gamma \in \Gamma $ and zero everywhere else. Since $(x_\gamma
    )=\sum_{\gamma \in \Gamma}^{}x_\gamma e_\gamma $ and $0\le
    x_\gamma e_\gamma $ for every $\gamma \in \Gamma $, there must
    exist some $\delta \in \Gamma $ for which $\phi (x_\delta e_\delta
    )>0$. Now if $y_\gamma \in (X_\gamma )_+$ with $\gamma \neq \delta $, then $(x_\delta
    e_\delta )\wedge (y_\gamma e_\gamma )=0$, and since $\phi $ is a
    lattice homomorphism, $\phi (x_\delta e_\delta )\wedge \phi
    (y_\gamma e_\gamma )=0$. Thus $\phi (y_\gamma e_\gamma )=0$. Now
    for a general $(y_\gamma )\in X$:
    \[
    \phi (y_\gamma )=\phi \bigg( \sum_{\gamma \in \Gamma }^{} y_\gamma
    e_\gamma \bigg) =\sum_{\gamma \in \Gamma }^{} \phi (y_\gamma
    e_\gamma )=\phi (y_\delta e_\delta ).
    \]
    If we define $\phi '\colon X_\delta \to \R$ by $\phi '(x)=\phi (x
    e_\delta )$, then it is immediate to check that $\phi '$ is a
    lattice homomorphism satisfying $\phi =\phi 'p_\delta $.
\end{proof}

\begin{lem}\label{lem:nonzero_hom}
    Let $A$ be a Banach \falg. Suppose that the product on $A$ is not
    identically zero. Then there exists a non-zero lattice-algebra
    homomorphism $\phi \colon A\to \R$.
\end{lem}
\begin{proof}
    By \cref{thm:structure} there exist a non-empty compact Hausdorff
    space $K$, a Banach lattice $X$, and a contractive injective
    lattice-algebra homomorphism $R\colon A\to X_0\oplus_\infty C(K)$.
    Projecting onto the second coordinate, we obtain a lattice-algebra
    homomorphism $R_2\colon A\to C(K)$. Since $K$ is non-empty (because the
    product on $A$ is not identically zero), there exists a non-zero
    lattice-algebra homomorphism $\delta _t\colon C(K)\to \R$ (a point
    evaluation at some $t \in K$). Hence the composition $\delta
    _tR_2\colon A\to \R$ is a non-zero lattice-algebra homomorphism.
\end{proof}

\begin{proof}[Proof of \cref{thm:no_injectives}.]
    Suppose that $I$ is a non-zero $\lambda $-injective Banach \falg, for
    some $\lambda \ge 1$. First note that the product in $I$ cannot be
    identically zero, for if it were, then $I$ would be $\lambda
    $-injective for lattice homomorphisms (if $X$ is a Banach lattice, and
    $T\colon X\to I$ is a lattice homomorphism, then $T\colon X_0\to
    I_0$ is a lattice-algebra homomorphism). But, as explained in
    \cref{sec:background_projective}, there are no non-zero $\lambda
    $-injective objects in the category whose objects are Banach
    lattices and whose morphisms are lattice homomorphisms.

    Choose $n \in \N$ with $n> \lambda $. Let
    $B$ be the $\el 1$-sum of $n$ copies of $I$, and let
    \[
    A=\{\, (x,\ldots ,x) : x \in I \, \}.
    \]
    It is immediate to check that $A$ is a closed sublattice-algebra
    of $B$. Let $\phi \colon A\to I$ be defined by $\phi (x,\ldots
    ,x)=x$ for every $x \in I$. This is a
    lattice-algebra homomorphism with $\|\phi \|=1/n$; indeed, if
    $\|(x,\ldots ,x)\|_1=1$, then $\|x\|=1/n$.

    Since $I$ is $\lambda $-injective, this map extends to a
    lattice-algebra homomorphism $\tilde \phi \colon B\to I$ with
    $\|\tilde \phi \|\le \lambda /n$. Since the product in $I$ is not
    identically zero, by \cref{lem:nonzero_hom}
    there exists a non-zero lattice-algebra homomorphism $\psi \colon
    I\to \R$. Then the composition $\psi \tilde \phi \colon B\to \R$ is a
    lattice homomorphism which, by \cref{lem:lat_hom_sum}, is of the
    form $\psi \tilde \phi =\phi 'p_i$ for some lattice homomorphism
    $\phi '\colon I\to \R$, where $p_i\colon B\to I$ denotes the
    projection onto the $i$ coordinate. Now note that, for every $x
    \in I$,
    \[
    \phi '(x)=(\psi \tilde \phi )(x,\ldots ,x)=(\psi \phi )(x,\ldots
    ,x)=\psi (x).
    \]
    Thus $\psi \tilde \phi = \psi p_i$. It is now immediate to check
    that $\|\psi p_i\|=\|\psi \|$, so that
    \[
    \|\psi \|=\|\psi \tilde \phi \|\le \|\psi \| \|\tilde \phi \|\le
    \|\psi \| \frac{\lambda }{n}.
    \]
    Since $\|\psi \|\neq 0$, $n\le \lambda $, which contradicts the
    initial choice of $n \in \N$.
\end{proof}

\begin{rem}
    According to \cref{thm:extension}, in Banach \falg s we can always
    extend real-valued lattice-algebra homomorphisms but, according to
    \cref{thm:no_injectives}, we can never do so with a uniform
    control over the norm.
\end{rem}

\chapter{The Free Banach \texorpdfstring{$f$-algebra}{f-algebra}
Generated by a Banach Space}
\label{sec:FBFA}

The notion of a free Banach lattice algebra has been considered
before \cite{wickstead2017_open, de_jeu2021}. Through abstract
arguments from universal algebras, it is possible to show that such
an object exists. Yet no explicit construction is known, and it is
not clear at all whether such an explicit description is even
possible.

The goal of this chapter is to study free objects in the much better
behaved category of Banach \falg s and to get as close as we can to an
explicit description of these spaces. We shall see that, even though a
large theory of free Banach lattices is known (see
\cref{sec:background_FBL} and the references therein), free Banach
\falg s are still challenging to understand. New ideas are required,
since the introduction of a product invalidates arguments that rely
on positive homogeneity. In fact, a first subtlety due to the lack of
positive homogeneity is already present in the definition of free
Banach \falg.

\begin{defn}
    Let $E$ be a Banach space. The \emph{free Banach \falg\ generated
    by $E$} is a Banach \falg\ $\FBFA E$ together with a linear isometric
    embedding $\eta _E \colon E\to \FBFA E$ such that, for every Banach
    \falg\ $A$ and every contractive operator $T\colon E\to A$, there exists a
    unique lattice-algebra homomorphism $\hat{T}\colon \FBFA E\to A$
    satisfying $\hat{T}\circ  \eta _E=T$ and $\|\hat{T}\|=\|T\|$.
\end{defn}

\begin{rem}
    The universal property defining $\FBFA E$ only extends
    contractive operators. The reason for this is that non-contractive
    maps cannot, in general, be extended to bounded algebra
    homomorphisms. For example, if the map
    $T\colon \R\to \R$, $Tx=2x$, extended to a lattice-algebra
    homomorphism $\hat{T}\colon \FBFA \R\to \R$, then
    \[
    \hat{T}(\eta _\R(1)^{n})=T(1)^{n}=2^{n}
    \]
    while $\|\eta _\R(1)^{n}\|\le \|\eta _\R(1)\|^{n}\le 1$. This
    would imply that $\hat{T}$ is unbounded. Another way to see this
    obstruction is that, while scaling a lattice homomorphism by a
    positive number gives another lattice homomorphism, scaling an
    algebra homomorphism does not give an algebra homomorphism (unless
    the scalar is $1$).
\end{rem}

\begin{rem}
    Let $S$ be a set. We shall denote $\FBFA{\el 1(S)}$ by $\FBFAs S$
    and call it the \emph{free Banach \falg\ generated by $S$}. The
    reason for this is that for every Banach \falg\ $A$ and every map
    $T\colon S\to A$ satisfying $\|T\|=\sup_{s \in S}\|Ts\|\le 1$
    there exists a unique lattice-algebra homomorphism $\hat{T}\colon
    \FBFAs S\to A$ with $\|\hat{T}\|=\|T\|$ satisfying $\hat{T}\eta
    _{e_s}=Ts$.
\end{rem}

The existence of the free Banach \falg\ generated by a Banach space follows from that of $\FAFA S$
using abstract general arguments that are already known. The details
of this construction are carried out in
\cref{sec:abstract_construction}. But, as we have anticipated, these arguments do not provide an explicit description of the space. In
\cref{sec:explicit_description} we give a more explicit and useful description of free
Banach \falg s using \cref{thm:falgYudin_ball}. In
\cref{sec:fin_dim} a
characterization, up to isomorphism, of the free Banach \falg\
generated by a finite-dimensional Banach space is derived. Finally,
\cref{sec:free_norm} is devoted to the study of the free norm; it will
be particularly useful when we study the properties of free Banach
\falg s.

Although the following objects are not central, it is
convenient to have them at hand, since
they appear as intermediate steps in the construction.

\begin{defn}
    Let $E$ be a vector space. The \emph{free Archimedean \falg\
    generated by $E$} is an Archimedean \falg\ $\FAFAv E$ together with a
    linear map $\delta_E \colon E\to \FAFAv E$ such that, for every
    Archimedean
    \falg\ $A$ and every linear map $T\colon E\to A$, there exists a
    unique lattice-algebra homomorphism $\hat{T}\colon \FAFAv E\to A$
    satisfying $\hat{T}\circ \delta _E=T$.
\end{defn}

\begin{defn}
    Let $E$ be a normed space. The \emph{free normed \falg\ generated
    by $E$} is a normed \falg\ $\FNFA E$ together with a linear isometric
    embedding $\eta _E \colon E\to \FNFA E$ such that, for every
    normed
    \falg\ $A$ and every contractive operator $T\colon E\to A$, there exists a
    unique lattice-algebra homomorphism $\hat{T}\colon \FNFA E\to A$
    satisfying $\hat{T}\circ \eta _E=T$ and $\|\hat{T}\|=\|T\|$.
\end{defn}

The reason to have that much notation is the following: through the
universal properties, it is not difficult to check that $\FAFAv E$ is
nothing more than $\FAFA S$, where $S$ is a Hamel basis of $E$. We will
equip $\FAFAv E$ with an appropriate seminorm, and quotient out by its kernel
to obtain a norm. This will yield $\FNFA E$. It is easy to check,
again using universal properties, that
the completion of $\FNFA E$ is precisely $\FBFA E$. With these
notations we are aware, at each step of the process, which universal
property each object has. Furthermore, we will also
exhibit an explicit description of these free objects.

\section{Abstract construction}\label{sec:abstract_construction}

We shall follow the previous outline to show that $\FBFA E$ exists. In
doing so, we follow closely the expositions of V.\ Troitsky
\cite{troitsky2019} and M.\ de Jeu \cite{de_jeu2021}.
The construction of $\FAFAv E$ that we use is essentially the
same, with the appropriate adaptations, as that of \cite{troitsky2019} for
the free vector lattice over a vector space.

Let $E$ be a vector space, and let $E^{\#}$ be its algebraic dual. The
space of functions $\R^{E^{\#}}$, equipped with pointwise linear
operations, order and product, is an Archimedean \falg. For
every $x \in E$, define $\delta _x\colon E^{\#}\to \R$ by $\delta
_x(\omega )=\omega (x)$ for all $\omega \in E^{\#}$. Denote $\FAFAv
E=\VLA \{\, \delta _x : x \in E \, \} \subseteq \R^{E^{\#}}$, and
define the map $\delta \colon E\to \FAFAv E$ by $\delta (x)=\delta
_x$.

\begin{prop}\label{prop:fafav}
    The free Archimedean \falg\ over a vector space $E$ is $\FAFAv E$
    together with the map $\delta $.
\end{prop}
\begin{proof}
    Being a sublattice-algebra of the Archimedean \falg\
    $\R^{E^{\#}}$, $\FAFAv E$ is again an Archimedean \falg.
    Moreover, the map $\delta $ is linear. Indeed, if $x,y
    \in E$, $\lambda ,\mu \in \R$ and $\omega \in E^{\#}$, then
    \[
    \delta _{\lambda x+\mu y}(\omega )=\omega (\lambda x+\mu
    y)=\lambda \omega (x)+\mu \omega (y)=(\lambda \delta _x+\mu \delta
    _y)(\omega ).
    \]

    Let $A$ be an Archimedean \falg\ and let $T\colon E\to A$ be a
    linear map. For $f \in \FAFAv E$, say $f=\Phi (\delta
    _{x_1},\ldots ,\delta _{x_n})$ for some LLA expression $\Phi $ and
    $x_1,\ldots ,x_n \in E$, define $\hat{T}(f)=\Phi (Tx_1,\ldots,
    Tx_n)$. As in the proof of \cref{thm:FAFA}, the only non-trivial
    part is to check that $\hat{T}$ is well-defined. For this it
    suffices to check that, if
    $\Phi (\delta _{x_1},\ldots ,\delta _{x_n})=0$ in $\FAFAv E$, then
    $\Phi (Tx_1,\ldots ,Tx_n)=0$ in $A$.

    Let $\{z_1,\ldots ,z_m\}$ be a basis of the span of $\{x_1,\ldots
    ,x_n\}$. Write, for $i=1,\ldots ,n$,
    \[
    x_i=\sum_{j=1}^{m} \alpha _{ij} z_j\quad\text{for certain }\alpha
    _{ij}\in \R.
    \]
    Define the LLA expression
    \[
        \Psi [t_1,\ldots ,t_m]=\Phi \bigg[\sum_{j=1}^{m}\alpha
        _{1j}t_j,\ldots ,\sum_{j=1}^{m}\alpha _{nj}t_j\bigg].
    \]
    By the linearity of $\delta $,
    \[
        0=\Phi (\delta _{x_1},\ldots ,\delta _{x_n})=\Psi (\delta
        _{z_1},\ldots ,\delta _{z_m}),
    \]
    which implies
    \[
    \Psi (\omega (z_1),\ldots ,\omega (z_m))=0
    \]
    for every $\omega \in E^{\#}$. Since $z_1,\ldots ,z_m$ are
    linearly independent, $(\omega (z_1),\ldots ,\omega (z_m))$ ranges
    through $\R^{m}$ as $\omega$ ranges through $E^{\#}$. By \cref{thm:falgYudin},
    $\Psi $ vanishes on $A$. Using the linearity of $T$ it follows
    that
    \[
    \Phi (Tx_1,\ldots ,Tx_n)=\Phi \bigg(\sum_{j=1}^{m}\alpha
    _{1j}Tz_j,\ldots ,\sum_{j=1}^{m}\alpha _{nj}Tz_j\bigg)=\Psi
    (Tz_1,\ldots ,Tz_m)=0.\qedhere
    \]
\end{proof}

\begin{rem}\label{rem:not_injective}
    Note that, when $E$ is a normed space, the previous construction can
    be carried out in $\R^{E^{*}}$ instead of $\R^{E^{\#}}$. In
    contrast with the free vector lattice over a vector
    space, however, we cannot go one step further and
    restrict the elements of $\FAFAv E$ to $B_{E^{*}}$ so as to obtain
    a representation of $\FAFAv E$ in $C(B_{E^{*}},w^{*})$. The reason
    is that LL expressions are positively homogeneous, while
    LLA expressions need not be (for instance, in $\FAFAv \R$ the
    non-zero element $((\delta _1)_+^2-(\delta _1)_+)_+$ vanishes on
    $B_{\R^{*}}=[-1,1]$; see \cref{rem:LLA vanishing}). Despite
    this, it is remarkable that in
    \cref{sec:explicit_description,sec:representation} we will still
    be able to represent $\FNFA E$ and, in some cases, $\FBFA E$
    inside $C(B_{E^{*}})$.
\end{rem}

Next we proceed to the abstract construction of $\FNFA E$ for a normed
vector space $E$. In $\FAFAv E$ (regarded as the sublattice-algebra of
$\R^{E^{*}}$ generated by $\{\, \delta_x : x \in E\,\}$), let $\mathcal{N}$ be the set
of lattice seminorms $\nu $ that are submultiplicative and satisfy $\nu (\delta _x)\le \|x\|$ for all $x \in E$. This set is not empty,
for if $x^{*} \in B_{E^{*}}$,
then $\nu _{x^{*}}(f)=|f(x^{*})|$ defines an element of
$\mathcal{N}$.

Let $\rho (f)=\sup_{\nu \in \mathcal{N}}\nu (f)$ for every $f \in \FAFAv E$. We claim
that $\rho $ defines a submultiplicative lattice seminorm on $\FAFAv E$.
First, we need to show that the quantity $\rho (f)$ is actually finite.
If $f=\Phi (\delta _{x_1},\ldots ,\delta _{x_n})$, for a certain
LLA expression $\Phi$ and $x_1,\ldots ,x_n \in E$, it is not difficult
to check, by induction on the complexity of $\Phi $, that there exists
a polynomial $p(t_1,\ldots ,t_n)\in \R_+[t_1,\ldots ,t_n]$ with
$p(0,\ldots ,0)=0$ such that
\begin{equation}\label{eq:polynomial_bound}
|f|\le p(|\delta _{x_1}|,\ldots ,|\delta _{x_n}|).
\end{equation}
Then for every $\nu \in \mathcal{N}$:
\[
\nu (|f|)\le \nu (p(|\delta _{x_1}|,\ldots ,|\delta _{x_n}|))\le p(
\nu (|\delta _{x_1}|),\ldots ,\nu (|\delta _{x_n}|))\le p(\|x_1\|,\ldots
,\|x_n\|),
\]
where in the second inequality we have used the triangle inequality
and that $\nu $ is
submultiplicative, and in both the second and the third inequalities
that the coefficients of $p$ are positive. So $\rho (f)=\sup_{\nu \in
\mathcal{N}}\nu (f) \le p(\|x_1\|,\ldots
,\|x_n\|)$ is indeed finite. That $\rho (f)$ is a submultiplicative
lattice seminorm follows from the definition.

The kernel of this seminorm
\[
\ker \rho =\{\, f \in \FAFAv E : \rho (f)=0 \, \}
\]
is an order and algebraic ideal. An elementary argument shows that the
quotient
$\FAFAv E/\ker \rho $ is again an \falg\ (see \cite[Proposition
3.2]{boulabiar2002}). Denote this quotient by $\FNFA E$, denote by
$q_E\colon \FAFAv E\to \FNFA E$ the quotient map, and set $\eta
_E=q_E\delta _E$. Then $\rho $ induces a norm in $\FNFA E$ that
makes it into a normed \falg. Moreover, the map $\eta _E$ is
contractive with respect to this norm (in fact, it is an isometry; see
below). Denote the norm simply by $\|{\cdot }\|$ and let $\FBFA E$ be the completion
of $\FNFA E$ with respect to this norm.
Abusing the notation, we shall denote also by $\eta _E$ the
composition of this map with the canonical embedding $\FNFA E\to \FBFA
E$.

\begin{prop}\label{prop:free_banach_falg}
    The free normed \falg\ generated by a normed space $E$ is $\FNFA E$
    together with the map $\eta _E$. The free Banach \falg\ generated
    by a Banach space $E$ is $\FBFA E$ together with the map $\eta
    _E$.
\end{prop}
\begin{proof}
    We need to check that $\eta _E$ is a linear isometric embedding
    and that $\FNFA E$, together with $\eta _E$, satisfies
    the universal property defining the free normed \falg\ generated by $E$.
    Let $(A,\|{\cdot }\|_A)$ be a normed \falg, and let $T\colon E\to
    A$ be a contractive operator. Let $\tilde T\colon \FAFAv E\to
    A$ be the unique lattice-algebra homomorphism satisfying $\tilde
    T\delta _E=T$. It is direct to check that, since $\|T\|\le 1$,
    $f\mapsto \|\tilde Tf\|_A/\|T\|$ defines an element of
    $\mathcal{N}$.
    This implies that $\|\tilde Tf\|_A\le \|T\|\rho (f)$ for all $f \in \FAFAv
    E$. Hence there exists a unique lattice-algebra homomorphism
    $\hat{T}\colon \FAFAv E/\ker \rho \to A$ satisfying
    $\hat{T}q_E=\tilde T$ and $\|\hat{T}\|\le \|T\|$. This implies $\hat{T}\eta
    _E=\hat{T}q_E\delta _E=\tilde T \delta _E=T$ and, since $\eta _E$
    is contractive, also $\|\hat{T}\|=\|T\|$.

    It only remains to show that $\eta _E$ is an isometric embedding.
    Let $T\colon E\to \el \infty (B_{E^{*}})$ be the standard
    isometric embedding $Tx=(x^{*}(x))_{x^{*}\in B_{E^{*}}}$. Since
    $\el \infty (B_{E^{*}})$, with coordinatewise order and product,
    is a Banach \falg, there exists a
    lattice-algebra homomorphism $\hat{T}\colon \FNFA E\to \el
    \infty (B_{E^{*}})$ satisfying $\|\hat{T}\|=1$ and $\hat{T}\eta
    _E=T$. Using that $\eta _E$ is contractive, it follows that
    \[
    \|x\|=\|Tx\|\le \|\hat{T}\|\|\eta _Ex\|=\|\eta _Ex\|\le \|x\|
    \]
    holds for all $x \in E$. Therefore $\|\eta _Ex\|=\|x\|$.

    It is direct that the completion of $\FNFA E$, namely
    $\FBFA E$, together with $\eta _E$, is the free Banach \falg\ generated by $E$.
\end{proof}

\begin{rem}\label{rem:FBLinFBFA}
    Let $Z$ be the closure of $\Lat \{\, \eta _x : x \in E \, \} $ in
    $\FBFA E$. Let $\phi _E\colon E\to Z$ be $\eta
    _E$ with the codomain restricted to $Z$. Then $Z$, together with
    $\phi _E$, is the free Banach lattice $\FBL E$. Indeed, let $X$ be a Banach
    lattice, and let $T\colon E\to X$ be an operator. Equip $X$ with
    the zero product, so that it becomes a Banach \falg. Then $S=T/\|T\|$
    extends to a lattice-algebra homomorphism $\hat{S}\colon \FBFA
    E\to X_0$ satisfying $\hat{S} \eta _E=S$ and $\|\hat{S}\|=1$. Let
    $\hat T\colon Z\to X$ be
    the restriction of the lattice homomorphism $\|T\| \hat{S}$ to
    $Z$. Then $\hat T$ is a lattice homomorphism
    satisfying $\hat T \phi _E=T$. In particular, $\|T\|\le \|\hat T\|\le
    \|T\|\|\hat{S}\|=\|T\|$, so $\|\hat T\|=\|T\|$. An immediate
    consequence of this is that the
    free norm in $\FBFA E$, when restricted to $Z$, is the norm of the free
    Banach lattice, for which an explicit formula was found in
    \cite{aviles_rodriguez_tradacete2018}. From now on, $Z$ will be
    denoted as $\FBL E$.

    Moreover, $\FBL E$ is contractively complemented in $\FBFA E$ by a
    lattice homomorphism. Indeed, the
    canonical inclusion $\eta_E \colon E\to {\FBL E}_0$ extends to a
    unique contractive lattice-algebra homomorphism $P\colon \FBFA
    E\to {\FBL E}_0$. Forget about the zero product in $\FBL E$, and
    consider $P$ merely as a lattice homomorphism. The restriction
    $P|_{\FBL E}\colon \FBL E\to \FBL E$ is a lattice homomorphism of
    norm one satisfying $P \eta _E=\eta _E$. By uniqueness, $P|_{\FBL
    E}$ is the identity on $\FBL E$. Hence the composition of $P$ with
    the canonical inclusion $\FBL E\to \FBFA E$ is a projection of
    $\FBFA E$ onto $\FBL E$.
\end{rem}

\section{Explicit description of the
kernel}\label{sec:explicit_description}

Thus far we have shown, using only general arguments, that $\FNFA E$
and $\FBFA E$ exist for a general Banach space $E$. But since we have
quotiented out $\FAFAv E$ (which we know explicitly) by a mysterious
ideal $\ker \rho $, we have no explicit description for $\FNFA E$.
Fortunately, we can use \cref{thm:falgYudin_ball} to describe this
ideal.

\begin{thm}\label{thm:kernel}
    Let $E$ be a Banach space. Let $\rho $ be the greatest
    submultiplicative lattice seminorm on $\FAFAv E$ such that $\rho
    (\delta _x)\le \|x\|$ for all $x \in E$. Then
    \[
    \ker \rho =\{\, f\in \FAFAv E : f|_{B_{E^{*}}}=0 \, \}.
    \]
\end{thm}
\begin{proof}
    For every $x^{*} \in B_{E^{*}}$,
    $\nu _{x^{*}}(f)=|f(x^{*})|$ defines a lattice seminorm that is
    submultiplicative and such that $\nu _{x^{*}}(\delta _x)\le \|x\|$ for all
    $x \in E$. If $f \in \ker \rho $, then
    \[
        |f(x^{*})|=\nu _{x^{*}}(f)\le \rho (f)=0.
    \]
    Hence every function in $\ker \rho $ vanishes on the unit ball.

    Conversely, suppose $f \in \FAFAv E$ is such that $f|_{B_{E^{*}}}=0$.
    Write $f=\Phi (\delta  _{x_1},\ldots ,\delta _{x_n})$ for a
    certain LLA expression $\Phi $ and some $x_1,\ldots ,x_n \in E$;
    by changing $\Phi $ as necessary, we may assume that $x_1,\ldots
    ,x_n$ are linearly independent and $\|x_i\|\le 1$ (see the proof
    of \cref{prop:fafav}). The goal is to show that
    \[
    q_Ef=\Phi (q_E \delta _{x_1},\ldots ,q_E \delta _{x_n})=\Phi (\eta
    _E(x_1),\ldots ,\eta _E(x_n))
    \]
    is zero in $\FNFA E$. By \cref{thm:structure} there exists an injective and
    contractive lattice-algebra homomorphism $R\colon \FNFA E\to
    X_0\oplus_\infty C(K)$, for a certain Banach lattice $X$ and
    a certain compact Hausdorff space $K$. Denote by $R_1$ (resp.\ $R_2$) the composition of $R$ with
    the projection onto the first (resp.\ second) coordinate. We will
    show that $R_i(q_Ef)=0$ for $i=1,2$.

    For every $x^{*}\in B_{E^{*}}$:
    \[
    0=f(x^{*})=\Phi (x^{*}(x_1),\ldots ,x^{*}(x_n)).
    \]
    We claim that there exists an $\varepsilon >0$ such that $(x^{*}(x_1),\ldots ,x^{*}(x_n))$ takes all values in
    $[-\varepsilon ,\varepsilon ]^{n}$ as $x^{*}$ ranges through
    $B_{E^{*}}$. Indeed, set $F=\barespn\{x_1,\ldots ,x_n\}$, and let
    $\|{\cdot }\|_1$ denote the $\el 1^{n}$-norm in $F$. Since $F$ is
    finite-dimensional, there exists $C>0$ such that $\|{\cdot
    }\|_1\le C\|{\cdot }\|$, where $\|{\cdot }\|$ stands for the
    norm in $E$. Define $\varepsilon =1/C$, and let $x^{*}\colon F\to
    \R$ be a linear functional satisfying $(x^{*}(x_1),\ldots
    ,x^{*}(x_n))\in [-\varepsilon ,\varepsilon ]^{n}$. If
    $x=\sum_{i=1}^{n}\lambda _i x_i$ is an arbitrary element of $F$,
    then
    \[
    |x^{*}(x)|\le \sum_{i=1}^{n}|\lambda _i| |x^{*}(x_i)|\le
    \varepsilon \|x\|_1\le \|x\|.
    \]
    This shows that $x^{*}\in B_{F^{*}}$. Use Hahn--Banach to extend
    this functional to an element $x^{*}\in B_{E^{*}}$. Since $(x^{*}(x_1),\ldots
    ,x^{*}(x_n))$ can be chosen arbitrarily in $[-\varepsilon
    ,\varepsilon ]^n$, the desired claim follows.
    It is a consequence of this claim and \cref{lem:calc_zero} that $\Phi $
    vanishes on every Banach lattice with the zero product; in
    particular, it vanishes on $X_0$, and therefore $R_1(q_E(f))=0$.

    Since $R_2$ is a
    contractive lattice-algebra homomorphism, $R_2=\hat{T}$
    for a certain contractive operator $T\colon E\to C(K)$ satisfying
    $R_2\eta _E=T$. For $t \in K$, denote by $\phi _t \in C(K)^{*}$ the
    evaluation functional at $t$. Recall that $\phi _t$ is a
    lattice-algebra homomorphism of norm one. Compute:
    \begin{align*}
        R_2(q_E(f))(t)&=\Phi (R_2(\eta _E(x_1))(t),\ldots ,R_2(\eta
        _E(x_n))(t))\\
                      &=\Phi ((Tx_1)(t),\ldots ,(Tx_n)(t))\\
                      &=\Phi ( (T^{*}\phi _t)(x_1),\ldots
                      ,(T^{*}\phi  _t)(x_n))\\
                      &=\Phi (\delta _{x_1}(T^{*}\phi _t),\ldots,
                      \delta _{x_n}(T^{*}\phi _t))\\
                      &=\Phi (\delta _{x_1},\ldots ,\delta
                      _{x_n})(T^{*}\phi _t)\\
                      &=f(T^{*}\phi _t)=0
    \end{align*}
    since $\|T^{*}\phi _t\|\le 1$ and $f$ vanishes on the unit ball.
    This shows that $R_2(q_E(f))=0$. Hence $R(q_E(f))=0$ and, since $R$
    is injective, $q_E(f)=0$.
\end{proof}

Since $\delta _x\colon E^{*}\to \R$ is weak$^*$ continuous, so is every
element of $\FAFAv E$. Therefore there is a restriction map:
\[
\begin{array}{cccc}
        & \FAFAv E & \longrightarrow & C(B_{E^{*}},w^{*}) \\
        & f & \longmapsto & f|_{B_{E^{*}}} \\
\end{array}.
\]
Operations in $\FAFAv E$ are computed pointwise, so this is a
lattice-algebra homomorphism. However, as noted in
\cref{rem:not_injective}, this map is not injective, so
we cannot view the elements of $\FAFAv E$ as continuous functions on
the unit ball of the dual.

It is remarkable that this map does induce an injective
lattice-algebra homomorphism $\iota \colon \FNFA E\to C(B_{E^{*}})$
because, according to \cref{thm:kernel}, its kernel is precisely $\ker
\rho $. Note that, by definition of the norm in $\FNFA E$, $\iota $ is
contractive. We have thus found an explicit lattice-algebraic
representation of $\FNFA E$.

\begin{cor}\label{cor:iota}
    Let $E$ be a Banach space. There exists a unique injective and
    contractive lattice-algebra homomorphism
    \[
    \iota_E \colon \FNFA E\to C(B_{E^{*}})
    \]
    satisfying $(\iota_E \eta _x)(x^{*})=x^{*}(x)$ for all $x \in E$ and
    $x^{*} \in B_{E^{*}}$.
\end{cor}

From now on, $\iota _E$ (or simply $\iota $ when $E$ is clear from
the context) will denote the map in the above corollary. The
following is a first application of this representation result.

\begin{prop}\label{prop:strong_unit}
    Let $E$ be a Banach space and let $f \in \FNFA E$. There exist $\lambda >0$ and $x_1,\ldots ,x_n
    \in E$ such that
    \[
    |f|\le \lambda |\eta _{x_1}|\vee \cdots \vee |\eta _{x_n}|.
    \]
\end{prop}
\begin{proof}
    As already mentioned in \cref{sec:abstract_construction}, there exists a polynomial $p(t_1,\ldots ,t_n)\in
    \R_+[t_1,\ldots ,t_n]$ with $p(0,\ldots ,0)=0$ such that $|f|\le p(|\eta  _{x_1}|,\ldots
    ,|\eta  _{x_n}|)$. Without loss of generality, $x_1,\ldots ,x_n \in B_E$. It follows
    that for every $x^{*} \in B_{E^{*}}$:
    \begin{align*}
        \iota (|f|)(x^{*})&\le p(|x^{*}(x_1)|,\ldots ,|x^{*}(x_n)|)\\&\le p(1,\ldots
    ,1)\, |x^{*}(x_1)|\vee \cdots \vee |x^{*}(x_n)|\\
                          &=p(1,\ldots ,1)\iota (|\eta _{x_1}|\vee
                          \cdots \vee |\eta _{x_n}|)(x^{*}),
    \end{align*}
    where in the second inequality we have used the fact that
    $|x^{*}(x_i)|\le 1$. Since $\iota $ is an injective lattice
    homomorphism, the result follows.
\end{proof}

\begin{rem}
    Since, a priori, $\iota \colon \FNFA E\to C(B_{E^{*}})$ is not an embedding
    (in fact, we will show in
    \cref{rem:notembed} that it is never an embedding), it is not
    clear at all whether its extension to
    $\FBFA E$ is injective. In \cref{sec:representation} we will
    use a completely different technique to show that the extension is indeed
    injective for some classes of Banach spaces.
\end{rem}

\begin{rem}
    Recall from \cref{rem:FBLinFBFA} that $\FBL E$ can be seen
    as a closed sublattice of $\FBFA E$; moreover, $\FBL E$ is
    $1$-complemented by a lattice homomorphic projection $P$. Using the
    representation of $\FNFA E$ inside $C(B_{E^{*}})$, we can now give an explicit
    expression for $P$ on $\FNFA E$. We claim that, for every $f \in
    \FNFA E$ and every $x^{*}\in B_{E^{*}}$,
    \[
        \iota (Pf)(x^{*})=\lim_{t\to 0^{+}} \frac{\iota (f)(tx^{*})}{t}.
    \]
    Indeed, consider the map $f\mapsto Pf$ defined by the above
    formula.
    According to \cref{lem:lla_limit} the limit exists, and $\iota
    (Pf)$ belongs to $\iota (\Lat \{\, \eta _x : x \in
    E\, \})$; that is, $Pf$ belongs to $\FVLv E\subseteq \FBL E$ (see
    \cref{rem:FVLinFAFA} for further details). The operator $P$, as defined by the formula above,
    is a lattice-algebra homomorphism $\FNFA E\to \FBL E_0$
    satisfying $P \eta _x =\eta _x$. Therefore, $P$ is the unique
    extension of $\eta \colon E\to {\FBL E}_0$ to $\FNFA E$. The extension of
    $P$ by density to the whole $\FBFA E$ must be, again by
    uniqueness, the contractive lattice projection described in
    \cref{rem:FBLinFBFA}.
\end{rem}

\begin{prop}\label{prop:FNFAsemiprime}
    The free normed \falg\ $\FNFA E$ generated by a Banach space $E$
    is semiprime.
\end{prop}
\begin{proof}
    We will use the representation $\iota \colon \FNFA E\to
    C(B_{E^{*}})$.
    Let $f \in N(\FNFA E)$. For every $x^{*}\in B_{E^{*}}$ different
    from $0$ there exists an $x \in B_E$ such that $\iota (\eta
    _x)(x^{*})=x^{*}(x)\neq 0$. Then
    \[
    0=\iota (f\eta _x)(x^{*})=\iota (f)(x^{*})\,x^{*}(x).
    \]
    It follows that $\iota (f)(x^{*})=0$. Since also $\iota (f)(0)=0$,
    it must be $\iota (f)=0$. But $\iota $ is injective, so $f=0$.
\end{proof}

Although the lattice-algebraic representation of $\FNFA E$ provided in
\cref{cor:iota} has proven useful, not much else can be said about
free normed \falg s unless we investigate the properties of its
norm. The following section is a first step in this direction: it
shows in particular that, when $E$ is finite-dimensional, the free norm
is equivalent to an AM-norm.

\section{The finite-dimensional case}\label{sec:fin_dim}

The goal of this section is to prove the following isomorphic
description of the space $\FBFA E$ when $E$ is a finite-dimensional
Banach space.

\begin{thm}\label{thm:finite}
    Let $E$ be a finite-dimensional Banach space. The free Banach \falg\
    $\FBFA E$ is lattice-algebra isomorphic to the space
    $C([0,1]\times S_{E^{*}})$ equipped
    with pointwise order and product
    \begin{equation}\label{eq:free_prod}
        (f\star g)(r,u)=r f(r,u) g(r,u)
    \end{equation}
    where $(r,u)\in [0,1]\times S_{E^{*}}$ and $f,g \in C([0,1]\times
    S_{E^{*}})$. Under this isomorphism, $\eta _x(r,u)=u(x)$ for every
    $x \in E$.
\end{thm}

\begin{figure}[htp]
    \centering
    \begin{subfigure}[b]{0.3\textwidth}
        \centering
        \includegraphics[width=\textwidth]{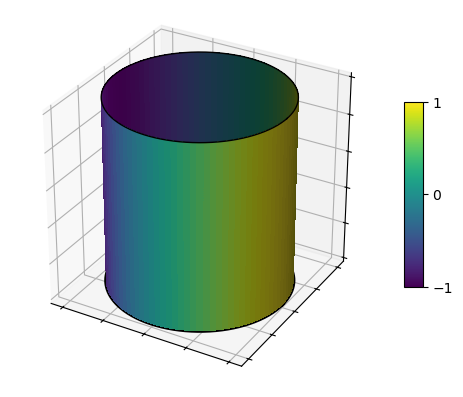}
        \caption{$\eta _{e_1}(r,u)=u(e_1)$}
    \end{subfigure}%
    \begin{subfigure}[b]{0.3\textwidth}
        \centering
        \includegraphics[width=\textwidth]{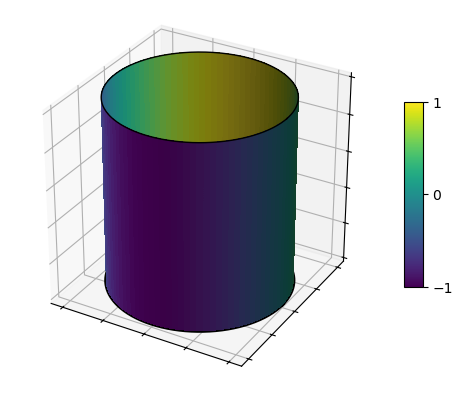}
        \caption{$\eta _{e_2}(r,u)=u(e_2)$}
    \end{subfigure}%
    \begin{subfigure}[b]{0.3\textwidth}
        \centering
        \includegraphics[width=\textwidth]{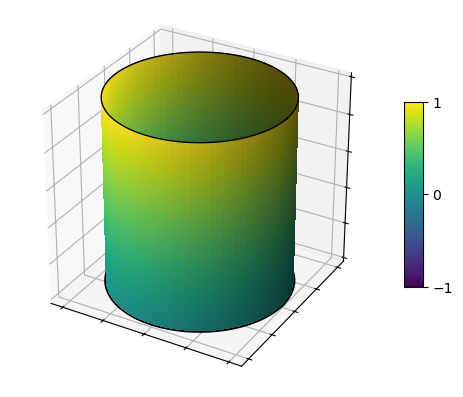}
        \caption{$(\one \star \one)(r,u)=r$}
    \end{subfigure}%
    \caption{Representation of three functions in $C([0,1]\times
        S_{\el 2^2})$. In the plots, $[0,1]$ is identified with $0\le
        z\le 1$ and $S_{\el 2^2}$ is identified with the unit circle
        in the plane $z=0$. The elements $e_1,e_2$ denote the
        canonical basis of $\el 2^2$.}
    \label{fig:FBFA}
\end{figure}

\Cref{fig:FBFA} illustrates the generators and the ``weight'' $\one
\star \one$ (where $\one$ denotes the constant one function) of the free Banach \falg\ generated by $\el 2^2$.

To make the steps in the proof clearer, we have divided it into several
lemmas. First, we want to show that there exists an element in $\FNFA
E$ that is a strong unit in $\FBFA E$, and that the norm it induces is
equivalent to the free norm.

\begin{lem}\label{lem:finite_prelim1}
    Let $E$ be a finite-dimensional Banach space. There exists a
    finite subset $F\subseteq S_E$ for which $e'=\sup_{x \in F}|\eta
    _x|$ is a strong unit in $\FBFA E$ and
    such that the norm it induces $\|{\cdot }\|_{e'}$ is equivalent to
    the free norm.
\end{lem}
\begin{proof}
    We shall first show that $\FBFA E$ has a strong unit $e \in \FBFA
    E$. Denote by $\Fin (S_{E})$ the family of finite subsets of
    $S_E$. For every $F \in \Fin (S_{E})$, define the element
    $e_F=\sup_{x \in F}|\eta _x| \in \FNFA E$. In fact, this element
    belongs to $\FVLv E$, where we are identifying $\FVLv E$ inside
    $\FNFA E$ as in \cref{rem:FBLinFBFA}. Keep in mind that, with this
    identification, the norm of $\FNFA E$, when restricted to $\FVLv
    E$, is precisely the norm of the free Banach lattice. Hence
    we can refer to the free norm on $\FVLv E$ without ambiguity.

    Let $\iota \colon \FNFA E\to C(B_{E^{*}})$
    be the injective lattice-algebra homomorphism from
    \cref{cor:iota}. For every $x^{*}\in B_{E^{*}}$, the net
    $(\iota (e_F)(x^{*}))_{F \in \Fin(S_E)}$ converges to $\|x^{*}\|_{E^{*}}$. Since $E$ is
    finite-dimensional, $\|{\cdot }\|_{E^{*}}$ is a continuous
    function on $B_{E^{*}}$ (with respect to the weak$^*$ topology). Moreover, $(e_F)_{F \in \Fin (S_{E})}$
    is an increasing net. By Dini's theorem, $\iota (e_F)$ converges
    uniformly to $\|{\cdot }\|_{E^{*}}$ in $C(B_{E^*})$.

    Recall from \cite[Section 9.1]{oikhberg_etal2022} that, because
    $E$ is finite-dimensional, the free norm is equivalent to the
    uniform norm on $\iota (\FVLv E)$; that is, $\iota $ is an
    embedding when restricted to $\FVLv E$. We have shown above that the net
    $(\iota (e_F))$ is Cauchy in the uniform norm. Hence the net
    $(e_F)$ is Cauchy in the free norm.

    Let $e \in \FBFA E$ be the limit of $(e_F)$. Since $(e_F)$ is an
    increasing net, $e=\sup_{F \in \Fin(S_E)}e_F$. We are
    going to show that $I_e$, the order ideal generated by $e$, is the whole
    $\FBFA E$ (and therefore that $e$ is a strong unit). For this, we
    first need to show that the free norm $\|{\cdot }\|$ and $\|{\cdot
    }\|_e$ are equivalent. It is
    clear that $\iota (e_F)^2\le \iota (e_F)$ holds pointwise in
    $B_{E^{*}}$ for every $F \in \Fin(S_E)$. Hence $e_F^2\le e_F$ and,
    taking limits, $e^2\le
    e$. It follows that, for every $f,g \in I_e$,
    \[
    |fg|\le |f| |g|\le \|f\|_e \|g\|_e e^2\le \|f\|_e \|g\|_e e.
    \]
    Therefore $fg \in I_e$ and $\|fg\|_e\le \|f\|_e \|g\|_e$. Also
    note that, if $x \in S_{E}$, then $|\eta _x|= e_{\{x\}}\le e$.
    In particular, $\FNFA E \subseteq I_e$ and $\|\eta _x\|_e\le
    \|x\|$ for every $x \in E$. We have thus shown that $\|{\cdot
    }\|_e$, when restricted to $\FNFA E$, is a submultiplicative
    lattice norm satisfying $\|\eta _x\|_e\le \|x\|$ for all $x \in
    E$. By construction of the free norm, $\|{\cdot }\|_e\le \|{\cdot
    }\|$ on $\FNFA E$.

    On the other hand, $|f|\le \|f\|_e e$ holds for every $f \in \FNFA
    E$. By taking the free norm on both sides of the inequality,
    $\|f\|\le \|f\|_e \|e\|$. This completes the proof that the free
    norm is equivalent to $\|{\cdot }\|_e$ on $\FNFA E$. Now if $f \in
    \FBFA E$ is arbitrary, and $f_n\to f$ in the free norm, with $f_n
    \in \FNFA E$, then $(f_n)$ is also a Cauchy sequence in the norm
    $\|{\cdot }\|_e$. Since this norm is complete, $(f_n)$ converges
    to an element of $I_e$; since this norm is equivalent to the free
    norm, the limit must be $f$. Hence, $f \in I_e$. This shows that
    $e$ is a strong unit, and that the norm $\|{\cdot }\|_e$ is
    equivalent to the free norm on $\FBFA E$.

    To finish the proof, it only remains to show that there exists a
    strong unit $e' \in \FNFA E$ of the desired form equivalent to
    $e$. For this, it
    suffices to find $e' \in \FNFA E$ of the desired form satisfying $e/2\le e'\le e$.
    (The factor $1/2$ does not play any special role in the proof; any scalar in $(0,1)$ would also work.) Let
    $F \in \Fin(S_E)$ be such that
    \[
        \|\iota (e_F)-\|{\cdot }\|_{E^{*}}\|_\infty \le 1/2
    \]
    and put $e'=e_F$. Certainly, $e'\le e$. For every $u \in S_{E^{*}}$,
    $|\iota (e')(u)-1|\le 1/2$, which implies $\iota (e')(u)\ge 1/2$.
    Since $\iota (e')$ is
    positively homogeneous, it follows that
    \[
        \iota (e')(x^{*})\ge 1/2 \|x^{*}\|_{E^{*}}\ge
        \iota (e_{F'}/2)(x^{*})
    \]
    for all $x^{*}\in B_{E^{*}}$ and all $F' \in \Fin(S_{E})$. Hence
    $e'\ge e_{F'}/2$ in $\FBFA E$ for all $F' \in \Fin(S_E)$. Taking
    the limit in $F' \in \Fin(S_{E})$, it follows that $e'\ge e/2$.
\end{proof}

\begin{rem}
    The subset $F$ in the previous lemma is far from unique. In fact, if
    $F\subseteq F'$, then $F'$ also satisfies the desired properties.
\end{rem}

From the previous lemma and Kakutani's theorem it already follows that
$\FBFA E$ is lattice isomorphic to $C(K)$, for a certain compact
Hausdorff space $K$.

\begin{lem}\label{lem:finite_prelim2}
    Let $E$ be a finite-dimensional Banach space. Let $C([0,1]\times
    S_{E^{*}})$ be equipped with pointwise order and product $\star$ given by
    \eqref{eq:free_prod}. Then the map $T\colon E\to C([0,1]\times
    S_{E^{*}})$, defined by
    \[
        (Tx)(r,u)=u(x)\quad\text{for all }x \in E, (r,u)\in
        [0,1]\times S_{E^{*}},
    \]
    extends to a contractive lattice-algebra homomorphism
    $\hat{T}\colon \FBFA E\to C([0,1]\times S_{E^{*}})$ whose explicit
    expression is
    \[
        (\hat{T}f)(r,u)=
    \begin{cases}
        \displaystyle\frac{\iota (f)(ru)}{r}&\text{if }r\neq 0\\
        \displaystyle\lim_{r\to 0^{+}}\frac{\iota (f)(ru)}{r}&\text{if }r=0
    \end{cases}
    \]
    for every $f \in \FNFA E$.
\end{lem}
\begin{proof}
    It is direct to check that $C([0,1]\times S_{E^{*}})$
    with the usual lattice structure, the product $\star$, and the
    uniform norm, is a Banach \falg. Note that
    \[
    \|Tx\|_\infty =\sup_{u \in S_{E^{*}}} |u(x)|=\|x\|.
    \]
    Hence $T$ extends in a unique way to a norm-one lattice-algebra
    homomorphism $\hat{T}\colon \FBFA E\to C([0,1]\times S_{E^{*}})$.

    Denote momentarily
    \[
    Sf(r,u)=
    \begin{cases}
        \displaystyle\frac{\iota (f)(ru)}{r}&\text{if }r\neq 0\\
        \displaystyle\lim_{r\to 0^{+}}\frac{\iota (f)(ru)}{r}&\text{if }r=0
    \end{cases}
    \]
    for every $(r,u) \in [0,1]\times S_{E^{*}}$ and $f \in \FNFA E$.
    To show that $\hat{T}|_{\FNFA E}=S$, we will check that the map
    $S\colon \FNFA E\to C([0,1]\times S_{E^{*}})$, $f\mapsto Sf$, is a
    lattice-algebra homomorphism extending $T$. Uniqueness of the
    extension will then imply that $\hat{T}|_{\FNFA E}=S$.

    First of all, note that the limit defining $Sf(0,u)$ exists by
    \cref{lem:lla_limit}. We need to check that $Sf$ is continuous on
    $[0,1]\times S_{E^{*}}$. It is clear that it
    is continuous on $\{\,(r,u)\in [0,1]\times S_{E^{*}}:r\neq 0\,\}$.
    Let $(r_n,u_n)\in [0,1]\times S_{E^{*}}$ be a sequence converging
    to a certain $(0,u) \in [0,1]\times S_{E^{*}}$. The limit
    defining $Sf(0,u)$ coincides with
    $f_0(u)=\lim_{n}f(r_nu)r_n^{-1}$ for every $u \in
    S_{E^{*}}$. Since $E$ is finite-dimensional,
    $S_{E^{*}}$ is a compact metric space under the metric induced by
    the norm. Moreover, $S_{E^{*}}$ can be embedded isometrically as a compact subspace of $\R^{n}$.
    Then $f_0$ is, by \cref{lem:lla_limit}, the uniform limit of
    continuous functions. It is, therefore, continuous. Furthermore,
    since $S_{E^{*}}$ is compact, $f_0$ is uniformly continuous on
    $S_{E^{*}}$.

    Fix $\varepsilon >0$. Choose $\delta >0$ such that
    $|f_0(u)-f_0(u')|<\varepsilon/2$ whenever $\|u-u'\|_{E^{*}}<\delta $ and
    such that $|f(r_nu)r_n^{-1}-f_0(u)|<\varepsilon/2$ whenever $0\le
    r_n<\delta $ and $u \in S_{E^{*}}$ (recall the limit
    defining $f_0(u)$ is uniform). For $n \in \N$
    big enough we will have $r_n,\|u-u_n\|<\delta $, and then,
    whenever $r_n\neq 0$,
    \begin{align*}
        |Sf(r_n,u_n)-f_0(u)|&=|f(r_n u_n)r_n^{-1} - f_0(u)|\\&\le |f(r_n u_n)r_n^{-1} -
    f_0(u_n)|+|f_0(u_n)-f_0(u)|<\varepsilon.
    \end{align*}
    When $r_n=0$, $Sf(r_n,u_n)=f_0(u_n)$ and
    $|f_0(u_n)-f_0(u)|<\varepsilon /2$.
    This proves that $Sf(r_n,u_n)\to f_0(u)=Sf(0,u)$, and therefore that $Sf$
    is a continuous function.

    It is clear
    from the definition that $S$ is a lattice homomorphism. It is also
    an algebra homomorphism because, for $r\neq 0$,
    \[
    S(fg)(r,u)=\iota (fg)(ru)r^{-1}=r Sf(r,u) Sg(r,u)=(Sf\star Sg)(r,u).
    \]
    The equality also holds for general $(r,u)\in [0,1]\times
    S_{E^{*}}$ by continuity. Certainly, $(S\eta _x)(r,u)=u(x)=(Tx)(r,u)$. Hence $Sf=\hat{T}f$ for all $f \in \FNFA E$.
\end{proof}

With the explicit expression we can now show that the map $\hat{T}$ is
injective on $\FNFA E$.

\begin{lem}\label{lem:finite_prelim3}
    Let $E$ be a finite-dimensional Banach space, and let $\hat{T}$ be
    as in \cref{lem:finite_prelim2}. If, for some $f,g \in \FNFA E$,
    $\hat{T}f\le \hat{T}g$, then $f\le g$.
\end{lem}
\begin{proof}
    Let us first check that $\hat{T}$ is injective on $\FNFA E$.
    Suppose $\hat{T}f=0$ for some $f \in \FNFA E$. Then, for every
    $u \in S_{E^{*}}$ and $0<r\le 1$,
    \[
        0=(\hat{T}f)(r,u)=\iota (f)(ru)/r.
    \]
    This implies that $\iota (f)$ vanishes on $B_{E^{*}}\setminus \{0\}$.
    Since every element of $\FNFA E$ vanishes also at $0$, it follows
    that $\iota (f)=0$. Injectivity of $\iota $ implies that $f=0$. Hence the map $\hat{T}|_{\FNFA E}\colon \FNFA E\to
    \hat{T}(\FNFA E)$ is a lattice isomorphism. In particular, if
    $\hat{T}f\le \hat{T}g$, then by applying the inverse it follows
    that $f\le g$.
\end{proof}

\begin{proof}[Proof of \cref{thm:finite}]
    Let $e' =\sup_{x \in F}|\eta _x|\in \FNFA E$ be as in \cref{lem:finite_prelim1}, where
    $F \in \Fin(S_E)$, and let
    $\hat{T}\colon \FBFA E\to C([0,1]\times S_{E^{*}})$ be as in
    \cref{lem:finite_prelim2}. First note that
    \[
    (\hat{T}e')(r,u)=\sup_{x \in F} |u(x)|\ge 1/2
    \]
    (see the end of the proof of \cref{lem:finite_prelim1}).
    Therefore $\hat{T}e'$ is a strong unit in $C([0,1]\times
    S_{E^{*}})$. We are going to show that
    $\|\hat{T}f\|_{\hat{T}e'}=\|f\|_{e'}$ for all $f \in \FNFA E$.

    Since $|f|\le \|f\|_{e'}e'$, applying $\hat{T}$ to this inequality
    yields immediately $\|\hat{T}f\|_{\hat{T}e'}\le \|f\|_{e'}$.
    Similarly, $\hat{T}|f|=|\hat{T}f|\le
    \|\hat{T}f\|_{\hat{T}e'}\hat{T}e'$ and from
    \cref{lem:finite_prelim3} it follows that $|f|\le
    \|\hat{T}f\|_{\hat{T}e'}e'$. Hence also $\|f\|_{e'}\le
    \|\hat{T}f\|_{\hat{T}e'}$.

    Since the free norm is equivalent to $\|{\cdot }\|_{e'}$ and the
    uniform norm is equivalent to $\|{\cdot }\|_{\hat{T}e'}$, the
    equality $\|\hat{T}f\|_{\hat{T}e'}=\|f\|_{e'}$ must hold for all
    $f \in \FBFA E$. In particular, $\hat{T}$ is injective and its
    image is closed under the uniform norm.

    To complete the proof it only remains to show that $\hat{T}$ is
    surjective. This will follow as a consequence of the Stone--Weierstrass
    theorem once we show that the image of $\hat{T}$ contains the
    constant one function and separates points.

    Define, as in the proof of \cref{lem:finite_prelim1}, the element $e_F=\sup_{x \in F}|\eta _x|$ for every $F
    \in \Fin(S_E)$. We
    already showed in the lemma that $(e_F)$ converges to a certain $e \in \FBFA E$ in
    the free norm. Note that
    \[
        (\hat{T}e_F)(r,u)=\sup_{x \in F}|Tx|(r,u)=\sup_{x \in
        F}|u(x)|
    \]
    converges to $1$. Therefore
    $\hat{T}e=\lim_{F}\hat{T}e_F=\one$, where $\one$ denotes the
    constant one function.

    Let $(r,u), (r',u') \in [0,1]\times S_{E^{*}}$. If $u\neq u'$,
    then there exists $x \in S_E$ such that $u(x)\neq u'(x)$. In other
    words:
    \[
        (\hat{T}\eta _x)(r,u)=u(x)\neq u'(x)=(\hat{T}\eta _x)(r',u').
    \]
    If instead $r\neq r'$, then we have
    \[
        (\hat{T}e^2)(r,u)=(\hat{T}e\star \hat{T}e)(r,u)=r
    \]
    and therefore $(\hat{T}e^2)(r,u)\neq (\hat{T}e^2)(r',u')$. Thus
    the image of $\hat{T}$ separates points of $[0,1]\times S_{E^{*}}$. From the Stone--Weierstrass
    theorem it follows that $\hat{T}$ is surjective.
\end{proof}

\begin{rem}\label{rem:notembed}
    The representation $\iota \colon \FNFA E\to C(B_{E^{*}})$ is never
    an embedding. Indeed, when $E$ is infinite-dimensional, the norm
    of the free Banach lattice cannot be equivalent to an
    AM-norm (this is a consequence of \cite[Proposition 9.30]{oikhberg_etal2022}).
    Suppose that $\iota $ were an embedding for some
    finite-dimensional $E$; we will show
    that this contradicts \cref{thm:finite}. In this case, the image
    of $\iota $ would be a closed sublattice-algebra of $J_0$, the
    closed (lattice and algebraic) ideal of
    $C(B_{E^{*}})$ formed by the functions that vanish at $0$. It is immediate to check that
    $\iota (\FNFA E)$ separates the points of $B_{E^{*}}$, and that
    for every $x^{*} \in B_{E^{*}}$ different from zero there exists
    $f \in \FNFA E$ such that $\iota (f)(x^{*})\neq 0$. As a
    consequence of the
    Stone--Weierstrass theorem, $\iota (\FNFA E)=J_0$. But $J_0$ does
    not have a strong order unit; this contradicts
    \cref{thm:finite}.
\end{rem}

As we shall see in \cref{sec:representation}, determining whether
$\FBFA E$ is semiprime or not is very important and yet far from
trivial.

\begin{cor}\label{cor:finite_semiprime}
    If $E$ is a finite-dimensional Banach space, then $\FBFA E$ is
    semiprime.
\end{cor}
\begin{proof}
    By \cref{thm:finite},
    we only need to check that the product \eqref{eq:free_prod} is
    semiprime. Suppose $f\star f=0$. Then
    \[
    0=(f\star f)(r,u)=r f(r,u)^2
    \]
    for all $(r,u)\in [0,1]\times S_{E^{*}}$. This means that $f$
    vanishes on the dense set $\{(r,u) \in [0,1]\times S_{E^{*}}:
    r\neq 0\}$. By continuity, $f=0$.
\end{proof}

In the next result, $\Sinfn$ denotes the unit sphere of $\el \infty ^{n}$.

\begin{cor}
    Let $E$ be a Banach space of dimension $n$. The free Banach \falg\
    $\FBFA E$ is lattice-algebra isomorphic to $C([0,1]\times \Sinfn)$
    with pointwise order, and product
    \[
        (f\star g)(r,u)=r f(r,u) g(r,u)
    \]
    where $(r,u)\in [0,1]\times \Sinfn$  and $f,g \in C([0,1]\times
    \Sinfn)$.
\end{cor}
\begin{proof}
    Identify, by taking coordinates, $S_{E^{*}}$ and $\Sinfn$ with
    their homeomorphic copies in $\R^{n}$. The map
    \[
    \begin{array}{cccc}
        \phi \colon& [0,1]\times S_{E^{*}} & \longrightarrow & [0,1]\times \Sinfn \\
            & (r,u) & \longmapsto & (r,u/\|u\|_\infty)  \\
    \end{array}
    \]
    is a homeomorphism. Consider the associated composition operator
    \[
    \begin{array}{cccc}
        T\colon& C([0,1]\times \Sinfn) & \longrightarrow &
        C([0,1]\times S_{E^{*}}) \\
            & f & \longmapsto & f\circ \phi \\
    \end{array}.
    \]
    This map is a lattice isometry. Moreover, if $f,g \in C([0,1]\times \Sinfn)$ and
    $(r,u)\in [0,1]\times S_{E^{*}}$, then
    \begin{align*}
        T(f\star g)(r,u)&=(f\star g)(\phi (r,u))\\
                        &=(f\star g)(r,u/\|u\|_\infty )\\
                        &=rf(r,u/\|u\|_\infty)g(r,u/\|u\|_\infty)\\
                        &=r (f\circ \phi )(r,u) (g\circ \phi )(r,u)\\
                        &=(Tf\star Tg)(r,u).
    \end{align*}
    This shows that $T$ is also an algebra homomorphism. The stated
    result now follows from \cref{thm:finite}.
\end{proof}

It might be pertinent to recall here that the free Banach lattice over a Banach space of dimension $n$ is lattice isomorphic to $C(\Sinfn)$, but the isomorphism constant increases without bound as the dimension grows (see \cite[Section 9]{oikhberg_etal2022}).

\begin{cor}\label{cor:dim_isomorphic}
    Let $E$ and $F$ be Banach spaces of dimension $n$. Then $\FBFA E$
    and $\FBFA F$ are lattice-algebra isomorphic.
\end{cor}

\begin{rem}
    Even though isomorphic Banach spaces generate isomorphic free
    Banach lattices, this is not clear at all for free Banach
    \falg s. In fact, in \cref{ex:extension_of_bijection} we will
    provide an isomorphism between Banach spaces that does not induce an
    isomorphism between the corresponding free Banach \falg s.
    It is clear that isometric Banach spaces generate lattice-algebra
    isometric free Banach \falg s. The isomorphic case, however,
    remains open in infinite dimensions.
\end{rem}

\begin{question}
    Let $E$ and $F$ be isomorphic Banach spaces. Are $\FBFA E$ and
    $\FBFA F$ lattice-algebra isomorphic?
\end{question}

The converse of the previous question is open even for the free Banach
lattice generated by a Banach space (see \cite[Section
10]{oikhberg_etal2022}). The isometric version of the converse will be
addressed in \cref{sec:isometries}.

We have completely identified the free Banach \falg\ generated by a
finite-dimensional Banach space, except for the free norm. The goal of
the next section is to better understand this norm.

\section{On the free norm}\label{sec:free_norm}

It follows from the abstract construction in \cref{sec:abstract_construction} that the norm in $\FBFA E$ is the supremum of the lattice
seminorms that are submultiplicative and that make the canonical map $\eta
_E\colon E\to \FBFA E$ contractive. There is an equivalent way of
defining this norm: for every $f \in \FBFA E$, $\|f\|$ is the
smallest number such that, whenever $A$ is a Banach
\falg\ and $T\colon E\to A$ is a contractive operator,
$\|\hat{T}f\|\le \|f\|$. Moreover, if $f \in \FBL E$, then
\[
\|f\|=\sup \{\, \|\hat{T}f\| \;:\; T\colon E\to \el 1^{n}\text{ with } \|T\|\le 1 \, \}.
\]
This important fact was first proved in
\cite{aviles_rodriguez_tradacete2018} and is at the center of the
theory of free Banach lattices. For this reason, in this section we
try to reach a similar result for free Banach \falg s. Only
considering operators that arrive at $\el 1^{n}$ is too ambitious;
instead, we will take the supremum over all operators that arrive at
finite-dimensional semiprime \falg s. Although the semiprimeness
condition may look artificial at first, it will turn out to be
very useful in \cref{sec:representation}.

\begin{defn}
    Let $E$ be a Banach space. For every $f \in \FBFA E$, define
    $\tau_E (f)$ to be the least positive number such that, if $A$ is a
    semiprime finite-dimensional Banach \falg, and
    $T\colon E\to A$ is a contractive operator, then $\tau_E (f)\ge \|\hat{T}f\|$.
    When there is no ambiguity on the underlying Banach space, $\tau _E$ is simply
    denoted by $\tau $.
\end{defn}

\begin{rem}
    Note that $\tau (f)$ is finite since $\|\hat{T}f\|\le
    \|f\|$ for every contractive operator $T\colon E\to A$. In
    particular, $\tau (f)\le \|f\|$. It is not difficult to check from
    the definition that $f\mapsto \tau (f)$ defines a
    submultiplicative lattice seminorm on $\FBFA E$.
\end{rem}

It will be very valuable for the discussion of
\cref{sec:representation} (and in particular for
\cref{prop:tau_norm_implies_semiprime}) to determine when $\tau$ defines a norm on $\FBFA E$ (i.e., when $\tau (f)=0$ implies
$f=0$). The following result is an important step in this direction.

\begin{thm}\label{thm:norm_fin_dim}
    The norm in $\FBFA{\el 1^{n}}$ is $\tau $ for every $n \in \N$.
\end{thm}

The main computation of the theorem is encapsulated in the next
lemma. Since it is quite technical, let us explain the ideas that
motivate it. From \cref{thm:finite} we know that $\FBFA{\el 1^{n}}$
is isomorphic to $C(K)$, for a certain compact Hausdorff space $K$,
with a product of the form $f\star g=wfg$, where $w \in C(K)_+$ and
juxtaposition denotes the pointwise product. (We are not being
explicit with the values of $K$ and $w$ on purpose; the argument will
not depend on this.) Let $f=\Phi (\eta _{e_1},\ldots ,\eta _{e_n}) \in
\FBFA E$, for a certain LLA expression
$\Phi $. Approximate the generators
$\eta _{e_1},\ldots ,\eta _{e_n}$ and the weight $w$ by simple
functions (i.e., functions that take only finitely many values).
Disjointify the family of characteristic functions involved in these
simple functions. Then their span will be a finite-dimensional
sublattice $A$ that contains discrete approximations of $\eta _{e_1},\ldots
,\eta _{e_n}$. Define a product in this space that consists in
multiplying pointwise by the discrete approximation of $w$. If the
approximations are good enough, this new product will be almost
submultiplicative for the free norm, and $\Phi $, when evaluated at
the discrete approximations of the $\eta _i$, will be very close to
$f$. This is the statement of the lemma. Afterwards, in the proof of the theorem, we will use this
finite-dimensional \falg\ $A$, with an appropriate norm, to construct
an operator $T\colon \el 1^{n}\to A$ for which $\|\hat{T}f\|_A$ is
very close to the free norm of $f$. One caveat though: since $K$ is not sufficiently disconnected, we may not be
able to approximate the functions by simple functions. For this reason
we shall work in $B(K)$, the Banach lattice of bounded
Borel-measurable functions on $K$.

\begin{lem}\label{lem:fin_dim_falg}
    Let $K$ be a compact Hausdorff space, and let $C(K)$ be equipped
    with the usual lattice structure and a certain product $\star$ that makes
    it a Banach \falg. Let $f_1,\ldots ,f_n \in C(K)$, and let $f=\Phi
    (f_1,\ldots ,f_n)$ for some LLA expression $\Phi $. Fix
    $\varepsilon, \varepsilon '>0$. Then there exists a finite-dimensional sublattice $A$
    of $B(K)$, together with a semiprime \falg\ product $\circ$ in
    $A$, and
    $(f_1)_d,\ldots ,(f_n)_d \in A$ such that $|(f_i)_d|\le |f_i|$ and, if $f_d=\Phi
    ((f_1)_d,\ldots ,(f_n)_d)$ (where now $\Phi $ is evaluated using
    the product of $A$), then $\|f-f_d\|_\infty <\varepsilon $.
    Moreover, $\circ$ can be chosen so that
    \[
    |x\circ y|\le |x\star y| + \varepsilon '\one\footnote{Here $\star$ denotes the Arens extension of the
    product in $C(K)$ to $B(K)$, see the proof below.}
    \]
    holds for all $x,y \in A$ with $\|x\|_\infty ,\|y\|_\infty \le
    1$.
\end{lem}
\begin{proof}
    Denote by $\star$ the product in $C(K)$. Since $C(K)$ with the
    product $\star$ is a Banach \falg, there exists a weight $w \in
    C(K)_+$, $\|w\|_\infty \le 1$, such that $(f\star
    g)(t)=w(t)f(t)g(t)$ for all $f,g \in C(K)$ (see \cite[Korollar
    1.4]{scheffold1981}).

    Consider the bidual $C(K)^{* *}$
    equipped with the Arens product. This is again a Banach \falg\
    (see \cite{scheffold1991}), and
    we have a lattice-algebra isometry $C(K)\to C(K)^{* *}$.
    Consider the Banach lattice $B(K)$ of bounded Borel-measurable
    functions on $K$.
    Note that $C(K)\subseteq B(K)\subseteq C(K)^{* *}$ in a natural way,
    since every bounded Borel-measurable function can be integrated against an
    element of $M(K)=C(K)^{*}$. We are going to show that $B(K)$ is a
    subalgebra of $C(K)^{* *}$, thus becoming a Banach \falg. Even
    better, the Arens product of $F, G \in B(K)$ is $(F\star
    G)(t)=w(t)F(t)G(t)$ because, if $f,g \in C(K)$ and $\mu \in M(K)$,
    then
    \[
        \left< \mu f,g \right> = \left< \mu ,f\star g \right> =
        \int_{K}^{} g(t) w(t) f(t) d \mu (t)
    \]
    implies $d(\mu f)(t)=w(t) f(t) d \mu (t)$, and
    \[
    \left< G \mu ,f \right> = \left< G, \mu f \right> = \int_{K}^{}
        G(t) d(\mu f)(t)=\int_{K}^{} f(t)w(t)G(t)d \mu (t)
    \]
    implies $d(G \mu )(t)=w(t)G(t)d \mu (t)$. Therefore
    \[
    \left< F\star G, \mu  \right> = \left< F,G \mu \right>
    =\int_{K}^{} F(t)w(t)G(t)d \mu (t)
    \]
    and, since this holds for all $\mu  \in M(K)$, it follows that
    $(F\star G)(t)=w(t)F(t)G(t)$ for all $t \in K$.

    Fix $\varepsilon >0$. For convenience, we shall assume that
    $\|f_1\|_\infty ,\ldots ,\|f_n\|_\infty \le 1$; this can always be
    achieved by changing $\Phi $ if necessary. Let $0=c_0<c_1<\cdots
    <c_N<c_{N+1}=1+\delta$ be
    a partition of the interval $[0,1+\delta ]$ with $\max_i
    |c_{i+1}-c_i|<\delta  $, this $\delta  >0$ small enough and to be
    determined later. (We use $1+\delta $ instead of $1$ to force
    $c_N>1$, since our functions can take the value $1$.)
    Consider the algebra $\mathcal{C}$ of subsets of $K$ generated by
    \[
        \{c_i\le (f _s)_\sigma < c_{i+1}\}\text{ for }s =1,\ldots
        ,n,
        i=0,\ldots ,N, \sigma \in \{+, -\}\text{ and }
    \]
    \[
        \{c_i\le w < c_{i+1}\}\text{ for }i=0,\ldots ,N.
    \]
    Note that this algebra is contained in the Borel sets of
    $K$. Also, since it is finitely generated, $\mathcal{C}$ is
    actually finite. And as a finite algebra, every set can be
    written as the union of the atoms inside that set. Denote by
    $\{a_1,\ldots, a_l\}\subseteq B(K)$ the characteristic functions of
    the atoms in $\mathcal{C}$.

    This observation implies that
    \[
    \chi _{\{c_i\le (f _s)_\sigma < c_{i+1}\}} = \sum_{j \in
    J(s,\sigma ,i)} a_j,
    \]
    for certain $J(s,\sigma ,i)\subseteq \{1,\ldots ,l\}$, and
    similarly one defines $J(w,i)$ so that
    \[
    \chi _{\{c_i\le w < c_{i+1}\}} = \sum_{j \in
    J(w,i)} a_j.
    \]
    Since $0\le (f_s)_\sigma ,w\le
    \one$,
    we have
    \[
        (f _s)_\sigma =\sum_{i=0}^{N} (f _s)_\sigma \sum_{j
        \in J(s,\sigma ,i)}^{} a_j, \quad w=\sum_{i=0}^{N} w \sum_{j
        \in J(w,i)}^{}a_j.
    \]

    Let $A=\barespn\{a_1,\ldots ,a_l\}$. We are going to define
    elements of $A$ that approximate $(f_s)_\sigma $ and $w$. For
    every $s \in \{1,\ldots ,n\}$, $\sigma \in \{+,-\}$, and $j \in \{1,\ldots
    ,l\}$, define $c_{i(s,\sigma,j)}\in \{c_0,\ldots ,c_{N}\}$ to be
    the unique element of the set that satisfies
    \[
        a_j \le \chi _{\{c_{i(s,\sigma,j)}\le (f_s)_\sigma
        <c_{i(s,\sigma ,j)+1}\}}.
    \]
    To make the notation shorter we will write $c_{i(s,\sigma ,j)}$ as
    $c_{i(j)}^{s,\sigma }$ most of the time. Define
    \[
    (f _s)_{\sigma ,d}=\sum_{j=1}^{l}c_{i(j)}^{s,\sigma
        }a_j.
    \]
    In other words, $(f_s)_{\sigma ,d}$ is a discrete version of
    $(f_s)_\sigma $ that is constant on the sets corresponding to the
    characteristic functions $a_1,\ldots ,a_l$. More precisely,
    $(f_s)_{\sigma ,d}$ is identically $c_{i(s,\sigma ,j)}$ on the
    support of $a_j$, where $c_{i(s,\sigma ,j)}\le (f_s)_\sigma
    <c_{i(s,\sigma ,j)+1}$ holds. In particular, $0\le (f_s)_{\sigma
    ,d}\le (f_s)_{\sigma }$.

    Similarly, for $j \in \{1,\ldots ,l\}$ define $c_{t(j)}\in
    \{c_1,\ldots ,c_N\}$ to be the unique element of the set that
    satisfies
    \[
        a_j \le \chi _{\{ c_{t(j)}\le w<c_{t(j)+1}\}}
    \]
    unless $a_j\le \chi _{\{0\le w<c_1\}}$, in which case
    $c_{t(j)}=c_1$. This way we prevent the discrete version of $w$
    from vanishing at any point. With this
    change it is still true that
    $|c_{t(j)}-w|<\delta $ holds in the support of $a_j$ for all
    $j=1,\ldots ,l$. Define
    \[
        w_d=\sum_{j=1}^{l} c_{t(j)} a_j.
    \]

    Since the $a_i$ are disjoint, $A$ is a sublattice of $B(K)$.
    Define a product $\circ$ on $A$ as follows: $a_i\circ
    a_j=0$ if $i\neq j$, and $a_j\circ a_j =c_{t(j)}a_j$. It is
    immediate that this extends by linearity to a semiprime \falg\ product on $A$.

    Fix $\varepsilon '>0$. We are going to show that, for $\delta
    <\varepsilon '$,
    the product in $A$ satisfies $|x\circ y|\le
    |x\star y|+\varepsilon '\one$ for all $x,y \in A$ with $\|x\|_\infty
    ,\|y\|_\infty \le 1$. Let
    $x=\sum_{j=1}^{l}\lambda _ja_j$ and $y=\sum_{j=1}^{l}\mu _j a_j$
    for some $\lambda _1,\ldots ,\lambda _l,\mu _1,\ldots ,\mu _l \in
    [-1,1]$. Then
    \begin{align*}
        |x\circ y|&=\bigg| \sum_{j=1}^{l}\lambda _j\mu
        _jc_{t(j)}a_j\bigg|\\&\le
    \bigg|\sum_{j=1}^{l}\lambda _j \mu _j w a_j\bigg| + \sum_{j=1}^{l}|\lambda
        _j| |\mu _j| |wa_j-c_{t(j)}a_j|\\&\le |x\star y|+\one (\sup_{j=1,\ldots
    ,l}\sup_{\supp(a_j)} |w-c_{t(j)}|)\\
                                         &\le |x\star y| + \delta \one.
    \end{align*}
    Let $k_1,\ldots ,k_n\ge 0$ be such that $k=k_1+ \cdots +k_n\ge 1$.
    By definition:
    \begin{align*}
        (f  _1)_{\sigma _1}^{k_1}\star \cdots \star (f
        _n)_{\sigma _n}^{k_n}&=\sum_{j=1}^{l}(f
        _1)_{\sigma _1}^{k_1}\cdots (f _n)_{\sigma _n}^{k_n}w^{k-1}a_j\\
        (f  _1)_{\sigma _1,d}^{k_1}\circ \cdots \circ (f
        _n)_{\sigma _n,d}^{k_n}&=\sum_{j=1}^{l} (c_{i(j)}^{1,\sigma
        _1})^{k_1}\cdots (c_{i(j)}^{n,\sigma _n})^{k_n} a_j\circ
        \cdots \circ a_j\\
                               &=\sum_{j=1}^{l} (c_{i(j)}^{1,\sigma
        _1})^{k_1}\cdots (c_{i(j)}^{n,\sigma _n})^{k_n}
        (c_{t(j)})^{k-1}a_j
    \end{align*}
    so that
    \begin{multline*}
        \|(f  _1)_{\sigma _1}^{k_1}\star \cdots \star (f
        _n)_{\sigma _n}^{k_n}-(f  _1)_{\sigma _1,d}^{k_1}\circ \cdots \circ (f
        _n)_{\sigma _n,d}^{k_n}\|_\infty =\\
        \sup_{j=1,\ldots ,l}\sup_{\supp a_j} | (f
        _1)_{\sigma _1}^{k_1}\cdots (f _n)_{\sigma
        _n}^{k_n}w^{k-1} - (c_{i(j)}^{1,\sigma
        _1})^{k_1}\cdots (c_{i(j)}^{n,\sigma _n})^{k_n}
        (c_{t(j)})^{k-1}|.
    \end{multline*}
    But in $\supp a_j$, $|(f _s)_{\sigma
    _s}-c_{i(j)}^{s,\sigma _s}|<\delta $ and $|w-c_{t(j)}|<\delta $.
    Using the fact that the function
    \[
    \begin{array}{cccc}
    & [0,1]^{n+1} & \longrightarrow & \R \\
            & (x_1,\ldots ,x_n,y) & \longmapsto & x_1^{k_1}\cdots
            x_n^{k_n}y^{k-1} \\
    \end{array}
    \]
    is uniformly continuous, we can make the distance
    \begin{equation}\label{eq:monomials}
    \|(f  _1)_{\sigma _1}^{k_1}\star \cdots \star (f
        _n)_{\sigma _n}^{k_n}-(f  _1)_{\sigma _1,d}^{k_1}\circ \cdots \circ (f
        _n)_{\sigma _n,d}^{k_n}\|_\infty
    \end{equation}
    arbitrarily small by choosing an appropriate $\delta >0$.

    According to
    \cref{prop:simplificationLLA}, there exists an LL
    expression $\Psi $ and $c \in \N$ such that
    \[
    \Phi (\lambda _1,\ldots ,\lambda _n)=\Psi ((\lambda  _1)_{\sigma
    _1}^{j_1}\cdots (\lambda  _n)_{\sigma _n}^{j_n}:j_1+ \cdots +j_n\le c, \sigma_i \in \{+,-\})
    \]
    for all $\lambda _1,\ldots ,\lambda _n \in \R$.
    This identity will also hold when evaluating the previous
    expression in arbitrary Banach \falg s (\cref{thm:falgYudin}).
    In particular, in $B(K)$,
    \[
    f=\Psi ((f  _1)_{\sigma
    _1}^{j_1}\star\cdots \star(f  _n)_{\sigma _n}^{j_n}:j_1+ \cdots +j_n\le c,
    \sigma_i \in \{+,-\})
    \]
    and, in $A$,
    \[
    \Phi ((f_1)_d,\ldots ,(f_n)_d)=\Psi ((f  _1)_{\sigma
    _1,d}^{j_1}\circ\cdots \circ(f _n)_{\sigma _n,d}^{j_n}:j_1+ \cdots
    +j_n\le c, \sigma_i \in \{+,-\}),
    \]
    where we have defined $(f_s)_d=(f_s)_{+,d}-(f_s)_{-,d}$. Note that, since
    $(f_s)_{+,d}\wedge (f_s)_{-,d}\le (f_s)_+\wedge (f_s)_{-}=0$,
    $(f_s)_{+,d}$ and $(f_s)_{-,d}$ are the positive and negative
    parts of $(f_s)_d$, respectively. In particular,
    \[
        |(f_s)_d|=(f_s)_{+,d}+(f_s)_{-,d}\le (f_s)_++(f_s)_-=|f_s|.
    \]
    From the norm continuity of lattice-linear function calculus, and the fact
    that the quantities \eqref{eq:monomials} can be made arbitrarily
    small for all $k_1+ \cdots +k_n\le c$ by choosing $\delta >0$
    appropriately (keep in mind there are only finitely many
    non-negative integers
    $k_1,\ldots ,k_n$ satisfying this condition), it follows that
    \begin{equation*}
    \|f-\Phi ((f_1)_d,\ldots ,(f_n)_d)\|_\infty <\varepsilon.\qedhere
    \end{equation*}
\end{proof}

\begin{proof}[Proof of \cref{thm:norm_fin_dim}]
    Let $\{e_1,\ldots ,e_n\}$ be the canonical basis of $\el 1^{n}$.
    According to \cref{thm:finite}, $\FBFA{\el 1^{n}}$ can be isomorphically identified with
    $C(K)$ equipped with an appropriate product $\star$, where $K=[0,1]\times
    \Sinfn$.
    With this identification, the free
    norm $\|{\cdot }\|$ is equivalent to the uniform norm $\|{\cdot
    }\|_\infty$. Since $\FNFA {\el 1^{n}}$ is dense in $\FBFA {\el
    1^{n}}$, it suffices to show the result for $g
    \in \FBFA {\el 1^{n}}$ living in $\FNFA {\el 1^{n}}$, the vector lattice algebra generated by
    $\eta _1=\eta (e_1),\ldots ,\eta _n=\eta (e_n)$. Indeed, suppose the
    result holds in $\FNFA {\el 1^{n}}$, and let $g \in \FBFA {\el
    1^{n}}$. Let
    $\varepsilon >0$ and let $f \in \FNFA {\el 1^{n}}$ be such that
    $\|g-f\|<\varepsilon /3$. Let $T\colon \el 1^{n}\to A$ be a contractive
    operator, where $A$ is a semiprime finite-dimensional Banach
    \falg, such that $\|f\|\le
    \|\hat{T}f\|+\varepsilon /3$. Then
    \[
    \|\hat{T}f\|\le \|\hat{T}(f-g)\|+\|\hat{T}g\|\le \varepsilon
    /3+\|\hat{T}g\|
    \]
    and therefore
    \[
    \|\hat{T}g\|\le \|g\|\le \|f\|+\varepsilon /3\le
    \|\hat{T}f\|+2/3\varepsilon \le \varepsilon +\|\hat{T}g\|.
    \]
    Since $\varepsilon >0$ is arbitrary, the result also follows for
    $g$.

    So assume $f \in \FNFA {\el 1^{n}}$. There exists an LLA
    expression $\Phi $ such that $f =\Phi (\eta_1,\ldots ,\eta_n)$.
    Fix $\varepsilon >0$. Then there exists $\varepsilon '>0$ small
    enough such that, if $f_i=\eta _i/(1+\varepsilon ')$, then
    \[
    \|f-\Phi (f_1,\ldots ,f_n)\|<\varepsilon.
    \]
    Without loss of generality, we can assume $\varepsilon
    '<\varepsilon $; this will be useful later.
    Put $f'=\Phi (f_1,\ldots ,f_n)$ and apply the previous lemma to $f'$:
    there exist a finite-dimensional
    sublattice $A$ of $B(K)$, together with a semiprime \falg\ product
    $\circ$ in $A$, and $(f_1)_d,\ldots ,(f_n)_d \in A$ satisfying
    $|(f_i)_d|\le |f_i|$ such that
    \[
    \|f' - (f')_d\|_\infty <\varepsilon \quad\text{where }(f')_d=\Phi ((f_1)_d,\ldots ,(f_n)_d)
    \]
    (the latter expression being evaluated inside $A$). Moreover, we can
    choose $\circ $ so that $|x\circ y|\le |x\star y|+(\varepsilon
    '/\|\one\|)\one$ holds for all $x,y \in A$ with $\|x\|_\infty
    ,\|y\|_\infty \le 1$.

    Note that the biduals of $C(K)$ with respect to either the free or
    the supremum norm coincide up to equivalence of norms, with the same equivalence constants. We shall denote
    the norms in the bidual the same way as the original ones. Therefore, making $\varepsilon
    >0$ small as necessary, we can assume that
    $\|f'-(f')_d\|<\varepsilon $.

    The sublattice $A$ is a
    Banach lattice with respect to both the free and supremum norms.
    For every $x,y \in A$ with $\|x\|
    ,\|y\|\le 1$, in particular we have $\|x\|_\infty ,\|y\|_\infty
    \le 1$, and therefore the inequality $|x\circ y|\le |x\star y| +
    (\varepsilon '/\|\one\|)\one$ holds. Taking free norms in this inequality:
    \[
        \|x\circ y\|\le \|x\star y\| + \varepsilon '\le 1+\varepsilon
        ',
    \]
    where in the second inequality we are using that the free norm is
    submultiplicative.
    Hence $\|x\circ y\|\le (1+\varepsilon ')\|x\| \|y\|$ holds for
    general $x,y \in A$. Define on $A$ the norm $\trinor{\cdot
    }=(1+\varepsilon ')\|{\cdot }\|$. This is certainly a lattice
    norm, and it is submultiplicative by the previous inequality.
    Hence $A$ is a semiprime Banach \falg\ with respect to this norm.

    Define a linear map $T\colon \el 1^{n}\to A$ by $Te_i=(f _i)_d$
    for $i=1,\ldots ,n$. Since
    \[
        \trinor{(f_i)_d}=(1+\varepsilon ')\|(f_i)_d\|\le
        (1+\varepsilon ')\|f_i\|\le 1
    \]
    this map is contractive.
    Therefore it extends to a unique contractive
    lattice-algebra homomorphism $\hat{T}\colon \FBFA {\el 1^{n}}\to
    A$ satisfying
    \[
    \hat{T}f=\Phi ( (f_1)_d,\ldots ,(f_n)_d)=(f')_d.
    \]
    Putting all the estimations together:
    \begin{align*}
        \trinor{\hat{T}f}&=(1+\varepsilon ') \|(f')_d\|\\
                         &\le(1+\varepsilon ')(\varepsilon +\|f'\|)\\
                         &\le(1+\varepsilon )(2\varepsilon +\|f\|).
    \end{align*}
    where in the last inequality we are also using that $\varepsilon
    '<\varepsilon $. Since $\varepsilon >0$ was arbitrary, $\trinor{\hat{T}f}$ can be made arbitrarily close to
    $\|f\|$ by choosing appropriate $A$ and $T$. This shows that $\tau (f)=\|f\|$.
\end{proof}

Next we want to go one step further and show that $\tau $ coincides
with the free norm in $\FBFA {L_1(\mu )}$ for an arbitrary measure
$\mu $.

When $T\colon E\to F$ is a contractive operator between Banach spaces $E$ and
$F$, denote by $\bar{T}\colon \FBFA E\to \FBFA F$ the unique
lattice-algebra homomorphism with $\|\bar{T}\|=\|T\|$ making the
diagram
\[
\begin{tikzcd}
    \FBFA E\arrow[r, "\bar{T}"]&\FBFA F\\
    E\arrow[u, "\eta _E"]\arrow[r, "T"]&F\arrow[u, "\eta _F"]
\end{tikzcd}
\]
commute. Note that $\bar{T}=\widehat{\eta _F T}$. The following is a
standard result relating $T$ and $\bar{T}$.

\begin{prop}\label{prop:complemented}
    Let $E$ be a Banach space, let $F$ be a subspace of $E$ with a
    contractive projection map
    $P\colon E\to F$, and
    let $\iota \colon F\to E$ be the inclusion map. Then
    $\bar{P}\colon \FBFA E\to \bar\iota (\FBFA F)$ is a contractive projection. In
    particular,
    $\bar{\iota}\colon \FBFA F\to \FBFA E$ is an isometric
    embedding.
\end{prop}
\begin{proof}
    Since $P \iota =I_F$, uniqueness implies $\bar{P}\bar{\iota
    }=I_{\FBFA F}$, so that $\bar
    \iota $ is injective and $\bar{P}$ is onto $\bar \iota (\FBFA
    F)$. From $P (\iota P)=P$ it follows that $\bar{P}$ is a
    projection, and since $\|\bar \iota\|, \|\bar{P}\|\le 1$,
    $\bar\iota $ is an isometric embedding.
\end{proof}

\begin{prop}\label{prop:tau_is_norm}
    Let $E$ be a Banach space. Suppose
    there exists a net $\{E_\lambda \}_\lambda $ of
    subspaces of $E$ such that:
    \begin{enumerate}
        \item the closure of $\bigcup_{\lambda}E_\lambda $ is $E$,
        \item there exist contractive projections $P_\lambda \colon
        E\to E_\lambda $,
        \item the free norm in $\FBFA{E_\lambda }$ coincides with $\tau
        _{E_\lambda }$.
    \end{enumerate}
    Then $\tau _{E}$ coincides with the free norm in $\FBFA{E}$.
\end{prop}
\begin{proof}
    The result will follow once we show that $\tau (f)\ge \|f\|$ for
    every $f \in \FBFA E$. Fix $\varepsilon >0$ and let $g \in \FNFA
    E$ be such that $\|f-g\|<\varepsilon/2$. There exists an LLA
    expression $\Phi $ and $x_1,\ldots ,x_n \in E$ satisfying $g=\Phi
    (\eta _{x_1},\ldots ,\eta _{x_n})$. Choose $z_1,\ldots ,z_n \in
    \bigcup_{\lambda } E_\lambda $ with $\|z_i-x_i\|<\delta $, this
    $\delta >0$ to be determined later. Since $\{E_\lambda
    \}_\lambda $ is a net, there exists some $G =E_{\lambda _0}$ for
    which $z_1,\ldots ,z_n \in G$.
    Let $h=\Phi (\eta _{z_1},\ldots ,\eta _{z_n})$. Since the
    operations in $\FBFA E$ are norm continuous,
    $\|g-h\|<\varepsilon/2 $ as long as we choose $\delta >0$ small enough. Hence
    $\|f-h\|<\varepsilon$.

    According to \cref{prop:complemented},
    one can view $\FBFA G$ as a closed sublattice-algebra of $\FBFA
    E$; with this identification, $h \in \FBFA G$, and $\|h\|_{\FBFA
    E}= \|h\|_{\FBFA G}$. By
    assumption, there exists a semiprime finite-dimensional Banach
    \falg\ $A$, and a contractive operator
    $T\colon G\to A$, such that $\|\hat{T}h\|_A\ge\|h\|
    -\varepsilon$. Then
    $TP_{\lambda_0}\colon E\to A$ is a contractive operator such that
    \[
        \widehat{TP_{\lambda_0} }h=\Phi (TP_{\lambda_0} z_1,\ldots
        ,TP_{\lambda_0}
    z_n)=\Phi (Tz_1,\ldots ,Tz_n)=\hat{T}h
    \]
    and
    \[
        \|\hat{T}h\|_A=\|\widehat{TP_{\lambda_0} }h\|_A\le \|f-h\|+\|\widehat{TP_{\lambda_0}
        }f\|_A<\varepsilon +\|\widehat{TP_{\lambda_0} }f\|_A
    \]
    so that
    \[
    \tau (f)\ge \|\widehat{TP_{\lambda_0} }f\|_A> \|\hat{T}h\|_A-\varepsilon\ge
    \|h\|-2\varepsilon \ge \|f\|_{\FBFA E}-3\varepsilon.
    \]
    Since $\varepsilon >0$ was arbitrary, this shows that $\tau (f)=
    \|f\|$.
\end{proof}

\begin{rem}
    Suppose that $F$ is a contractively complemented subspace of $E$ and
    that $\tau _E$ coincides with the free norm in $\FBFA E$. We are
    going to show that $\tau _F$ also coincides with the free norm in
    $\FBFA F$. Let $P\colon
    E\to F$ be a contractive projection. By
    \cref{prop:complemented}, the norms $\|{\cdot }\|_{\FBFA {F}}$
    and $\|{\cdot }\|_{\FBFA E}$ coincide on $\FBFA{F}$ when the
    latter is viewed as a sublattice-algebra of $\FBFA E$. So it suffices to
    check that $\tau _{F}(f)\ge \tau _E(f)$ for every
    $f \in \FBFA {F}$. Since $\tau $ is continuous with respect to the
    free norm (this is easy to see using that $\tau $ is a seminorm
    and that $\tau \le \|{\cdot }\|$), it suffices to
    check the inequality for $f \in \FNFA {F}$. Let $A$ be a semiprime
    finite-dimensional Banach \falg, and let $T\colon E\to A$
    be a contractive operator. Its restriction $S\colon F\to A$ is
    also a contractive operator such that $\hat{S}f=\hat{T}f$ for all
    $f \in \FNFA {F}$. Therefore $\|\hat{T}f\|=\|\hat{S}f\|\le \tau
    _{F}(f)$. Since $T$ and $A$ were arbitrary, $\tau _E(f)\le \tau
    _{F}(f)$.
\end{rem}

\begin{cor}\label{cor:tau_is_norm_example}
    The norm in $\FBFA{L_1(\mu )}$ is $\tau $ for every measure $\mu $.
\end{cor}
\begin{proof}
    In $L_1(\mu )$ there exists a net of subspaces, each isometric to
    $\el 1^{n}$ for some $n \in \N$, satisfying the conditions in
    \cref{prop:tau_is_norm}. The result then follows from
    \cref{thm:norm_fin_dim}.
\end{proof}

\chapter[Representation of Free Banach \texorpdfstring{$f$-algebras}{f-algebras}]{Representation of Free Banach
\texorpdfstring{$f$-algebras}{f-algebras} in Spaces of Continuous
Functions}
\label{sec:representation}

Recall from \cref{cor:iota} that, for every Banach space $E$, we have an
injective and contractive lattice-algebra homomorphism $\iota \colon
\FNFA E\to C(B_{E^{*}})$. This map can be extended to the completion
$\hat{\iota }\colon \FBFA E\to C(B_{E^{*}})$. The extension is
certainly contractive and a lattice-algebra homomorphism. But is it
injective? When it is, we will say that
$\FBFA E$ is \emph{representable in $C(B_{E^{*}})$}. Since $\iota $ is
never an embedding (\cref{rem:notembed}), it is not clear at all
whether $\FBFA E$ is representable in $C(B_{E^{*}})$ or not.

Compare this with the situation for the free vector lattice $\FVLv
E$, whose restriction map $\FVLv E\to C(B_{E^{*}})$ extends to an
injective and contractive lattice-algebra homomorphism $\FBL E\to
C(B_{E^{*}})$ (see \cref{sec:background_FBL}). Once $\FBL E$ has
been represented inside $C(B_{E^{*}})$, many of its properties
follow, as they depend only on this representation rather than on
the explicit expression of the norm (see \cite{oikhberg_etal2022}).
Therefore, even though no explicit expression for the norm in $\FBFA
E$ is available, we will be able to recover several properties
analogous to those of $\FBL E$ as long as we assume that $\FBFA E$ is
representable in $C(B_{E^{*}})$ (this will be done in
\cref{sec:properties}).

In this chapter we explore when $\FBFA E$ is representable. The
space $\FBFA E$ is quite difficult to understand, as it is the
completion of $\FNFA E$ with respect to a certain abstract norm.
First, we characterize those $E$ for which $\FBFA E$ is
representable in $C(B_{E^{*}})$ as those for which $\FBFA E$ is
semiprime. Representability is thus reduced to a purely algebraic
question. However, whether or not $\FBFA E$ is semiprime is still
not obvious, since the completion of a semiprime normed \falg\ need
not be semiprime. At this point, the technical computations
involving the seminorm $\tau $ carried out in the previous chapter
come to rescue: we show that, if $\tau $ defines a norm on $\FBFA
E$, then $\FBFA E$ is semiprime.

\begin{thm}\label{thm:semiprime_iff_representable}
    Let $E$ be a Banach space. The map $\iota \colon \FNFA E\to
    C(B_{E^{*}})$ extends to an injective lattice-algebra homomorphism
    $\hat{\iota }\colon \FBFA E\to C(B_{E^{*}})$ if and only if $\FBFA E$ is semiprime.
\end{thm}
\begin{proof}
    Suppose first that the map $\iota \colon \FNFA E\to
    C(B_{E^{*}})$ extends to an injective
    lattice-algebra homomorphism $\hat{\iota }\colon \FBFA E\to
    C(B_{E^{*}})$. Let $f \in N(\FBFA E)$. For every $x^{*} \in
    B_{E^{*}}$, $x^{*}\neq 0$, there exists some $x \in B_{E}$ for
    which $\iota (\eta _x)(x^{*})=x^{*}(x)\neq 0$. Since $f\eta _x=0$, it
    follows that
    \[
    0=\hat{\iota }(f \eta _x)(x^{*})=\hat{\iota }(f)(x^{*})
    \iota (\eta _x)(x^{*}).
    \]
    Thus $\hat{\iota }(f)(x^{*})=0$ for every $x^{*}\in
    B_{E^{*}}\setminus \{0\}$, and by continuity $\hat{\iota }(f)(0)=0$. So $\hat{\iota } (f)=0$ and, since
    we are assuming $\hat{\iota }$ to be injective, $f=0$.

    Conversely, suppose $\FBFA E$ is semiprime. The left
    regular representation $L\colon \FBFA E\to \mathcal{Z}(\FBFA E)$,
    which sends $f \in \FBFA E$ to the left multiplication operator
    $L_f(g)=fg$, $g \in \FBFA E$, is a
    contractive lattice-algebra homomorphism. Since $\FBFA E$ is
    semiprime, $L$ is also injective. Let $Z$ be the closure
    of $\VLA(\{\, L_f : f \in \FBFA E \, \}\cup \{I\})$ in
    $\mathcal{Z}(\FBFA E)$ (that is, the closure with respect to the
    uniform norm that induces $I$, which coincides with the operator
    norm in $\mathcal{Z}(\FBFA E)$, see \cite{wickstead1977}). Being a Banach \falg\ with
    identity and having the uniform norm with respect to this
    identity, $Z$ is lattice-algebra isometric to $C(K)$ for a certain
    compact Hausdorff space $K$ (\cite{martignon1980}).
    Restricting the range of $L$ yields a contractive and injective
    lattice-algebra homomorphism $L\colon \FBFA E\to C(K)$.

    Along this proof, $\phi _t$ for $t \in K$ shall denote the
    point-mass measure at $t$. The composition map $\phi _tL\colon
    \FBFA E\to \R$ is a contractive lattice-algebra homomorphism.
    Hence it is uniquely determined by $\phi _tL\eta \colon E\to \R$,
    which is an element of $B_{E^{*}}$. Define the map
    \[
    \begin{array}{cccc}
    \psi \colon& K & \longrightarrow & (B_{E^{*}},w^{*}) \\
            & t & \longmapsto & \phi _tL\eta  \\
    \end{array}.
    \]
    Note that $\psi (t)=(L \eta
    )^{*}\phi _t$. If
    we identify $K$ homeomorphically with $\{\, \phi _t : t \in K \,
    \} $, then the map $\psi $ is nothing but the restriction of
    $(L\eta )^{*}$ to $\{\, \phi _t : t \in K \, \} $. Since $L \eta $
    is continuous, its adjoint $(L\eta )^{*}$ is weak$^*$ to weak$^*$
    continuous. Hence $\psi $ is continuous.
    It is also injective: if $\psi (t)=\psi (s)$, then for
    all $x \in E$, $\phi _t(L_{\eta _x})=\phi _s(L_{\eta _x})$, and
    since the point-mass measures are lattice-algebra homomorphisms,
    $\phi _t$ and $\phi _s$ coincide on $L(\FNFA E)$. Since
    $L(\FNFA E)$ is dense in $L(\FBFA E)$, by continuity they must
    also coincide on $L(\FBFA E)$. Also, $\phi _t(I)=1=\phi
    _s(I)$. Since $L(\FBFA E)$ and $I$ generate, as a Banach lattice
    algebra, the whole $C(K)$, it follows that $\phi _t=\phi _s$ and
    therefore that $t=s$.

    Being $K$ compact, $\psi $ is actually an embedding. That is, $K$
    is homeomorphic to its image $\psi (K)\subseteq B_{E^{*}}$.
    Note that $(L \eta )(x)=\iota (\eta _x)\circ \psi $. In other
    words: the map
    \[
    \begin{array}{cccc}
    S\colon&  E & \longrightarrow & C(\psi (K)) \\
            & x & \longmapsto & \iota (\eta _x)|_{\psi (K)} \\
    \end{array}
    \]
    extends to a lattice-algebra homomorphism $\hat{S}\colon \FBFA
    E\to C(\psi (K))$ which is injective. Indeed, $C_{\psi ^{-1}}\circ L$, where $C_{\psi
    ^{-1}}\colon C(K)\to C(\psi (K))$ denotes the composition
    map, is an injective and contractive lattice-algebra homomorphism satisfying
    $(C_{\psi ^{-1}}\circ L)\eta(x)=\iota (\eta _x)|_{\psi (K)}$. By
    uniqueness, $\hat{S}$ is precisely $C_{\psi ^{-1}}\circ L$.

    Consider the natural restriction map $R\colon C(B_{E^{*}})\to
    C(\psi (K))$, which is also a contractive lattice-algebra homomorphism. The diagram
    \[
    \begin{tikzcd}[column sep=small]
        \FBFA E \arrow[rr, "\hat{S}"]\arrow[dr, "\hat{\iota }"']&
                                                            &C(\psi
                                                            (K))\\
        &C(B_{E^{*}})\arrow[ru, "R"']&
    \end{tikzcd}
    \]
    commutes, because $\hat{S} (\eta _x)=\iota (\eta _x)|_{\psi (K)}=(R
    \hat{\iota })(\eta
    _x)$ and all arrows are contractive lattice-algebra homomorphisms. Since
    $\hat{S}$ has been shown to be injective, so must be $\hat{\iota }$.
\end{proof}

\begin{cor}\label{cor:finite_representable}
    If $E$ is a finite-dimensional Banach space, then $\FBFA E$ is
    representable.
\end{cor}
\begin{proof}
    This is a consequence of \cref{cor:finite_semiprime}.
\end{proof}

Even though $\FNFA E$ is semiprime for every Banach space $E$, it
is not clear at all whether its completion $\FBFA E$ is also
semiprime. In general, it is not true that the completion of a
semiprime normed \falg\ is semiprime. The following lemma provides a
general procedure to construct examples of such algebras.

\begin{lem}\label{lem:comp_falg}
    Let $B$ be a semiprime Banach \falg. Let $A\subseteq \R_0\oplus_\infty B$
    be a sublattice-algebra such that $(1,0)\not\in A$ but $(1,0)\in
    \bar{A}$. Then $A$ is a semiprime normed \falg\ whose completion
    is not semiprime.
\end{lem}
\begin{proof}
    Since $(1,0) \in \bar{A}$, and $(1,0)^2=(0,0)$, $\bar{A}$ is not
    semiprime. To show that $A$ is semiprime, let $(\lambda,a)\in A$
    be such that $(\lambda,a )^2=(0,0)$. Then $a^2=0$ and, since $B$
    is semiprime, $a=0$. Thus $(\lambda ,a)=\lambda (1,0) \in A$, but
    we are assuming $(1,0) \not\in A$, so that $\lambda =0$.
\end{proof}

Below we present two explicit examples of this construction. The first
uses a limiting process to guarantee that $(1,0)\not\in A$, while the
second uses the freeness of the generators.

\begin{example}\label{ex:comp_falg_simple}
    Let
    \[
    A=\{\, (\lambda ,(x_n)) \in \R_0\oplus c : \lim nx_n=\lambda \, \}.
    \]
    It is clear that this is a sublattice of $\R_0\oplus c$. It is
    also a subalgebra, for if $(\lambda ,(x_n)),(\mu ,(y_n)) \in A$,
    then their product is $(0,(x_ny_n))$ and, since $\lim x_n=0$,
    $\lim x_ny_n=0$. Clearly $(1,0) \not\in A$. For every $n \in \N$,
    the sequences
    \[
    x^{(n)}=(0,\ldots ,0,\frac{1}{n+1},\frac{1}{n+2},\ldots ),
    \]
    where the zeros are in the first $n$ coordinates,
    are such that $x^{(n)}\to 0$ in $c$. Moreover, $(1,x^{(n)}) \in A$.
    Therefore $\lim (1,x^{(n)})=(1,0) \in \bar{A}$. By
    \cref{lem:comp_falg}, $A$ is a semiprime normed \falg\ whose
    completion is not semiprime.
\end{example}

\begin{example}\label{ex:comp_falg}
    Consider the Banach \falg\ $B=\R_0\oplus_\infty \FBFA{\el 1}$.
    Denote by $e_n$, $n\ge 1$, the canonical basis of $\el 1$. Let
    \[
    A=\VLA\{(1,\eta _{e_n/n}): n \in \N\}\subseteq B.
    \]
    It is clear that $A$ is a sublattice algebra. Its completion
    contains $(1,0)$, because $\|\eta _{e_n/n}\|=\|e_n/n\|_1=1/n$ tends
    to zero.

    Let us show that $(1,0) \not\in A$. Suppose, in order to reach a
    contradiction, that $(1,0) \in A$. Then there
    exists an LLA expression $\Phi [t_1,\ldots ,t_n]$ such that
    $(1,0)=\Phi ((1,\eta _{e_1}),\ldots ,(1,\eta _{e_n/n}))$. Since the
    operations are computed coordinatewise:
    \begin{equation}\label{eq:comp_falg}
        (1,0)=(\Phi (1,\ldots ,1),\Phi (\eta _{e_1},\ldots ,\eta
    _{e_n/n})).
    \end{equation}
    We are going to show that $\Phi $ vanishes on
    $[-1/n,1/n]^{n}\subseteq \R^{n}$. Fix $(\lambda _1,\ldots ,\lambda
    _n)\in [-1/n,1/n]^{n}$. The element
    \[
        x^{*}=(\lambda _1,2\lambda _2,\ldots ,n\lambda _n,0,0,\ldots )\in
        B_{\el \infty}
    \]
    seen as a contractive functional on $\el 1$ extends to a
    lattice-algebra homomorphism $\widehat{x^{*}}\colon \FBFA{\el
    1}\to \R$. Therefore
    \[
    0=\widehat{x^{*}}(\Phi (\eta _{e_1},\ldots ,\eta _{e_n/n}))=\Phi
    (x^{*}(e_1),\ldots ,x^{*}(e_n)/n)=\Phi (\lambda _1,\ldots ,\lambda
    _n).
    \]
    This proves that $\Phi $ vanishes on $[-1/n,1/n]^{n}$. By \cref{lem:calc_zero},
    $\Phi $ must also vanish on $\R_0$. This contradicts the fact
    that, according to \eqref{eq:comp_falg}, $1=\Phi (1,\ldots ,1)$ in
    $\R_0$. Hence $(1,0)\not\in A$. By \cref{lem:comp_falg}, $A$ is an
    example of a semiprime normed \falg\ whose completion is not
    semiprime.
\end{example}

The previous example illustrates why determining whether the
completion of a semiprime normed \falg\ is again semiprime or not is far
from trivial. The next proposition shows how $\tau $ can be very
useful in that matter.

\begin{prop}\label{prop:tau_norm_implies_semiprime}
    Let $E$ be a Banach space. Suppose $\tau $ defines a norm on $\FBFA E$. Then $\FBFA E$ is semiprime.
\end{prop}
\begin{proof}
    This proof is based on the observation that every semiprime
    finite-dimensional \falg\ $A$ has an
    identity. Indeed, if $\{e_1,\ldots ,e_n\}$ is a basis of pairwise
    disjoint elements of $A$, then $e_i^2=w_ie_i$ for some $w_i\ge 0$,
    and since $A$ is semiprime, $w_i>0$. The identity is
    the element $\sum_{i=1}^{n} e_i/w_i$.

    Let $f \in N(\FBFA E)$. Let $A$ be a semiprime finite-dimensional Banach
    \falg\ and let $T\colon E\to A$ be a
    contractive operator. Since $\range \hat{T}$ is a
    sublattice-algebra of $A$, in particular is semiprime,
    and by the initial observation has an identity. Let $g$ be in the
    preimage of the identity. Then
    \[
    \hat{T}f=\hat{T}f \hat{T}g=\hat{T}(fg)=0.
    \]
    Since $T$ and $A$ were arbitrary, $\tau   (f)=0$, and since we are
    assuming $\tau $ defines a norm on $\FBFA E$, $f=0$.
\end{proof}

\begin{cor}\label{cor:representable_examples}
    Let $F$ be a contractively complemented subspace of a Banach space $E$. Suppose
    there exists a net $\{E_\lambda \}_\lambda $ of
    subspaces of $E$ such that:
    \begin{enumerate}
        \item the closure of $\bigcup_{\lambda}E_\lambda $ is $E$,
        \item there exist contractive projections $P_\lambda \colon
        E\to E_\lambda $,
        \item the free norm in $\FBFA{E_\lambda }$ coincides with $\tau
        _{E_\lambda }$.
    \end{enumerate}
    Then $\FBFA F$ is representable in $C(B_{F^{*}})$.
\end{cor}
\begin{proof}
    By \cref{prop:tau_is_norm}, $\tau $ is a
    norm on the Banach \falg s generated by these Banach spaces. The
    result follows from a successive application of
    \cref{prop:tau_norm_implies_semiprime} and
    \cref{thm:semiprime_iff_representable}.
\end{proof}

\begin{cor}
    $\FBFA{L_1(\mu )}$ is representable in $C(B_{L_\infty (\mu )})$
    for every measure $\mu $.
\end{cor}

\begin{rem}
    In particular, everything we have
    done so far applies to $\FBFAs S=\FBFA{\el 1(S)}$ for every set
    $S$. Hence, from the point of view of the free Banach \falg\
    generated by a set, we are being completely general.
\end{rem}

After this discussion, we have to leave the following open.

\begin{question}
    Is $\FBFA E$ semiprime for every Banach space $E$?
\end{question}

\chapter[Properties of the free Banach $f$-algebra]{Properties of the Free Banach \texorpdfstring{$f$-algebra}{f-algebra} Generated by a
Banach Space}
\label{sec:properties}

In this chapter we investigate the lattice and algebraic properties of the
free Banach \falg\ generated by a Banach space.

\section{Basic properties}

\begin{prop}
    Let $E$ be a Banach space. Neither $\FNFA E$ nor $\FBFA E$ have an
    approximate identity.
\end{prop}
\begin{proof}
    Suppose $(h_\alpha )$ were an approximate identity in $\FBFA E$.
    Let $x^{*} \in B_{E^{*}}$ be a non-zero functional. Consider it as
    an operator $x^{*}\colon E\to \R_0$, and extend it to a contractive
    lattice-algebra homomorphism
    \[
        \widehat{x^{*}}\colon \FBFA E\to \R_0.
    \]
    Let $x \in E$ be such that $x^{*}(x)\neq 0$. Since $\eta _x
    h_\alpha \to \eta _x$ in $\FBFA E$, it follows
    that
    \[
        \widehat{x^{*}}(\eta _x) \widehat{x^{*}}(h_\alpha )\to
        \widehat{x^{*}}(\eta _x)=x^{*}(x)\neq 0.
    \]
    This is impossible
    because in $\R_0$ the product of any two elements is 0. The same
    argument works for $\FNFA E$.
\end{proof}

The following two properties are analogous to
\cref{prop:FBL_strong_unit}. In this context, however, one has to
be careful with the fact that we do not know whether $\FBFA E$ is
representable in $C(B_{E^{*}})$ or not.

\begin{prop} Let $E$ be a Banach space.
    \begin{enumerate}
        \item $E$ is finite-dimensional if and only if $\FBFA E$ has a
            strong unit.
        \item $E$ is separable if and only if $\FBFA E$ has a quasi-interior point.
    \end{enumerate}
\end{prop}
\begin{proof}
    \begin{enumerate}
        \item When $E$ is finite-dimensional, it was shown in
            \cref{thm:finite} that $\FBFA E$ has a strong unit.
            Conversely, suppose $e \in \FBFA E$ is a strong unit but
            $\dim E=\infty$. Scaling $e$, we may assume that $|\eta _x|\le
            e$ for all $x \in B_E$. Let
            $f \in \FNFA E$ be such that $\|e-f\|<1/2$. Then $\|\hat{\iota }e-\hat{\iota }f\|_\infty <1/2$.
            Since $f \in \FNFA E$, $\hat{\iota }f=\Phi (\iota (\eta
            _{x_1}),\ldots ,\iota (\eta
            _{x_n}))$ for some LLA expression $\Phi $ and some
            $x_1,\ldots ,x_n \in E$. We are assuming $E$ to be
            infinite-dimensional, so there exists $x^{*}\in S_{E^{*}}$ that
            vanishes on $x_1,\ldots ,x_n$. Hence $(\hat{\iota }f)(x^{*})=0$. For every $x \in B_{E}$:
            \begin{align*}
                |x^{*}(x)|&=\iota (|\eta _x|)(x^{*})\\&\le
                \hat{\iota }e(x^{*})\\&=|(\hat{\iota
                }e)(x^{*})-(\hat{\iota} f)(x^{*})|\\&\le
                \|\hat{\iota }e-\hat{\iota }{f}\|_\infty
                <1/2.
            \end{align*}
            This contradicts the fact that $\|x^{*}\|=1$.
        \item Suppose $E$ is separable. Then $\eta (E)$ is
            separable and, since $\FBFA E$ is the closed
            sublattice-algebra generated by $\eta (E)$, it is also
            separable. In particular, it has a quasi-interior point.

            Conversely, suppose $e \in \FBFA E$ is a quasi-interior
            point. If $x^{*}\in B_{E^{*}}$ is such that
            $\hat{\iota }e(x^{*})=0$, then $\hat{\iota }f(x^{*})=0$ for all $f
            \in I_e$. Since $I_e$ is dense in $\FBFA E$,
            it follows that $\hat{\iota }f(x^{*})=0$ for all
            $f \in \FBFA E$. In particular, $x^{*}(x)=\hat{\iota }\eta
            _x(x^{*})=0$ for all $x \in E$. This shows that $\hat{\iota }e$ only
            vanishes at $0$.

            For every $n \in \N$, let $U_n=\{\, x^{*}\in B_{E^{*}} :
            \hat{\iota }e(x^{*})<1/n \, \} $. Then $U_n$ is a weak$^*$ open
            subset of $B_{E^{*}}$ and $\bigcap_{n=1} ^{\infty
            }U_n=\{0\}$. For every $n$, let $A_n\subseteq E$ be a
            finite set such that $U_n$ contains
            \[
                V_n=\{\, x^{*}\in
                B_{E^{*}} : |x^{*}(x)|<1\text{ for all }x \in A_n \,
            \}.
            \]
            Let $F$ be the closed span of $\bigcup_{n=1}^{\infty } A_n
            $. We claim that $F=E$. Indeed, if $F$ were not $E$, then there
            would exist some $x^{*}\in S_{E^{*}}$ vanishing on $F$. On
            the other hand, since $\bigcap_{n=1} ^{\infty }V_n=\{0\}$, there
            exists some $n$ for which $x^{*}\not\in V_n$, that is,
            $|x^{*}(x)|\ge 1$ for some $x \in A_n$. Since $A_n
            \subseteq F$, this is a contradiction.\qedhere
    \end{enumerate}
\end{proof}

\section{Order density}

In the previous section, we insisted that the representability of
$\FBFA E$ in $C(B_{E^{*}})$ (equivalently, the semiprimeness of
$\FBFA E$) is key in showing many properties analogous to those
studied for the free Banach lattice in \cref{sec:background_FBL}.
The upcoming results will show clearly why this is the case. The
next proposition, which is the analogue of \cref{prop:FBL_basic},
follows from a topological argument identical to the proof of
\cref{prop:lat_prop_not_one} once we use the representation
$\hat{\iota }\colon \FBFA E\to C(B_{E^{*}})$.

\begin{prop}
    Let $E$ be a Banach space for which $\FBFA E$ is semiprime.
    \begin{enumerate}
        \item For every $x \in E$, $x\neq 0$, $|\eta _x|$ is a weak
            order unit in $\FBFA E$.
        \item If $E$ has dimension strictly greater than one, then $\FBFA E$
            has no non-trivial projection bands.
        \item If $E$ has dimension strictly greater than one, then $\FBFA E$
            is not $\sigma $-order complete and has no atoms.
    \end{enumerate}
\end{prop}

Even though the free vector lattice is not order dense in
the free Banach lattice unless the base space is
finite-dimensional \cite{azouzi_dhifaoui2025}, the free normed \falg\ is always order dense
in the free Banach \falg\ whenever the latter is representable.
Showing this requires first the following technical fact.

\begin{lem}\label{lem:bump}
    Let $E$ be a Banach space and let $U\subseteq B_{E^{*}}$ be a
    weak$^*$ open set. Then there exists an element $f \in \iota
    (\FNFA E)$, $f>0$,
    such that $\supp f\subseteq U$ and $\|f\|_\infty \le 1$.
\end{lem}
\begin{proof}
    For every $a<b$, the support of the real function $x\mapsto [x_+(x_+ -a)(x_+ -
    b)]_-$ is contained in $[a,b]\cap
    [0,\infty )$. Define, depending on the interval $(a,b)$, the LLA
    expression
    \[B_{(a,b)}[t]=
    \begin{cases}
        [t_+(t_+ - a)(t_+ - b)]_-&\text{if }b>0,\\
        [t_-(t_- + a)(t_- + b)]_-&\text{if }b\le 0.
    \end{cases}
    \]
    Note that $x \mapsto B_{(a,b)}(x)$ defines a non-zero real
    function whose support is $[a,b]\cap [0,\infty )$ if $b>0$ and
    $[a,b]$ otherwise. Normalize it so as to have $0\le
    B_{(a,b)}(x)\le 1$ for all $x \in \R$.

    Making
    $U$ smaller, if necessary, we may assume that there exist $x_1,\ldots
    ,x_n \in B_E$ linearly independent, $\delta >0$ and $x^{*}\in
    E^{*}$ with $\|x^{*}\|<1$ such that
    \[
    U=\{\, y^{*}\in B_{E^{*}} : |y^{*}(x_i)-x^{*}(x_i)|<\delta  \text{
    for }i=1,\ldots ,n\, \}.
    \]
    For a fixed $0<\delta '\le \delta $ to be determined later,
    define
    \[
    f=\bigwedge_{i=1} ^{n} B_{(x^{*}(x_i)-\delta'/2,x^{*}(x_i)+\delta'/2
    )}(\eta _{x_i}).
    \]

    The most difficult part is to show that
    $f\neq 0$. For this, define a
    functional $y^{*}$ on $F=\barespn\{x_1,\ldots ,x_n\}$ in the
    following way. If $|x^{*}(x_i)|\ge \delta '/2$, then set
    $\varepsilon _i=0$ and $y^{*}(x_i)=x^{*}(x_i)+\varepsilon _i$. If instead $-\delta '/2<x^{*}(x_i)<\delta
    '/2$, then there exists $0\le \varepsilon _i<\delta'/2$ such that
    $x^{*}(x_i)+\varepsilon _i> 0$; in this case, set
    $y^{*}(x_i)=x^{*}(x_i)+\varepsilon _i$.

    We are going to show that, for $\delta '$ small enough,
    $\|y^{*}\|\le 1$. Let $C>0$ be such that $\|{\cdot }\|_\infty \le
    C\|{\cdot }\|$ in $F$, where $\|\sum_{i=1}^n \lambda_i x_i\|_\infty=\max\{|\lambda_i|:1\leq i\leq n\}$. For $x \in B_F$, $x=\sum_{i=1}^{n}\lambda
    _ix_i$, compute
    \begin{align*}
        |y^{*}(x)|&=\left| \sum_{i=1}^{n}\lambda
        _i(x^{*}(x_i)+\varepsilon _i) \right|\\
                  &=\left| x^{*}(x)+\sum_{i=1}^{n}\lambda
                  _i\varepsilon _i \right| \\
                  &\le \|x^{*}\|\|x\|+\|x\|_\infty
                  \sum_{i=1}^{n}\varepsilon _i\\
                  &\le \|x\|(\|x^{*}\|+Cn \delta '/2).
    \end{align*}
    Recall that $\|x^{*}\|<1$ was fixed at the beginning, and that $n$
    and $C$ are determined by $x_1,\ldots ,x_n$, which were also chosen
    before $\delta '$. Thus we may set $\delta '$ small enough so as
    to have $\|y^{*}\|\le 1$. Use Hahn--Banach to extend it to a
    functional $y^{*}\in B_{E^{*}}$. This functional satisfies
    \[
    B_{(x^{*}(x_i)-\delta'/2,x^{*}(x_i)+\delta'/2)}(y^{*}(x_i))>0
    \]
    for every $i=1,\ldots ,n$, and so $f(y^{*})>0$. This shows that
    $f\neq 0$. Now we check the remaining properties of $f$. By
    construction, for every $z^{*}\in B_{E^{*}}$, $f(z^{*})\ge 0$, and
    $f(z^{*})\neq 0$ implies $z^{*}(x_i) \in (x^{*}(x_i)-\delta
    /2,x^{*}(x_i)+\delta /2)$ for $i=1,\ldots ,n$. Therefore $f\ge 0$
    and $\supp f \subseteq U$.
\end{proof}

\begin{prop}\label{prop:ord_density}
    Let $E$ be a Banach space. Then $\iota (\FNFA E)$ is order dense in
    $C(B_{E^{*}})$.
\end{prop}
\begin{proof}
    Let $f \in C(B_{E^{*}})$, $f>0$, and let $x^{*}\in B_{E^{*}}$ be
    such that $f(x^{*})\neq 0$. By continuity, there exists an
    $\varepsilon >0$ and a weak$^*$ open neighbourhood $U$ of $x^{*}$
    such that $f(y^{*})>\varepsilon $ whenever $y^{*}\in U$. Let $g
    \in \iota (\FNFA E)$ be an element satisfying the conditions of
    \cref{lem:bump} for the open set $U$. Then $0<\varepsilon g\le f$.
\end{proof}

\begin{cor}\label{cor:ord_dense_FBFA}
    Let $E$ be a Banach space for which $\FBFA E$ is semiprime. Then
    $\FNFA E$ is order dense in $\FBFA E$.
\end{cor}
\begin{proof}
    Let $f \in \FBFA E$, $f>0$. By \cref{thm:semiprime_iff_representable}, we have that $\hat{\iota }\colon \FBFA
    E\to C(B_{E^{*}})$ is injective, so $\hat{\iota } (f)>0$. By
    \cref{prop:ord_density} there exists $g \in \FNFA E$ such that
    $0<\hat{\iota} (g)\le \hat{\iota }(f)$. Again the
    injectivity of $\hat{\iota }$ implies $0<g\le f$.
\end{proof}

\begin{prop}
    For every Banach space $E$, $\FNFA E$ has the countable chain
    condition. If $\FBFA E$ is semiprime, then it also has the
    countable chain condition.
\end{prop}
\begin{proof}
    That $\FNFA E$ has the countable chain condition follows from the
    existence of the representation $\iota _E\colon \FNFA E\to
    C(B_{E^{*}})$ and a topological argument as in \cref{rem:CCC}. If
    $\FBFA E$ is semiprime, then $\FNFA E$ is order dense in $\FBFA
    E$, and therefore $\FBFA E$ also has the countable chain
    condition.
\end{proof}

\section{Extension of operators}

The following is an analogue of \cref{lem:FBL_extension} with
essentially the same proof.

\begin{lem}\label{lem:extension}
    Let $E$ and $F$ be Banach spaces for which $\FBFA E$ and $\FBFA F$
    are semiprime. Let $T\colon E\to F$ be a contractive operator.
    Then the
    extension operator $\bar{T}\colon \FBFA E\to \FBFA F$ satisfies
    $\hat{\iota}_F(\bar{T}f)=\hat{\iota}_E (f)\circ T^{*}$.
\end{lem}
\begin{proof}
    The composition operator
    \[
    \begin{array}{cccc}
    C_{T^{*}}\colon& \FBFA E & \longrightarrow & C(B_{F^{*}}) \\
            & f & \longmapsto & \hat{\iota }_E(f)\circ T^{*} \\
    \end{array}
    \]
    is a lattice-algebra homomorphism. Since
    $C_{T^{*}}(\eta _x)=\iota_{F}(\eta _{Tx})$, $C_{T^{*}}$ maps $\FNFA E$ to
    $\iota_{F} (\FNFA F)$ contractively, so the range of $C_{T^{*}}$ is contained
    in $\hat{\iota }_F(\FBFA F)$. By uniqueness,
    $C_{T^{*}}=\hat{\iota }_F\circ\bar{T}$.
\end{proof}

Recall from \cref{prop:FBL_extension_properties} that an operator
between Banach spaces $T\colon E\to F$ is injective, onto or has
dense range if and only if the induced operator between the free
Banach lattices $\bar{T}\colon \FBL E\to \FBL F$ is injective, onto
or has dense range, respectively. The following example shows that
for free Banach \falg s the situation is different.

\begin{example}\label{ex:extension_of_bijection}
    Let $T\colon \el 1^{2}\to \el \infty ^2$ be the formal identity.
    Identify both $B_{\el 1^2}$ and $B_{\el \infty ^2}$ as subsets of
    the Euclidean plane $\R^2$. Identify also $\FBFA{\el 1^2}$ and
    $\FBFA{\el \infty ^2}$ as sublattice-algebras of $C(B_{\el \infty
    ^2})$ and $C(B_{\el 1^2})$ through $\hat{\iota }_{\el 1^2}$
    and $\hat{\iota }_{\el \infty ^2}$, respectively
    (\cref{cor:finite_representable} ensures that this is possible). The induced operator $\bar{T}\colon
    \FBFA{\el 1^2}\to \FBFA{\el \infty ^2}$ is, according to
    \cref{lem:extension}, nothing more than the restriction
    map $\bar{T}f=f|_{B_{\el 1^2}}$. By \cref{lem:bump},
    there exists a non-zero $f \in \FBFA{\el 1^2}$ such that $\supp f
    \subseteq [1/2,1]^2$. Then $\bar{T}f=0$, so that $\bar{T}$ is not
    injective, even though $T$ is.
\end{example}

Nevertheless, we have the following positive results.

\begin{prop}\label{prop:extension_properties}
    Let $E$ and $F$ be Banach spaces and let $T\colon E\to F$ be a contractive operator.
    \begin{enumerate}
        \item $\bar{T}$ has dense range if and only if $T$ has dense
            range.
        \item If $T$ is a quotient map, then so is $\bar{T}$.
        \item If $\bar{T}$ is a bijection, then so is $T$.
        \item With the additional assumption that $\FBFA E$ and $\FBFA
            F$ are semiprime, $\bar{T}$ is injective if and only if
            $T^{*}(B_{F^{*}})$ is $w^{*}$-dense in $B_{E^{*}}$.
    \end{enumerate}
\end{prop}
\begin{proof}
    The proofs of $(i)$ and $(ii)$ are very similar to those
    in
    \cite[Proposition 2.3]{oikhberg_etal2022}. It is the proofs of
    $(iii)$ and $(iv)$ that require a different approach.
    \begin{enumerate}
        \item  Suppose first that $T$ has
            dense range. For every $y \in F$ and $\varepsilon >0$
            there exists $x \in E$ such that $\|Tx-y\|<\varepsilon $.
            This is the same as $\|\bar{T}\eta _x-\eta 
            _y\|<\varepsilon $, which implies that $\overline{\range
            \bar{T}}$ contains $\{\, \eta _y : y\in F \, \} $. Since
            $\bar{T}$ is a lattice-algebra homomorphism,
            $\overline{\range \bar{T}}$ is a closed
            sublattice-algebra. Hence it must be $\FBFA F$.

            To show the converse implication, suppose $\range T$ is
            not dense, and let $y^{*}\in B_{F^{*}}$ be a non-zero
            functional that vanishes on it. Then the map
        $\widehat{y^{*}} \in \FBFA{F}^{*}$ vanishes on $\bar{T}\eta
            _x$ for every $x \in E$. Since $\widehat{y^{*}}$ is a
            lattice-algebra homomorphism, it vanishes on $\range
            \bar{T}$. Hence the range of $\bar{T}$ is not dense.

        \item Suppose $T$ is a quotient map. Let $Z=\FBFA E/\ker \bar{T}$ and
            let $Q\colon \FBFA E\to Z$ be the quotient map. Recall
            that $Z$ is a Banach \falg\ and that $Q$ is a
            lattice-algebra homomorphism. By standard factorization
            arguments there exists an injective and contractive
            operator $S\colon Z\to \FBFA F$ such that $\bar{T}=SQ$.
            Let us show that $S$ is also a lattice-algebra
            homomorphism. Fix $z_1,z_2 \in Z$ and let $x_1,x_2 \in
            \FBFA E$ be such that $Qx_i=z_i$ for $i=1,2$. Since both
            $\bar{T}$ and $Q$ are lattice-algebra homomorphisms:
            \[
            S|z_1|=S|Qx_1|=SQ|x_1|=\bar{T}|x_1|=|\bar{T}x_1|=|Sz_1|,
            \]
            and also
            \begin{align*}
                S(z_1z_2)=S(Qx_1Qx_2)&=SQ(x_1x_2)\\&=\bar{T}(x_1x_2)=\bar{T}(x_1)\bar{T}(x_2)=Sz_1
            Sz_2.
            \end{align*}
            Since $\ker T\subseteq \ker Q\eta _E$, $Q \eta _E$ is
            contractive and $T$ is a quotient map, there exists a
            contractive operator $R\colon F\to Z$ such that $Q\eta _E=RT$. Let
            $\hat{R}\colon \FBFA F\to Z$ be its canonical extension.
            Fix $y \in F$ and let $x \in E$ be such that $Tx=y$. Then
            \[
            S \hat{R}\eta_F(y)=SRy=SRTx=SQ\eta _E x=\bar{T}\eta
            _Ex=\eta _Fy.
            \]
            Since $S \hat{R}$ is a lattice-algebra homomorphism, it
            must be the identity. This shows that $S$ is surjective,
            and therefore an isometry with inverse $\hat{R}$. In
            particular, $\bar{T}=SQ$ is a quotient map.
        \item Suppose $\bar{T}$ is bijective. It is obvious that $T$
            is injective and, by $(i)$, $T$ also has dense range.
            Let $y \in F$ and let $(x_n)$ be a sequence in $E$ such
            that $Tx_n\to y$. Then $\bar{T}(\eta _E(x_n))\to
            \eta _F(y)$ in $\FBFA F$. Since $\bar{T}$ is surjective,
            there exists $f \in \FBFA E$ for which $\bar{T}f=\eta
            _F(y)$. From applying $\bar{T}^{-1}$ to $\bar{T}f=\lim
            \bar{T}(\eta _E(x_n))$ it follows that $\eta _E(x_n)\to f$
            in $\FBFA E$. Hence $(\eta _E(x_n))$ is a Cauchy sequence
            and, since $\eta _E$ is an isometric embedding, so is
            $(x_n)$. If $x$ denotes the limit of $(x_n)$ in $E$, then
            $Tx=y$.

        \item Suppose
            $T^{*}(B_{F^{*}})$ is not dense in $B_{E^{*}}$. Let
            $U\subseteq B_{E^{*}}$ be a weak$^*$ open set disjoint from
            $T^{*}(B_{F^{*}})$. By \cref{lem:bump}, there exists a non-zero
            $f \in \FNFA E$ such that $\supp \iota _E(f) \subseteq U$.
            Using \cref{lem:extension} one gets that, for every
            $y^{*}\in B_{F^{*}}$, $\hat{\iota
            }_F(\bar{T}f)(y^{*})=\hat{\iota }_E(f)(T^{*}y^{*})=0$.
            Hence $\bar{T}f=0$ and $\bar{T}$ is not injective.

            Conversely, suppose $T^{*}(B_{F^{*}})$ is dense in
            $B_{E^{*}}$. Let $f \in \FBFA E$ be such that
            $\bar{T}f=0$. Using again the formula from
            \cref{lem:extension},
            \[
                0=\hat{\iota }_F(\bar{T}f)(y^{*})=\hat{\iota }_E(f)(T^{*}y^{*})
            \]
            for every $y^{*}\in B_{F^{*}}$. Thus $\hat{\iota
            }_E(f)$ vanishes on
            $T^{*}(B_{F^{*}})$ and, since it is continuous, the
            density assumption implies that it vanishes on
            $B_{E^{*}}$. Hence $f=0$.\qedhere
    \end{enumerate}
\end{proof}

\Cref{ex:extension_of_bijection} shows that the converse implication
in $(iii)$ above is not true. However, two properties remain to be
clarified, namely:
\begin{question}
    \begin{enumerate}
        \item If $\bar{T}$ is a quotient map, is $T$ also a quotient
            map?
        \item Does $\bar{T}$ surjective imply that $T$ is surjective?
            What about the converse?
    \end{enumerate}
\end{question}

Isometries and isometric embeddings are treated in the next section.

\section{Isometries and isometric embeddings}
\label{sec:isometries}

It is direct from the universal property that two isometric Banach
spaces have isometric free Banach \falg s. Do isometric free Banach
\falg s correspond to isometric Banach spaces? The goal of this
section is to show that they do when the Banach spaces have smooth
dual. The proof follows the same steps as
\cref{thm:FBL_isometries}, but one has to proceed with more care:
every isometry between the free Banach \falg s induces a
$w^{*}$-continuous function between the unit balls of the duals that
need not be positively homogeneous. Therefore this function cannot
be extended to the whole dual, and the results in
\cite{ilisevic_turnesk2020} cannot be applied directly. Instead, a
refinement of the result in that paper is needed.

Recall that a Banach space $X$ is said to be \emph{smooth} if for
every $x \in X$ there exists a unique $f_x \in S_{X^{*}}$ such that
$f_x(x)=\|x\|$. This functional $f_x$ is called the \emph{support
functional} at $x$ and is given by
\[
f_x(y)=\lim_{t\to 0}\frac{\|x+ty\|-\|x\|}{t}.
\]
The quantity $[x,y]=\|y\|f_y(x)$ defines a semi-inner product with
$\|x\|=[x,x]^{1/2}$. We say that $x$ \emph{is orthogonal to} $y$, and
denote it by $x\perp y$, if $[y,x]=0$. Note that orthogonality is right
additive. Moreover, Birkhoff--James orthogonality coincides with this
notion of orthogonality; that is,
\[
x\perp y\quad\text{if and only if}\quad\|x+\lambda y\|\ge \|x\|\text{ for
all }\lambda \in \R.
\]

\begin{thm}\label{thm:wigner}
    Let $X$ and $Y$ be smooth normed spaces and
    let $f\colon B_X\to B_Y$ be a surjective mapping satisfying
    \[
        |[f(x),f(y)]|=|[x,y]|\quad\text{for all }x,y \in B_X.
    \]
    Then there exist a linear surjective isometry $U\colon X\to Y$
    and a map $\sigma \colon B_X\to \{-1,1\}$ such that $f=\sigma
    U|_{B_X}$.
\end{thm}
\begin{proof}
    The first step is to show that, for every $x \in B_X$ and $\lambda
    \in \R$ such that $\lambda x \in B_X$, there exists a scalar
    $\gamma =\gamma (\lambda ,x)$ with $|\gamma |=|\lambda |$ such
    that $f(\lambda x)=\gamma f(x)$. This $\gamma $ is defined by the
    identity
    \[
    \min_{\xi \in \R}\|f(\lambda x)-\xi f(x)\|=\|f(\lambda x)-\gamma
    f(x)\|.
    \]
    Note in particular that $\|f(\lambda x)-\gamma f(x)\|\le
    \|f(\lambda x)\|\le 1$. By definition
    \[
    \|f(\lambda x)-\gamma f(x)+\mu f(x)\|\ge \|f(\lambda x)-\gamma
    f(x)\|
    \]
    holds for every $\mu \in \R$, so by Birkhoff--James orthogonality, $f(\lambda x)-\gamma f(x) \perp f(x)$. Since $f$ is
    surjective and $f(\lambda x)-\gamma f(x) \in B_Y$, there exists $z
    \in B_X$ such that $f(z)=f(\lambda
    x)-\gamma f(x)$. From $f(z)\perp f(x)$ we get $z\perp x$ and then
    $z\perp \lambda x$ and $f(z)\perp f(\lambda x)$. Using that
    the semi-inner product is right additive we conclude that
    \[
    f(z)\perp (f(\lambda x)-\gamma f(x)).
    \]
    But the right hand side is precisely $f(z)$.
    Thus $f(z)=0$ and $f(\lambda x)=\gamma f(x)$. Moreover,
    \[
    |\lambda | \|x\|=\|\lambda x\|= \|f(\lambda x)\|=\|\gamma
    f(x)\|=|\gamma | \|f(x)\|= |\gamma | \|x\|
    \]
    implies $|\lambda |= |\gamma |$.

    Next, let $x,y \in B_X$ be linearly independent and such that $x+y
    \in B_X$. We want to show
    that $f(x+y)=\alpha f(x)+\beta f(y)$, where $\alpha =\alpha
    (x,y)$, $\beta =\beta (x,y)$ and $|\alpha |=1=|\beta |$. As
    before, define $\alpha ,\beta \in \R$ by the relation
    \[
    \min_{\xi ,\eta \in \R} \|f(x+y)-\xi f(x)-\eta
    f(y)\|=\|f(x+y)-\alpha f(x)-\beta f(y)\|
    \]
    and note that $f(x+y)-\alpha f(x)-\beta f(y) \in B_Y$.
    Then
    \[
    f(x+y)-\alpha f(x)-\beta f(y)\perp f(x)\quad\text{and}\quad
    f(x+y)-\alpha f(x)-\beta f(y)\perp f(y).
    \]
    Take $z \in B_X$ such that $f(z)=f(x+y)-\alpha f(x)-\beta f(y)$.
    Then $f(z)\perp f(x)$ implies $z\perp x$, $f(z)\perp f(y)$ implies
    $z\perp y$ and the right additivity of orthogonality implies
    $z\perp x+y$ and $f(z)\perp f(x+y)$. Hence $f(z)\perp f(z)$ and
    $f(z)=0$. Let us show that $|\alpha |= 1$. Let $\lambda _0\in
    \R$ be defined by
    \[
    \|x+\lambda _0y\|=\min_{\lambda \in \R} \|x+\lambda y\|.
    \]
    Once again, $x+\lambda _0y \in B_X$.
    Then $x+\lambda _0y \perp y$ and $x+\lambda _0y\not\perp x$, for
    otherwise $x+\lambda _0y=0$, contradicting the fact that $x$ and
    $y$ are linearly independent. Denote $w=x+\lambda _0y$. Since
    $w\perp y$ we also have $f(w)\perp f(y)$. Then
    $[f(x+y),f(w)]=\alpha [f(x),f(w)]$ and
    \[
        |\alpha | |[x,w]|= |\alpha |
        |[f(x),f(w)]|=|[f(x+y),f(w)]|= |[x+y,w]|=|[x,w]|.
    \]
    Since $[x,w]\neq 0$, it follows $|\alpha |= 1$. Similarly,
    $|\beta |= 1$.

    Let us consider the projective spaces $\mathbb{P}X$ and $\mathbb{P}Y$ obtained as respective quotients of $X$ and $Y$ by identifying proportional vectors. Note that the map
    \[
    \begin{array}{cccc}
    \tilde f\colon& \mathbb{P}X & \longrightarrow & \mathbb{P}Y \\
            & \left< x \right> & \longmapsto & \left< f(x/\|x\|) \right> \\
    \end{array}
    \]
    is well-defined and surjective. We want to check that it satisfies
    the assumptions of the Fundamental Theorem of Projective Geometry
    (see \cite[Theorem 2.1]{ilisevic_turnesk2020}) whenever $\dim X\ge 3$.

    First we need to check its image is not contained in a projective
    line. Let $x \in S_X$ be a unit vector. Since $\dim X\ge 3$ we can
    choose unit vectors $y \in \ker f_x$ and $z \in \ker f_x \cap \ker
    f_y$. From $x\perp y$, $x\perp z$ and $y\perp z$ it follows that
    $x,y,z$ are linearly independent. Then $f(x),f(y),f(z)$ are
    unit vectors such that $f(x)\perp f(y)$, $f(x)\perp f(z)$ and
    $f(y)\perp f(z)$. As before we conclude that they are linearly
    independent. So the image of $f$ is not contained in a
    two-dimensional subspace, thus the image of $\tilde f$ is not
    contained in a projective line.

    Second we need to check that $\tilde f$ preserves projective
    lines. Suppose $c=\lambda a+\mu b$ with $a,b,c \neq 0$. Without
    loss of generality we may assume $a,b \in S_X$. Set
    $M=\max\{|\lambda |, |\mu |, \|c\|\}$. Then
    \begin{align*}
        \tilde f(\left< c \right>)&=\left< f(\frac{c}{\|c\|}) \right>\\
                                  &=\left<
                                  f(\frac{M}{\|c\|}\frac{c}{M})
                                  \right>\\
                                  &=\left< f(\frac{\lambda
                                  }{M}a+\frac{\mu }{M}b)\right>
    \end{align*}
    where in the third equality we have used that both $c/M$ and
    $c/\|c\|$ belong to $B_X$. After the previous computation, using
    the properties of $f$ together with the fact that $a,b, \lambda
    /Ma, \mu /Mb \in B_X$, it follows that $\tilde f(\left< c
    \right>)$ is contained in the projective line that passes through
    $\left< a \right>$ and $\left< b \right>$.

    By the Fundamental Theorem of Projective Geometry there exists a
    bijective linear map $A\colon X\to Y$ such that
    \[
    \tilde f(\left< x \right>)=\left< Ax \right>\quad\text{for all }x \in
    X.
    \]
    This means in particular that for every $x \in S_X$ there exists a
    nonzero $\lambda (x)\in \R$ such that $f(x)=\lambda (x)Ax$. More
    generally, for nonzero $x \in B_X$, since $f(x/\|x\|)=\gamma
    f(x)$, one also has $f(x)=\lambda (x)Ax$ for some $\lambda (x)\neq
    0$. Note that, if $x \in B_X$ and $\mu \in \R$ are nonzero, and
    $\mu x \in B_X$, then
    \[
    \gamma \lambda (x)Ax=\gamma f(x)=f(\mu x)=\lambda (\mu x) \mu Ax
    \]
    and therefore
    \[
    |\lambda (\mu x)|=|\gamma/\mu | |\lambda (x)|=|\lambda (x)|.
    \]
    Let $x,y \in B_X$ be linearly independent and such that $x+y
    \in B_X$. Then
    \[
    \lambda (x+y)A(x+y)=f(x+y)=\alpha f(x)+\beta f(y)=\alpha \lambda
    (x)Ax+\beta \lambda (y) Ay
    \]
    and since $Ax$ and $Ay$ are linearly independent, it follows that
    \[
    |\lambda (x)|=|\alpha \lambda (x)|=|\lambda (x+y)|=|\beta \lambda
    (y)|=|\lambda (y)|.
    \]
    In conclusion, $|\lambda (x)|$ is constantly $|\lambda |$. By
    letting $\sigma (x)=\lambda (x)/|\lambda |$ and $U=|\lambda |A$
    one gets the desired result for $\dim X\ge 3$.

    For $\dim X=2$, let $x_0,y_0 \in X$ be such that $x_0 \perp y_0$
    and $\|x_0\|=\|y_0\|=1/2$. For every
    $\mu $ with $\|x_0+\mu y_0\|\le 1$ (in particular, for every $\mu
    \in [-1,1]$), there exist $\omega _1, \omega _2 \in \R$ such that
    \[
    f(x_0+\mu y_0)=\omega _1 f(x_0)+\omega _2f(y_0)
    \]
    with $|\omega _1|=1$ and $|\omega _2|=|\mu| $.

    Define $\omega (\mu)=\omega _1$ and $h(\mu )=\omega _2/\omega _1$.
    Once we show that $h(\mu )=\mu h(1)$ for all $\mu \in [-1,1]$, most of the work will be done.
    For every $\lambda ,\mu \in [-1,1]$ with $\|x_0+(\lambda +\mu
    )y_0\|\le 1$ we have
    \[
    f(x_0+(\lambda +\mu )y_0)=\omega (\lambda +\mu )(f(x_0)+h(\lambda +\mu
    )f(y_0))
    \]
    but also
    \begin{align*}
        f(x_0+(\lambda +\mu )y_0)&=f( (x_0+\lambda y_0) + \mu y_0)\\
                     &=\omega _4'f(x_0+\lambda y_0)+\omega _3 f(y_0)\\
                     &=\omega _4f(x_0)+\omega _4h(\lambda
                     )f(y_0)+\omega _3f(y_0)
    \end{align*}
    for certain real numbers $\omega _3,\omega _4',\omega _4$ with
    $|\omega _3|=|\mu |$, $|\omega _4'|=1$ and $|\omega _4|=1$. Since
    $f(x_0)$ and $f(y_0)$ are linearly independent, $\omega _4=\omega
    (\lambda +\mu )$ and
    \[
    \omega (\lambda +\mu )h(\lambda )+\omega _3=\omega (\lambda +\mu ) h(\lambda +\mu )
    \]
    which implies
    \begin{equation}\label{eq:h1}
    h(\lambda +\mu )=h(\lambda )+ \frac{\omega _3}{\omega}.
    \end{equation}
    Taking absolute values, we have
    \[
    |\lambda +\mu |=\left| h(\lambda )+\frac{\omega _3}{\omega } \right|.
    \]
    Since $|\omega _3/\omega |=|\mu |$, if $\mu \neq
    0$, then we can divide by $|\omega _3/\omega |$ to get
    \[
    \left| \frac{\lambda }{\mu }+1 \right| =\left| h(\lambda
    )\frac{\omega }{\omega _3} +1\right|.
    \]
    Using that $|h(\lambda )|=|\lambda |$ it easily follows
    \begin{equation}\label{eq:h2}
    \frac{\lambda }{\mu }=h(\lambda )\frac{\omega }{\omega _3}.
    \end{equation}
    Hence, if also $\lambda \neq 0$, then we can put together
    \eqref{eq:h1} and \eqref{eq:h2} to get
    \[
    h(\lambda +\mu )=h(\lambda ) \left( 1+\frac{\mu }{\lambda }
    \right).
    \]

    Fix $0\le \eta<1$. In this case, $\|x_0+\eta y_0\|\le 1$. We may set $\lambda =1$ and $\mu =\eta -1$ in
    previous equation to get $h(\eta )=\eta h(1)$. Similarly, $\lambda
    =-1$ and $\mu =1-\eta $ yield $h(-\eta )=\eta h(-1)$. Finally,
    setting $\lambda =1/2$ and $\mu =-1$ yields $h(-1/2)=h(1/2)$,
    which in turn implies $h(-1)=h(1)$. Hence $h(\mu)=\mu  h(1)$
    holds for all $\mu  \in [-1,1]$; that is,
    \[
    f(x_0+\mu  y_0)=\omega (\mu  )(f(x_0)+\mu h(1) f(y_0)),
    \]
    where $|\omega (\mu )|=1=|h(1)|$.

    By the same argument, for every $\lambda \in [-1,1]$ there exist signs $\omega' (\lambda )$ and $g(1)$ such that
    \[
    f(\lambda x_0+y_0)=\omega '(\lambda )(\lambda g(1)f(x_0)+f(y_0)).
    \]
    Since
    \[
    \omega '(1)(g(1)f(x_0)+f(y_0))=f(x_0+y_0)=\omega
    (1)(f(x_0)+h(1)f(y_0))
    \]
    we conclude that $h(1)g(1)=1$, that is, $h(1)=g(1)$.

    Define a linear map $U\colon
    X\to Y$ by $Ux_0=f(x_0)$ and $Uy_0=h(1)f(y_0)$. We are going to
    show that $f$ coincides with $U|_{B_X}$ up to a sign. Every
    element of the unit ball can be written uniquely as $\lambda
    x_0+\mu y_0$ for certain $\lambda ,\mu \in \R$. Suppose first that
    $|\lambda |\ge |\mu |$. If $\lambda =0$, then $f(0)=0=U(0)$, so
    assume also $\lambda \neq 0$. Then $x_0+\mu /\lambda y_0 \in B_X$,
    and therefore there exists a $\gamma \in \R$ with $|\gamma
    |=|\lambda |$ such that
    \begin{align*}
        f(\lambda x_0+\mu y_0)&=\gamma f(x_0+\mu /\lambda y_0)\\
                              &=\gamma \omega (\mu /\lambda
                              )(f(x_0)+\mu /\lambda h(1)f(y_0))\\
                              &=\gamma \omega U(x_0+\mu /\lambda
                              y_0)\\
                              &=\frac{\gamma }{\lambda }\omega
                              U(\lambda x_0+\mu y_0).
    \end{align*}
    Hence $f$ coincides with $U$ at $\lambda x_0+\mu y_0$ up to a sign
    $\omega \gamma /\lambda $. Similarly, if $|\mu |>|\lambda |$, then
    there exists a $\delta \in \R$ with $|\delta |=|\mu|$ such
    that
    \begin{align*}
        f(\lambda x_0+\mu y_0)&=\delta f(\lambda /\mu x_0+y_0)\\
                              &=\delta \omega '(\lambda /\mu )(\lambda
                              /\mu g(1)f(x_0)+f(y_0))\\
                              &=\delta \frac{\omega '}{\mu
                              }g(1)(\lambda f(x_0)+\mu h(1)f(y_0))\\
                              &=\delta \frac{\omega '}{\mu
                              }g(1)U(\lambda x_0+\mu y_0),
    \end{align*}
    and again $f$ coincides with $U$ at $\lambda x_0+\mu y_0$ up to a
    sign $\delta \omega '/\mu g(1)$.

    It only remains to check that $U$ preserves the norm. An element
    of $X$ can be written uniquely as $\lambda x_0+\mu y_0$.
    Suppose first that $|\lambda |>|\mu |$, so as to have $\|x_0+\mu
    /\lambda y_0\|\le 1$. In this case,
    \begin{align*}
        \|\lambda x_0+\mu y_0\|&=|\lambda | \|x_0+\mu /\lambda y_0\|\\
        &=|\lambda | \|f(x_0+\mu /\lambda
        y_0)\|\\& =|\lambda | \|U(x_0+\mu /\lambda y_0)\|\\
                               &=\|U(\lambda x_0+\mu y_0)\|.
    \end{align*}
    The case $|\mu |\ge |\lambda |$ is completely analogous. This
    completes the proof when $\dim X=2$.

    When $\dim X=1$, fix a unit vector $x_0\in X$, and define a linear
    map $U\colon X\to Y$ by $U(\lambda x_0)=\lambda f(x_0)$. This map
    clearly preserves the norm. For every $\lambda \in [-1,1]$,
    $\lambda \neq 0$, there
    exists $\gamma \in \R$ with $|\gamma |=|\lambda |$ such that
    \[
    f(\lambda x_0)=\gamma f(x_0)=\gamma /\lambda U(\lambda x_0).
    \]
    Hence, $f$ coincides with $U$ at $\lambda x_0$ up to a sign $\gamma
    /\lambda $.
\end{proof}

\begin{rem}
    The previous theorem is still true if $f\colon \delta B_X\to
    \delta B_Y$, for some $\delta >0$, and $f(\delta B_X)=\delta B_Y$. Indeed, let $X_\delta$ (resp.\
    $Y_\delta$) be the
    space $X$ (resp.\ $Y$) with the norm scaled by $1/\delta$.
    These are still smooth normed spaces satisfying $B_{X_\delta
    }=\delta B_X$ and $B_{Y_\delta }=\delta B_Y$. Moreover, one can
    easily check that the semi-inner products are related by
    $[x,y]_{X_\delta }=[x,y]_X/\delta ^2$, and the same is true for $Y$ and
    $Y_\delta $. Hence $f\colon B_{X_\delta }\to B_{Y_\delta }$ is a
    surjective mapping that preserves the semi-inner product. By
    the previous theorem, there exist a linear surjective isometry
    $U\colon X_\delta \to Y_\delta $ and a map $\sigma \colon
    B_{X_\delta }\to \{-1,1\}$ such that $f=\sigma U|_{B_{X_\delta
    }}$. It is then clear that $U\colon X\to Y$ is still a linear
    surjective isometry, and that $f=\sigma U|_{\delta B_X}$.
\end{rem}

With this we can now state and prove our main theorem in this section. The proof
follows closely that of \cref{thm:FBL_isometries}.

\begin{thm}\label{thm:isometries}
    Let $E$ and $F$ be Banach spaces with smooth dual. An operator
    $T\colon \FBFA E\to \FBFA F$ is a surjective lattice-algebra
    isometry if and only if $T=\bar{V}$, for some surjective isometry
    $V\colon E\to F$. In particular, $E$ and $F$ are isometric if and
    only if $\FBFA E$ and $\FBFA F$ are lattice-algebra isometric.
\end{thm}
\begin{proof}
    First we are going to construct a map $\phi _T\colon B_{F^{*}}\to
    B_{E^{*}}$. Let $y^{*}\in B_{F^{*}}$. Extend $y^{*}\colon F\to \R$
    to a unique contractive lattice-algebra homomorphism
    $\widehat{y^{*}}\colon \FBFA F\to \R$. Define $\phi_T
    (y^{*})=T^{*}(\widehat{y^{*}})\circ \eta _E$ (that is,
    $\widehat{\phi _T(y^{*})}=T^{*}(\widehat{y^{*}})$).

    We will use the
    lattice-algebraic homomorphisms $\hat{\iota}_E\colon \FBFA E\to
    C(B_{E^{*}})$ and $\hat{\iota }_F\colon \FBFA F\to
    C(B_{F^{*}})$ to express $\phi _T$ in
    another way. Note that, if $g \in \FBFA
    F$, then $\widehat{y^{*}}(g)=\hat{\iota }_F(g)(y^{*})$. Indeed,
    one can easily check that $g\mapsto \hat{\iota }_F(g)(y^{*})$
    is a contractive lattice-algebra homomorphism extending
    $y^{*}$. If $f \in \FBFA E$ and $y^{*}\in B_{F^{*}}$, then
    \[
        (\hat{\iota }_E(f)\circ \phi_T )(y^{*})=\hat{\iota
        }_E(f)(\phi _T(y^{*}))=\widehat{\phi_T
        (y^{*})}(f)=T^{*}(\widehat{y^{*}})(f)=\hat{\iota }_F(Tf)(y^{*}).
    \]
    Thus $\hat\iota _F(Tf)=\hat\iota _E(f)\circ \phi _T$. Plugging in $f=\eta _E(x)$ yields $\phi
    _T(y^{*})(x)=\hat{\iota }_F(T\eta _E(x))(y^{*})$ for all $y^{*}\in B_{F^{*}}$.
    In particular, if $(y_\alpha ^{*})$ is a net that weak$^{*}$
    converges to $y^{*}$ in $B_{F^{*}}$, then
    \[
    \phi _T(y_\alpha ^{*})(x)=\hat{\iota }_F(T\eta _E(x))(y_\alpha
    ^{*})\to \hat{\iota }_F(T\eta
    _E(x))(y^{*})=\phi _T(y^{*})(x).
    \]
    It follows that $\phi _T$ is weak$^{*}$ to weak$^{*}$ continuous.
    It is also direct to check that $\phi _T$ is invertible with $\phi
    _T^{-1}=\phi _{T^{-1}}$.

    Suppose that $B_{E^{*}}$ and $B_{F^{*}}$ are equipped with the
    corresponding semi-inner products. We are going to show that $\phi
    _T$ preserves the semi-inner product when restricted to
    $\frac{1}{2}B_{F^{*}}$. According to \cite[Lemma
    10.19]{oikhberg_etal2022}, showing that $\phi _T$ preserves the
    semi-inner product is equivalent to showing that
    \[
    \max_{\pm} \|\phi_T (x^{*})\pm \phi_T (y^{*})\|=\max_{\pm } \|x^{*}\pm
    y^{*}\|
    \]
    holds for every $x^{*}, y^{*} \in \frac{1}{2}B_{F^{*}}$. Consider the
    operator
    \[
    \begin{array}{cccc}
    S\colon& F & \longrightarrow & \el 1^{2} \\
            & y & \longmapsto & (x^{*}(y),y^{*}(y)) \\
    \end{array},
    \]
    where $\el 1^2$ is equipped with the pointwise product.
    It is direct to check that $\|S\|=\max_{\pm }\|x^{*}\pm y^{*}\|\le 1$.
    One can also check that its extension $\hat{S}\colon \FBFA F\to
    \el 1^2$ is $\hat{S}(f)=(\hat{\iota
    }_F(f)(x^{*}),\hat{\iota }_F(f)(y^{*}))$. Then $\hat{S}T\colon
    \FBFA E\to \el 1^2$ is a lattice-algebra homomorphism with
    $\|\hat{S}T\|=\|\hat{S}\|$. Moreover, the universal property
    of the free Banach \falg\ implies $\|\hat{S}T\eta
    _E\|=\|\hat{S}T\|$. Unfolding the definitions:
    \[
    \hat{S}T(\eta _E(x))=(\hat{\iota }_F((T\eta
    _E)(x))(x^{*}),\hat{\iota }_F((T\eta
    _E)(x))(y^{*}))=(\phi _T(x^{*})(x),\phi _T(y^{*})(x)).
    \]
    By direct computation, $\|\hat{S}T\eta _E\|=\max_{\pm }\|\phi
    _T(x^{*})\pm \phi _T(y^{*})\|$, and this coincides with
    $\|S\|=\max_{\pm}\|x^{*}\pm y^{*}\|$.

    For the same reasons, $\phi _{T^{-1}}$ preserves the semi-inner
    product when restricted to $\frac{1}{2}B_{E^{*}}$. This implies
    that
    $\phi _T|_{\frac{1}{2}B_{F^{*}}}\colon \frac{1}{2}B_{F^{*}}\to
    \frac{1}{2}B_{E^{*}}$ is surjective.
    Hence we can apply \cref{thm:wigner} and the remark after it to
    this map: there exists a linear surjective isometry
    $U\colon F^{*}\to E^{*}$ and a map $\sigma \colon \frac{1}{2}B_{F^{*}}\to
    \{-1,1\}$ such that $\phi _T=\sigma U|_{\frac{1}{2}B_{F^{*}}}$. We claim that
    $\sigma $ is continuous on $\frac{1}{2}B_{F^{*}}\setminus \{0\}$ for the norm
    topology. Suppose it was not. Then there would exist a convergent
    sequence $y_k^{*}\to y^{*}$ in $\frac{1}{2}B_{F^{*}}\setminus \{0\}$ such
    that $\sigma (y_k^{*})$ does not converge to $\sigma (y^{*})$.
    To fix ideas, suppose $\sigma (y^{*})=-1$ (the other case is
    identical). Then there exists a subsequence of $y_k^{*}$, which we
    denote the same way, for which $\sigma (y_k^{*})=1$. Note that
    $Uy_k^{*}\to Uy^{*}$. At the same time, since norm convergence
    implies weak$^*$ convergence, $Uy_k^{*}=\phi _T(y_k^{*})\to \phi
    _T(y^{*})=-Uy^{*}$. Hence $Uy^{*}=0$, a contradiction with the
    fact that $y^{*}\neq 0$.

    When $\dim F>1$, $\frac{1}{2}B_{F^{*}}\setminus \{0\}$ is connected, and
    therefore $\sigma $, being continuous, must be constant. Absorbing
    $\sigma $ in $U$, we have $\phi
    _T=U|_{\frac{1}{2}B_{F^{*}}}$, and therefore
    that $U$ is weak$^*$ continuous on bounded sets. This implies that
    $U$ is an adjoint operator; its pre-adjoint $V$ is an isometry
    such that $T=\bar{V}$ (the details are the same as in the end of the
    proof of \cite[Theorem 10.18]{oikhberg_etal2022}). When $\dim
    F=1$, but $\dim E>1$, we can use the same argument with $\phi
    _{T^{-1}}$.

    It only remains to consider the case $\dim E=\dim F=1$. Recall
    from \cref{thm:finite} that, in this case, $\FBFA E$ is isomorphic
    to
    $C([0,1]\times \{-1,1\})$ with the supremum norm, pointwise
    lattice structure, and product
    \[
        (f \star g)(r,u )=rf(r,u )g(r,u
        )
    \]
    for every $f,g \in C([0,1]\times \{-1,1\})$ and $(r,u
    )\in [0,1]\times \{-1,1\}$. The canonical map $\eta \colon \R\to
    C([0,1]\times \{-1,1\})$ is determined by $\eta (1)(r,u
    )=u $. Suppose
    \[
    T\colon C([0,1]\times \{-1,1\})\to C([0,1]\times \{-1,1\})
    \]
    is a surjective lattice-algebra isometry. Being a lattice
    isometry, it must be of the form $Tf=f\circ \phi $, where $\phi
    \colon [0,1]\times \{-1,1\}\to [0,1]\times \{-1,1\}$ is a
    homeomorphism (\cite[Theorem 4.25]{abramovich_aliprantis2002}).
    Put $\phi (r,u )=(\phi _1(r,u ),\phi
    _2(r,u ))$.
    Since $T$ is also an algebra homomorphism,
    \[
    T(fg)(r,u )=\phi _1(r,u )f(\phi (r,u
    ))g(\phi (r,u ))
    \]
    must be equal to
    \[
        (TfTg)(r,u )=r f(\phi (r,u ))g(\phi
        (r,u ))
    \]
    for every $f,g \in C([0,1]\times \{-1,1\})$. Thus $\phi
    _1(r,u )=r$. Being $\phi _2$ continuous, it
    must be constant on connected components. That is, there exists
    $\sigma \in \{-1,1\}$ such that $\phi _2(r,u )=\sigma
    u $. The isometry
    \[
    \begin{array}{cccc}
    V\colon& \R & \longrightarrow & \R \\
            & x & \longmapsto & \sigma x \\
    \end{array}
    \]
    is such that
    \[
        T(\eta _x)(r,u )=\eta _x(\phi
        (r,u ))=x \phi _2(r,u )=x \sigma u
        =\eta _{Vx}(r,u ).
    \]
    Hence $T=\bar{V}$, as wanted.
\end{proof}

What can be said about isometric embeddings? In the case of free
Banach lattices, the problem of whether an embedding between two
Banach spaces becomes an embedding between the free Banach lattices is
reduced to a certain Banach space question involving extensions of
operators to $L_1$ (see \cite[Theorem 3.7]{oikhberg_etal2022}). In one
direction, something similar can be said about free Banach \falg s.
But now, instead of considering operators with range in $L_1$, we have
to consider operators with range in arbitrary Banach \falg s.

More precisely, suppose that $i\colon E\to F$ is an isometric
embedding between Banach spaces, and that for every contractive
operator $T_E\colon E\to A_E$, where $A_E$ is a Banach \falg, there
exists a contractive operator $T_F\colon F\to A_F$, where $A_F$ is a
Banach \falg\ in which $A_E$ embeds isometrically through a
lattice-algebra homomorphism $j\colon
A_E\to A_F$, making the diagram
\begin{equation}\label{eq:extension_diagram}
\begin{tikzcd}
    E\arrow[r, "i"]\arrow[d, "T_E"']&F\arrow[d,"T_F"]\\
    A_E\arrow[r, "j"']&A_F
\end{tikzcd}
\end{equation}
commute. Then the map
\[
\bar{i}\colon \FBFA E\to \FBFA F
\]
is an isometric embedding. Indeed, fix $f \in \FBFA E$ and
$\varepsilon >0$, and let $T_E\colon E\to A_E$ be a contractive
operator, where $A_E$ is a Banach \falg, satisfying
$\|\hat{T}_Ef\| > \|f\|-\varepsilon$. Let $T_F\colon F\to A_F$ and
$j\colon A_E\to A_F$ be as in \eqref{eq:extension_diagram}. In the diagram
\[
\begin{tikzcd}
    \FBFA E \arrow[r, "\bar{i}"]\arrow[dd,
    "\hat{T}_E"', bend right]&\FBFA F\arrow[dd, "\hat{T}_F", bend
    left]\\
    E\arrow[u, "\eta _E"']\arrow[r, "i"]\arrow[d, "T_E"]&
    F\arrow[u, "\eta _F"]\arrow[d, "T_F"']\\
    A_E\arrow[r, "j"]&A_F
\end{tikzcd}
\]
we have that
\[
    \hat{T}_F\bar{i}\eta _E=\hat{T}_F\eta _F
    i=T_Fi=jT_E=j\hat{T}_E\eta _E.
\]
By the universal property of $\FBFA E$, it follows that
$\hat{T}_F\bar{i}=j\hat{T}_E$. Hence
\[
\|\bar{i}f\|\ge
\|\hat{T}_F\bar{i}f\|=\|j\hat{T}_Ef\|=\|\hat{T}_Ef\|>\|f\|-\varepsilon .
\]
Since $\varepsilon >0$ was arbitrary, $\|\overline{i}f\|\ge \|f\|$,
and therefore $\bar{i}$ is an isometric embedding.

\begin{question}
    Suppose $\bar{i}\colon \FBFA E\to \FBFA F$ is an isometric
    embedding. Is it true that for every contractive map $T_E\colon
    E\to A_E$, where $A_E$ is a Banach \falg, there exists a
    contractive operator $T_F\colon F\to A_F$, where $A_F$ is a
    Banach \falg\ in which $A_E$ embeds isometrically through $j$,
    making the diagram \eqref{eq:extension_diagram} commute?
\end{question}

\section{Free Banach \texorpdfstring{$f$-algebras}{f-algebras} with identity}
\label{sec:identity}

In this final section we discuss free objects in the category of
(Archimedean, normed, Banach) \falg s with an algebraic identity. We
always assume that the homomorphisms in these categories preserve
the identity.

Recall from \cref{sec:background_falg} that the identity of an
Archimedean \falg\ is positive and a weak order unit, and that every
Banach \falg\ with identity can be renormed to have an identity of
norm one. From now on, we shall assume that Banach \falg s with
identity have identities of norm one. In this case, every such
Banach \falg\ is lattice-algebra isometric to $C(K)$, for some
compact Hausdorff space $K$ (this is an immediate consequence of
\cref{thm:BLA}). Hence, the category of Banach \falg s
with identity is precisely that of $C(K)$-spaces, where the
morphisms are lattice-algebra homomorphisms preserving the
algebraic identity. This, in turn, is the category of AM-spaces with
unit whose morphisms are the lattice homomorphisms that preserve the
norming unit of the AM-space.

Let us now define precisely the object we are trying to construct.

\begin{defn}
    Let $E$ be a vector space. The \emph{free Archimedean \falg\
    with identity generated by $E$} is an Archimedean \falg\ $\FAFAov
    E$ with identity $1$ together with a
    linear map $\delta_E \colon E\to \FAFAov E$ such that, for every
    Archimedean
    \falg\ $A$ with identity $1_A$ and every linear map $T\colon E\to A$, there exists a
    unique lattice-algebra homomorphism $\hat{T}\colon \FAFAov E\to A$
    satisfying $\hat{T}(1)=1_A$ and $\hat{T}\delta _E=T$.
\end{defn}

To construct the free Archimedean \falg\ we used LLA expressions.
These were formal expressions obtained from finitely many formal
variables using: a 0-ary operation $0$, a unary operation $\lambda $
for every $\lambda \in \R$, and three binary operations $+$, $\vee $
and $\cdot $. If we add an additional 0-ary operation $1$, we obtain
expressions that can be evaluated at any vector lattice algebra with
identity. These we shall call \emph{$1$-LLA expressions}.

The following is a simple consequence of the representation theorem by
Henriksen and Johnson. The proof is essentially the same as that of
\cref{lem:calc_identity}.

\begin{lem}
    Let $\Phi $ be a $1$-LLA expression. If $\Phi $ vanishes on $\R$,
    then it also vanishes on every Archimedean \falg\ with identity.
\end{lem}

With this result, the construction of $\FAFAov E$ follows the same
steps as that of $\FAFAv E$ presented in \cref{prop:fafav}.

\begin{prop}
    Let $E$ be a vector space. The free Archimedean \falg\ with
    identity generated by $E$ is
    \[
        \FAFAov E=\VLA(\{\, \delta _x : x \in E \, \} \cup
        \{1\})\subseteq \R^{E^{\#}},
    \]
    where $1(\omega )=1$ for all $\omega \in E^{\#}$, together with
    the map $\delta_E \colon E\to \FAFAov E$, $\delta_E(x)=\delta _x$.
\end{prop}

At this point we could proceed as in \cref{sec:abstract_construction} to obtain an abstract
construction of the free Banach \falg\ with identity. But, as noted above,
this object is the same as the free AM-space with unit, which was
already considered in \cite[Theorem 5.4]{jardon-sanchez_etal2022}.

\begin{thm}[{\cite[Theorem 5.4]{jardon-sanchez_etal2022}}]
    Let $E$ be a Banach space. Let $\phi _E\colon E\to C(B_{E^{*}})$
    be the natural embedding. For every compact Hausdorff space $K$
    and every norm one operator $T\colon E\to C(K)$ there exists a
    unique lattice homomorphism $\hat{T}\colon C(B_{E^{*}})\to C(K)$
    satisfying $\hat{T}\circ \phi _E=T$ and
    $\hat{T}\one_{B_{E^{*}}}=\one_K$. Moreover, $\|\hat{T}\|=1$
    and $\hat{T}$ is an algebra homomorphism.
\end{thm}

The free Banach \falg\ with identity generated by $E$ is not only
contained in $C(B_{E^{*}})$, but is in fact the whole $C(B_{E^{*}})$,
and the free norm is the supremum norm. Note also that, if we want to
preserve the norm, we can only extend norm one operators, since every
extension has norm one (because it sends the identity to the
identity).


\backmatter
\bibliographystyle{amsalpha}
\bibliography{library.bib}


\end{document}